%% file: main.tex
\documentclass[11pt]{amsart}
 \pdfoutput=1

\usepackage{rotating}
\usepackage{multirow}
\usepackage{amsmath}

\usepackage{amsfonts,amssymb,amscd,bbm,booktabs,color,enumerate,float,graphicx,latexsym, multirow,mathrsfs,psfrag}

\usepackage[all]{xy}
\usepackage[margin=1in,marginparwidth=0.8in, marginparsep=0.1in]{geometry}

\usepackage[bookmarks=true, bookmarksopen=true,%
bookmarksdepth=3,bookmarksopenlevel=2,%
colorlinks=true,%
linkcolor=blue,%
citecolor=blue,%
filecolor=blue,%
menucolor=blue,%
urlcolor=blue]{hyperref}

\newcommand{\red}[1]{{\color{red}#1}}

\newtheorem{dummy}{dummy}[section]
\newtheorem{lemma}[dummy]{Lemma}
\newtheorem{theorem}[dummy]{Theorem}

\newtheorem{conjecture}[dummy]{Conjecture}
\newtheorem{corollary}[dummy]{Corollary}
\newtheorem{proposition}[dummy]{Proposition}
\theoremstyle{definition}
\newtheorem{definition}[dummy]{Definition}

\newtheorem{example}[dummy]{Example}
\newtheorem{remark}[dummy]{Remark}

\newcommand{\bC}{\mathbb{C}}

\newcommand{\bF}{\mathbf{F}}

\newcommand{\bP}{\mathbb{P}}

\newcommand{\bR}{\mathbb{R}}
\newcommand{\bZ}{\mathbb{Z}}

\newcommand{\cD}{\mathcal{D}}

\newcommand{\cI}{\mathcal{I}}

\newcommand{\cL}{\mathcal{L}}
\newcommand{\cM}{\mathscr{M}}

\newcommand{\cP}{\mathscr{P}}

\newcommand{\cS}{\mathcal{S}}
\newcommand{\cT}{\mathcal{T}}

\newcommand{\cX}{\mathscr {X}}

\newcommand{\op}{\operatorname}

\newcommand{\Z}{Z}

\newcommand{\PGL}{\mathrm{PGL}}

\newcommand{\Li}{\mathrm{Li}}

\numberwithin{equation}{subsection}
\numberwithin{figure}{subsection}

\newcommand{\leg}{S}

\newcommand{\necklace}{\Gamma^{\mathrm{neck}}}

\newcommand{\CC}{\mathbb{C}}

\newcommand{\G}{\mathrm{G}}

\newcommand{\Hom}{\mathrm{Hom}}

\renewcommand{\log}{{\op{log}}}
\newcommand{\lra}{\longrightarrow}

\usepackage{tikz}
\usetikzlibrary{matrix}
\usetikzlibrary{shapes.geometric,decorations.pathreplacing}
\usepackage{tikz-cd}

\title[The Chromatic Lagrangian]{The Chromatic Lagrangian:\\ \vskip.1in
Wavefunctions and Open Gromov-Witten Conjectures}
\dedicatory{To Steve Zelditch, in memoriam}
\author[Gus Schrader, Linhui Shen, and Eric Zaslow]{Gus Schrader,${}^*$ Linhui Shen,${}^{**}$  and Eric Zaslow${}^{*}$\\
\\
{\tiny ${}^{*}$  Department of Mathematics, Northwestern University}\\
{\tiny ${}^{**}$  Department of Mathematics, Michigan State University}\\
}

\begin{document}

\maketitle

\begin{abstract}
Inside a symplectic leaf of the cluster Poisson variety of Borel-decorated $PGL_2$
local systems on a punctured surface is an isotropic subvariety we will call
the \emph{chromatic Lagrangian}.  Local charts for the quantized cluster variety
are quantum tori defined by cubic planar graphs, and can be put in standard
form after some additional markings giving the notion of a \emph{framed seed}. 
The mutation structure is encoded as a groupoid.  The local description of the
chromatic Lagrangian defines a \emph{wavefunction} which, we
conjecture, encodes open Gromov-Witten invariants of a Lagrangian threefold in
threespace defined by the cubic graph and the other data of the framed seed.
We also find a relationship we call \emph{framing duality}:
for a family of ``canoe'' graphs,
wavefunctions for different framings encode DT invariants of symmetric quivers. 

\end{abstract}

\setcounter{tocdepth}{1}
\tableofcontents

\input{introduction.tex}

\input{Cluster.tex}

\input{groupoid.tex}

\input{cubic-graph} 

\input{foamsphasescones.tex}

\input{wavefunction.tex}

\input{analytic.tex}

\input{framing_and_quiver.tex}

\input{appendix.tex}

\input{references.tex}

\end{document}

%% file: introduction.tex
\section{Introduction}
\label{sec:intro}

This paper exploits cluster theory to compute wavefunctions for Lagrangian branes in
threespace
and to make explicit conjectures about their all-genus open Gromov-Witten invariants.  For certain branes, these numbers also relate to the cohomologies of twisted character varieties and Donaldson-Thomas invariants of quivers. 
Two structural tools in the schema are the behavior under mutation and the dependence of quantities on phases and framings.  
\emph{Let's now explain what we mean by all this.}

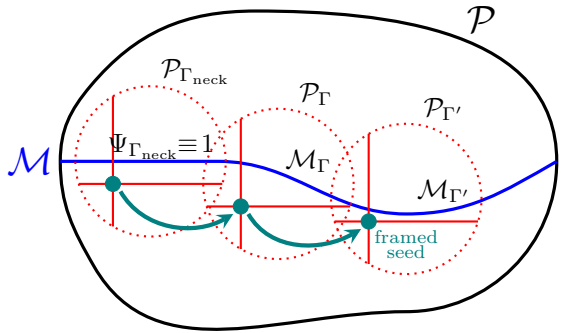
\begin{figure}[ht]
\begin{tikzpicture}[scale=1]
\pgfmathsetmacro{\a}{.3};
\pgfmathsetmacro{\bb}{({1-\a*\a})};
\pgfmathsetmacro{\b}{sqrt(\bb)};
\pgfmathsetmacro{\c}{.866};
\pgfmathsetmacro{\shiftonex}{4};
\pgfmathsetmacro{\shiftoney}{-.5};
\pgfmathsetmacro{\shifttwox}{.6};
\pgfmathsetmacro{\shifttwoy}{0};
\pgfmathsetmacro{\shiftthreex}{2.3};
\pgfmathsetmacro{\shiftthreey}{-.3};

\newcommand*{\boff}{10}; \newcommand*{\aoff}{35}; \newcommand*{\rad}{3};
\coordinate (a) at (-.15,-1.4);
\coordinate (b) at (4,-2);
\coordinate (c) at (6,0);
\coordinate (d) at (2,2);
\coordinate (e) at (-.2,1);
\coordinate (f) at (-.6,0);
\coordinate (g) at (1.5,0);
\coordinate (h) at (4,-.7);
\draw[very thick] (a) to[out=300,in=180] (b) to[out=180+180,in=270] (c) to[out=90,in=0] (d) to[out=180,in=50] (e) to[out=230, in=90] (f) to [out=270, in=120] (a);
\draw[very thick,blue] (c) to [out=210,in=0] (h) to[out=180,in=0] (g) to[out=180,in=0] (f);
\node at (5,1.9) {\Large $\mathcal P$};
\node[blue] at (-1,0) {\Large $\mathcal M$};
\node at (4.5,.7) {\small$\mathcal P_{\Gamma'}$};
\node at (2.8,.9) {\small$\mathcal P_{\Gamma}$};
\node at (1.2,1.2) {\small $\mathcal P_{\Gamma_{\mathrm{neck}}}$};
\node at (4.5,-.4) {\small $\mathcal M_{\Gamma'}$};
\node at (2.7,0) {\small $\mathcal M_{\Gamma}$};

\node[teal] at (\shiftonex,\shiftoney-.5) {\tiny framed};
\node[teal] at (\shiftonex,\shiftoney-.7) {\tiny seed};
{\draw[thick,red] (\shiftonex-.5,\shiftoney+\c) -- (\shiftonex-.5,\shiftoney-\c);
\draw[thick,red] (\shiftonex-\b,-.5-\a) -- (\shiftonex+\b,-.5-\a);
\draw[thick,dotted,red] (\shiftonex,\shiftoney) circle (1cm);
\draw[teal,fill] (\shiftonex - .5,\shiftoney - \a) circle (.1);}
{\draw[thick,red] (\shifttwox-.5,\shifttwoy+\c) -- (\shifttwox-.5,\shifttwoy-\c);
\draw[thick,red] (\shifttwox-\b,\shifttwoy-\a) -- (\shifttwox+\b,\shifttwoy-\a);
\draw[thick,dotted,red] (\shifttwox,\shifttwoy) circle (1cm);
\draw[teal,fill] (\shifttwox - .5,\shifttwoy - \a) circle (.1);}
{\draw[thick,red] (\shiftthreex-.5,\shiftthreey+\c) -- (\shiftthreex-.5,\shiftthreey-\c);
\draw[thick,red] (\shiftthreex-\b,\shiftthreey-\a) -- (\shiftthreex+\b,\shiftthreey-\a);
\draw[thick,dotted,red] (\shiftthreex,\shiftthreey) circle (1cm);
\draw[teal,fill] (\shiftthreex - .5,\shiftthreey - \a) circle (.1);}
\coordinate (ASE) at (\shiftthreex-.5+.1,\shiftthreey-\a-.1);
\coordinate (ASW) at (\shiftthreex-.5-.1,\shiftthreey-\a-.1);
\coordinate (BSW) at (\shiftonex-.5-.1,\shiftoney-\a-.1);
\coordinate (CSE) at (\shifttwox-.5+.1,\shifttwoy-\a-.1);
\draw[-stealth,ultra thick, teal] (ASE) to[out=300,in=210] (BSW);
\draw[-stealth,ultra thick, teal] (CSE) to[out=300,in=210] (ASW);
\node at (.74,.2) {\small $\Psi_{\Gamma_{\mathrm{neck}}}\!\!\equiv\! 1$};

\end{tikzpicture}
\caption{The cluster variety of decorated $PGL_2$-local systems $\cP$ and the
chromatic Lagrangian $\cM$ (in blue).  Each framed seed (teal dot) identifies the chart
$\cP_\Gamma$ with a quantum torus, in which the ideal $\cM_\Gamma$ is described by a
cyclic vector or \emph{wavefunction,} $\Psi_\Gamma$.  Arrows
in the \emph{framed seed groupoid} allow us to
determine $\Psi_{\Gamma'}$ from $\Psi_\Gamma$.
Any seed connected to the necklace graph $\Gamma_{\rm neck}$
by admissible mutations has a computable wavefunction,
conjectured to be the generating function of
all-genus open Gromov-Witten invariants of the corresponding Lagrangian.}
\label{fig:schema}
\end{figure}

Let $\cP$ be the symplectic cluster variety of Borel-decorated, $PGL_2$ local systems on a punctured
sphere $S$ with unipotent monodromy around the punctures.  There is a Lagrangian
subvariety $\cM\subset \cP$ of decorated
local systems with trivial monodromy at the punctures.
Cluster charts $\cP_\Gamma$ of $\cP$ are labeled by cubic graphs $\Gamma$ on $S$,
or dually ideal triangulations.  They are algebraic tori, and can be identified
with a torsor over
rank-one local systems on a
genus-$g$ Legendrian surface $S_\Gamma$ in the five-sphere:
after choosing a base point, we can write the chart as
$\cP_\Gamma \cong H^1(S_\Gamma;\bC^*)\cong (\bC^*)^{2g}.$
Then $\cM_\Gamma := \cM\cap \cP_\Gamma$ is a subspace of a torus
closely related to the space of graph colorings of $\Gamma,$ so we call
$\cM$ the \emph{chromatic Lagrangian}.
This symplectic torus chart has a canonical quantization to a quantum torus, and the chromatic Lagrangian $\cM_\Gamma$ quantizes to an ideal, ${\mathcal I}_\Gamma.$
Quantization is compatible with the cluster structure in the 
sense that the defining equations of ideals in mutated charts are related by mutations --- see Theorem \ref{thm:d-mod}.
The explicit description
of $\cM_\Gamma$ (and its quantization ${\mathcal I}_\Gamma$) will lead to enumerative predictions, but will also depend on
further choices: a \emph{phase}, a \emph{framing} and a \emph{cone}.

Concretely, as we describe in some more detail below,
$\cP_\Gamma$ quantizes to a standard quantum torus $\cT_q$, inside of which ${\mathcal I}_\Gamma$ is an ideal
defined by relations, one for each face of $\Gamma.$  Then
the wavefunction $\Psi_\Gamma$ is defined as
the unique power-series
solution to the equations having the form $1 + \cdots$.
The definition depends on additional data in each seed
allowing to construct the power-series representation of
$\cT_q$.
As we will show, these data depend on some choices,
after which the wavefunction is determined by its value in an initial cluster seed.  This gives an effective method to compute it.

Central to the strategy for calculation is to understand the effects
of a mutation $\Gamma \rightsquigarrow \Gamma'$,
which is dual to a flip of a triangulation, and to understand its interaction
with phases, framings and cones.
The entire structure is captured by the \emph{framed seed groupoid},
an enhancement of the cluster groupoid, whose arrows are either
mutations or changes of the various decorations ---
see Figure \ref{fig:schema}.

\subsection{Framed Seed Groupoid}
\label{sec:groupoid}

In a bit more detail, the edge lattice $\Lambda := \bZ^{E_\Gamma}$
of a cubic graph on an oriented surface (for us, the sphere)
has a natural skew form $(*,*)$
defined from the cyclic structure on edges meeting at a vertex.
Quotienting by its kernel $\Lambda_c$ defines a symplectic lattice
$\underline{\Lambda}$.  Roughly, a framed seed is an identification
of this lattice with the standard symplectic lattice $\bZ^g\oplus \bZ^g.$
More formally, it is a tuple
$(\mathbf{i},K,\widetilde{\mathbf{t}},\mathbf{f}),$
where $\mathbf{i}$ is a cluster seed (i.e. a basis for a lattice $\Lambda$ equipped with an integral skew form),
$K\subset \underline{\Lambda}$ is a maximal isotropic sublattice,
$\widetilde{\mathbf{t}}:\Lambda\rightarrow\mathbb{Z}$ is a character of $\Lambda$, and $\mathbf{f} = (\sigma,\{a_i\})$
is a pair of a splitting $\sigma: K^\vee\to \underline{\Lambda}$
of $0\to K\to \underline{\Lambda}\to K^\vee \to 0$ together with a basis $\{a_i\}$ of $K^\vee.$
Note we have $K^\vee \cong \underline{\Lambda}/K$
via the symplectic structure on $\underline{\Lambda},$
and we thereby obtain a dual basis $\{b_i\}$.  In total,
the data of the framed seed provides an identification of
$\underline{\Lambda}$ with the standard symplectic lattice
$\bZ^g \otimes \bZ^g$ with symplectic form $\sum_i du_i \wedge dv_i$.

We are interested in the set-up detailed in \cite{TZ},
i.e.~the construction of a Legendrian surface
$S_\Gamma \subset S^5$
from the data of $\Gamma\subset S$, and a
singular exact Lagrangian filling $L_0$ of $S_\Gamma$ as defined by an ideal
\emph{foam}, $\bF$, the combinatorial dual
of a tetrahedronization of a ball.  A smoothing $L$ can be defined by studying the
local model of the Harvey-Lawson special Lagrangian smoothing of the singular
Harvey-Lawson cone, and amounts to local choice of one of the
three possible face-matchings at each tetrahedron.
This geometry gives rise to a framed seed as follows: the group $H_1(S_\Gamma,\mathbb{Z})$ is identified with $\underline{\Lambda}$,
with its intersection form,
and $K$ is the kernel of the homology
push-forward of inclusion of the boundary $S_\Gamma \hookrightarrow L$.
Then $\{a_i\}$ is a basis for $H^1(L);$ the dual
basis $\{b_i\}$ for $H_1(L)$ and these give
coordinates $U_i$ and $V_j$ for $H^1(L,\bC^*)\cong \cP_\Gamma$, respectively.
The quantization then leads to an
isomorphism of $\cP_\Gamma$ with the quantum torus
$V_i U_j = q^{2\delta_{i,j}} U_j V_i.$
Each edge is then labeled by a monomial $X_e$ in the $U_i$ and $V_j$ with
$q$-dependent coefficient.

\subsection{Wavefunctions}
\label{sec:intro-wavefunction}

After quantization in
each chart $\cP_\Gamma$,
the Lagrangian subvariety $\cM_\Gamma \subset \cP_\Gamma$
becomes a left ideal $\cI_\Gamma$, and we can identify the
left $\cP_\Gamma$-module $\cP_\Gamma/\cI_\Gamma$ with
the principal ideal defined by cyclic vector $\Psi_\Gamma \in \bC[[\{X_i\}]],$
satisfying $\cI_\Gamma \Psi_\Gamma = 0$ in the standard representation
defined by exponentiating the Weyl representation:
$(U_i\cdot f)(X) = X_if(x),$ $(V_i\cdot f)(X) = f(q^2X_i),$ where ${q = e^{\pi i \hbar^2}}.$
The generators for $\cI_\Gamma$ are relations determined by the faces of $\Gamma$,
giving us concrete $q$-difference equations for $\Psi_\Gamma.$
For example, in the case where $\Gamma$ is the tetrahedron graph,
$S_\Gamma$ is a genus-one surface and
the quantum torus has generators $U$
and $V$ obeying ${VU = q^2UV}.$  For a certain choice of framed seed (see Figure~\ref{fig:alg-reframed-canoe} and Lemma~\ref{lem:reframed-canoe}),
the face equations are all equivalent to
$(1 -U -V)\Psi = 0$, and the unique power-series solution is $\Psi = \Phi(-q^{-1}X),$
where $\Phi(x) = \prod_{n\geq 0}(1+q^{2n+1}x)^{-1}$ is a
quantum dilogarithm.

The equations for $\cI_\Gamma$ are compatible with mutations
$\Gamma \rightsquigarrow \Gamma'$, meaning generators of $\cI_{\Gamma'}$ are
related to generators of $\cI_{\Gamma}$ by a cluster coordinate transformation,
and these are effected (up to a known basis change) by conjugation by a
quantum dilogarithm.  
The upshot is that graph mutations change
the wavefunction by the action of the quantum dilogarithm,
and as long as can make sense of this action on the ring of power series,
we may compute the resulting wavefunction.  We call such mutations \emph{admissible.}
In Lemma~\ref{lem:grpd-rep}, we show that the action on formal power series yields an \emph{algebraic representation} of the admissible part of the framed seed groupoid, meaning that in addition to  mutation, we can as well effect changes of other data of a framed seed
(phase, framing, basis) by known operators.
Moreover, the necklace graph $\Gamma_{\rm neck}$
(see Figure \ref{fig:necklacecanoe}) is a distinguished
base point for the framed seed groupoid, with known wavefunction 
$\Psi_{\Gamma_{\rm neck}} \equiv 1.$
So we can find any wavefunction for any point on the framed seed
groupoid connected to this basepoint by an admissible path.

One must check that the resulting wavefunction is independent of path, and
this amounts to checking that the cluster modular group (the automorphisms
of the standard quantum torus determined by loops in the groupoid) acts
trivially on the necklace wavefunction.  This can be verified explicitly
by observing that the necklace wavefunction is uniquely determined by
the defining equations for the ideal.

In this way, the cluster structure of the cluster modular groupoid
can be exploited to find wavefunctions.  Some have conjectural interpretations.

\subsection{Open Gromov-Witten Conjectures}
\label{sec:intro-ogw}

The cubic planar graphs $\Gamma$ that label cluster charts $\cP_\Gamma$
also describe Legendrian
surfaces $S_\Gamma$, which form asymptotic boundary conditions for
categories of A-branes,
by which we mean categories of constructible sheaves
with singular support on $S_\Gamma$ \cite{N,NZ}. 
Non-exact Lagrangian fillings $L\subset \bC^3$
asymptotic to $S_\Gamma$ have open Gromov-Witten invariants which
we conjecture, following the pioneering work of Aganagic-Vafa \cite{AV},
are predicted by the geometry of the brane moduli space $\cM_\Gamma
\subset \cP_\Gamma$.

The classical geometry conjecturally leads to open Gromov-Witten invariants.
The subvariety $\cM_\Gamma \subset \cP_\Gamma$ is Lagrangian.  Choosing
a framed seed $A$ and lifting to the universal
cover, we get $\widetilde{\cM}_{\Gamma} \subset \bC^{2g},$ and any connected
component determines a potential $W_\Gamma$ so that $\widetilde{\cM}_\Gamma$
is the graph of $dW_\Gamma$.
The instanton part of $W_\Gamma$ is conjectured to be the open Gromov-Witten
generating function.  

\begin{remark}
    Before stating the conjecture, we must acknowledge that there is no accepted definition of Open Gromov-Witten invariants, due to the
    necessity of choosing a homotopy class of nonvanishing section on the boundary of moduli space in order to construct a virtual fundamental
    cycle, hence a deformation-invariant count.    In \cite{L}, this was done in the case where the Lagrangian was invariant under a circle action,
    given a lift of that action to the torus boundary.  In \cite{ST}, an alternate approach involving bounding pairs was used to define
    open Gromov-Witten invariants for rational cohomology spheres.  Our conjecture below is subordinate to the construction of a rigorously
    defined invariant.
\end{remark}

\noindent{\bf Conjecture:}  
$W^{(A)}_\Gamma$ is the generating function
of disk invariants and obeys Ooguri-Vafa integrality: $W^{(A)}_\Gamma(X) =
\sum_{d\in \bZ_{\geq 0} \setminus \{0\}} n^{(A)}_d {\rm Li}_2(X^d),$
with $n^{(A)}_d \in \bZ.$ 
This conjecture appeared in essentially the same form in \cite[Section 1.2]{TZ}.

The cluster variety $\cP$ has a quantization, each chart $\cP_\Gamma$
of which can be identified, through a framing, with a quantum torus, $\cD$:
$V_i U_i = {q^2} U_i V_i,$ where ${q = e^{i\pi\hbar^2}}.$  Then
$\cM_\Gamma$ quantizes as an ideal $\cI$, and the left $\cD$-module
$\cD/\cI\cD$ is cyclic for a vector $\Psi_\Gamma$. Thanks to general results of Kontsevich-Soibelman~\cite{KS}, it follows that the wavefunctions $\Psi_\Gamma$ we construct satisfy the Ooguri-Vafa integrality property \cite{OV}: namely, they admit factorizations
\begin{equation}
\label{eq:OV-intro}
{\Psi_\Gamma^{(A)} =
\prod_{d \in \bZ_{\geq 0}^g} \prod_{k\in \bZ}\Phi((-q)^kX^d )^{n^{(A)}_{d,k}},}
\end{equation}
where $n^{(A)}_{d,k}\in \bZ$, and for fixed degree $d$ only finitely many of these integers $n^{(A)}_{d,k}$ are nonzero.


\begin{quote}
{\bf Conjecture:}  $\Psi^{(A)}_\Gamma$ is the generating function
of all-genus open Gromov-Witten invariants. 
(See Conjecture \ref{conj:ogw}
for details.)
\end{quote}
The conjecture implies the one above from \cite{TZ}
since
$\Psi_\Gamma^{(A)} \sim e^{-W^{(A)}_\Gamma/g_s}$ and $\Phi(x)\sim e^{-{\rm Li}_2(x)/g_s}$ as the string coupling constant $g_s = 2\pi i \hbar^2$ tends to $0$, and then $n^{(A)}_d = \sum_k n^{(A)}_{d,k}.$

Since all ideal triangulations are related by flips, every cubic planar graph of genus $g$ (meaning it has $2g+2$ vertices) can be obtained from $\Gamma^{\rm neck}_g$ through a sequence
of mutations. 
Our rubric therefore leads to conjectures for Lagrangian fillings for many
Legendrian surfaces.

\begin{remark}
As explained in~\cite{FG2}, the symplectic form on $\cP$ arises as the image under the regulator map $f\wedge g\mapsto d\log(f)\wedge d\log(g)$ of a canonical element $\mathcal{W}\in K_2(\mathbb{Q}(\cP))$ in the Milnor $K_2$-group of the field of rational functions on $\cP$. In \cite{DGGo} it is shown that the chromatic Lagrangian $\cM$ is in fact a $K_2$\emph{-Lagrangian}: the image of the K-theory class $\mathcal{W}$ under the restriction map $K_2(\mathbb{Q}(\cP))\rightarrow K_2(\mathbb{Q}(\cM))$ vanishes.
%
%
In fact, this
 $\mathrm{K}_2$ Lagrangianicity of $\cM$ is formally implied by the Ooguri-Vafa  integrality~\eqref{eq:OV-intro} of the wavefunction.\footnote{For example, locally in a cluster chart in the two-dimensional case, the regulator map sends $X\wedge Y$ to $d\log X\wedge d\log Y.$  Ooguri-Vafa integrality of the potential $W = \sum_d n^{(A)}_d \mathrm{Li}_2(X^d),$ says that $Y = \prod_d (1-X^d)^{-d n^{(A)}_d}.$ To see that this implies $\mathrm{K}_2$ Lagrangianicity, note $X^d\wedge (1-X^d)$ vanishes by the Steinberg relations.  Combined with the other relation $(ab)\wedge c = a\wedge c + b\wedge c$, it follows that $X\wedge Y$ restricts to zero.  The general case is proven similarly.}
\end{remark}

\subsection{Analytic Aspects}

A quantization in the physical sense would require that we construct,
in addition to wavefunctions
for each seed of the cluster modular groupoid,
a Hilbert space with arrows acting by
unitary isomorphisms.
Fock and Goncharov constructed such a quantization
depending on a parameter $\hbar \in \bR$, a central character for
the kernel of the skew form, with reality
being crucial for each logarithmic cluster variable $x$ to act in a unitary
way, and for mutations to be effected by a unitary action of
the Faddeev (noncompact quantum) dilogarithm $\varphi(x).$

Such an approach cannot work for us, as the unipotency condition defining
our cluster variety requires the central character to act as an \emph{imaginary}
number, ruling out self-adjointness in the na\"ive sense.
Nevertheless, in Section \ref{sec:analyticaspects} we present what we
think of as good evidence for the existence of a
quantization in the analytic sense, and for a well-defined wavefunction at each
seed.  Solutions are symmetric in $\hbar \leftrightarrow \hbar^{-1}$,
reflecting the symmetry of the ``squashed three-sphere" in the physical set-up
(see, e.g., \cite[Equation (2.16)]{CEHRV}).
In this set-up, all seed arrows would be admissible.  For example, mutating at all three strands
of the genus-two necklace graph $\Gamma_{\rm neck}^2$ would not be admissible
in the algebraic set-up of Section \ref{sec:intro-wavefunction}, but leads to an analytic wavefunction.  Indeed,
in Section \ref{sec:analyticaspects}
we show in this and several other
cases that different paths to the same framed seed
lead to the same wavefunction.  The identities needed to establish this path-independence (e.g.~\eqref{lem:23}) are consequences of the analytic properties
of the Faddeev dilogarithm and its Fourier self-duality. As an illustration of the analytic set-up, in Section~\ref{sec:cube} we show how it reproduces the all-genus analog of the proposal in~\cite{TZ} for the superpotential associated to the $g=3$ cubic graph given by the 1-skeleton of the cube.

\subsection{Framing Duality}



\begin{figure}
\begin{tikzpicture}
\pgfmathsetmacro{\g}{5}
\pgfmathsetmacro{\gmin}{\g-1}
\pgfmathsetmacro{\gmid}{\g/2}
\pgfmathsetmacro{\gone}{\g/3+1/3}
\pgfmathsetmacro{\gtwo}{\g/3+\g/3+2/3}
\pgfmathsetmacro{\size}{.33}
\pgfmathsetmacro{\nodesize}{.06}
\pgfmathsetmacro{\min}{.7}
\node at (\gmid,-\min/2) {$\Gamma_{\rm{neck}}^g$};
\draw[fill] (\g-\size,.5) circle (\nodesize cm);
\draw[fill] (\g+\size,.5) circle (\nodesize cm);
\foreach \i in {0,...,\gmin}
{\draw[thick] (\i,.5) circle (\size cm);
\draw[fill] (\i-\size,.5) circle (\nodesize cm);
\draw[fill] (\i+\size,.5) circle (\nodesize cm);
\draw[thick] (\i+\size,.5) -- (\i+1-\size,.5);}
\draw[thick] (\g,.5) circle (\size cm);
\draw[thick] (0-\size,.5) to[out=180,in=180] 
(0-\size/2,-\min) to[out=0,in=180] (\gone,-\min)
to[out=0,in=180] (\gtwo,-\min) 
to[out=0,in=180] (\g+\size/2,-\min) to[out=0,in=0] (\g+\size,.5);
\end{tikzpicture}
\qquad
\begin{tikzpicture}
\pgfmathsetmacro{\g}{5}
\pgfmathsetmacro{\seatsize}{.8}
\pgfmathsetmacro{\gmin}{\g-1}
\pgfmathsetmacro{\gmid}{\g*\seatsize/2+.5}
\pgfmathsetmacro{\size}{.06}
\node at (\gmid,-.4) {$\Gamma_{\rm{canoe}}^g$};
\draw[fill] (0,\seatsize/2) circle (\size cm);
\draw[thick] (0,\seatsize/2) -- (1,0);
\draw[thick] (0,\seatsize/2) -- (1,\seatsize);
\draw[fill] (1+\gmin*\seatsize+1,\seatsize/2) circle (\size cm);
\draw[thick] (1+\gmin*\seatsize+1,\seatsize/2) -- (1+\gmin*\seatsize,0);
\draw[thick] (1+\gmin*\seatsize+1,\seatsize/2) -- (1+\gmin*\seatsize,\seatsize);
\draw[thick] (1+\gmin*\seatsize,0) -- (1+\gmin*\seatsize,\seatsize);
\draw[fill] (1+\gmin*\seatsize,0) circle (\size cm);
\draw[fill] (1+\gmin*\seatsize,\seatsize) circle (\size cm);
\foreach \i in {1,...,\gmin}
{
\pgfmathsetmacro{\ival}{1+\i*\seatsize-\seatsize}
\draw[thick] (\ival,0)--(\ival,\seatsize);
\draw[thick] (\ival,0)--(\ival+\seatsize,0);
\draw[thick] (\ival,\seatsize)--(\ival+\seatsize,\seatsize);
\draw[fill] (\ival,0) circle (\size cm);
\draw[fill] (\ival,\seatsize) circle (\size cm);
}
\draw[thick] (0,\seatsize/2) to[out=270,in=180] (\gmid,-.8) to[out=0,in=270] (1+\gmin*\seatsize+1,\seatsize/2);
\end{tikzpicture}
\caption{Mutating the necklace graph $\Gamma_{\rm{neck}}^g$ (left) along the
$g$ short strands results in the canoe graph $\Gamma_{\rm{canoe}}^g$ (right).
Here $g=5.$  The Legendrian surfaces generalize the Chekanov and Clifford tori,
respectively, which arise when $g=1.$}
\label{fig:necklacecanoe}
\end{figure}
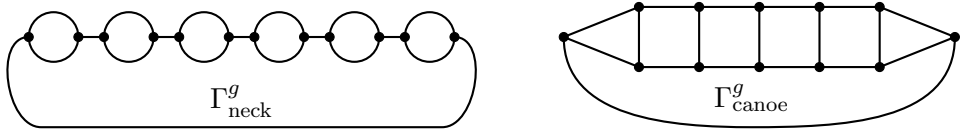

We notice a curious identity between wavefunctions and quiver invariants.
A special role is played by the Legendrian Clifford torus
and its higher-genus generalizations.
These Clifford surfaces of genus $g$ arise from ``canoe'' graphs (see
Figure \ref{fig:necklacecanoe}).
The Clifford surfaces arise from mutations of the higher-genus version of the Chekanov torus, a genus-$g$ Legendrian surface corresponding to a ``necklace'' graph (see again Figure \ref{fig:necklacecanoe}).  Each Chekanov surface
has a distinguished exact Lagrangian filling and therefore 
a distinguished phase and no holomorphic disks:  
$\Psi_{\Gamma_{\rm{neck}}^g}\equiv 1$ and $W_{\Gamma_{\rm{neck}}^g}\equiv 0$.
After mutation,
we get a distinguised phase for the Clifford surface $S_{\Gamma_{\rm{canoe}}^g}$,
i.e.~a Lagrangian filling $L$ with $\partial L = S_{\Gamma_{\rm{canoe}}^g}$
and $b_1(L)=g$, though
the different framings in this phase are parametrized by a $g\times g$ symmetric
integral matrix, $A$. 
The corresponding
wavefunction $\Psi^{(A)}_{\Gamma_{\rm{canoe}}^g}$
can be computed from cluster theory as in Section \ref{sec:groupoid},
and, as stated in Section \ref{sec:intro-ogw} above,
is conjecturally the partition function of the open topological string.

%


We can now state \emph{framing duality} in the following way.
Let $A$ be a $g\times g$ symmetric, integral matrix with non-negative entries.
Let $Q_A$ be the symmetric quiver with $g$ nodes and adjacency matrix $A$.
Recall that the DT series is
the generating function for cohomologies of quiver representation spaces $M_d$
(see Section \ref{sec:dtsersforsymmquivs} for
a precise definition)
over different dimension vectors, $d$.  Explicitly,
$\mathrm{DT}_{A} = \sum_{d\in \bZ_{\geq 0}} \sum_{s\in \bZ}(-1)^k H^s(M_d)X^d t^{k/2}.$
Then we have, after setting $t^{\frac{1}{2}} = -q$:
\begin{equation}
\label{eq:framingduality}
\text{The wavefunction is the DT series of $Q_A$: \qquad $\Psi_{\Gamma_{\rm canoe}^g}^{(A)} =
\;\mathrm{DT}_{A}$ }.
\end{equation}
\vskip0.1in

Further, explicit formulas show that the invariants $n^{(A)}_d$ relate to refined Kac polynomials of quivers, as defined in \cite{RV}.
Recall that the Kac polynomial $A_d(q)$ of a quiver $Q$ counts the number of isomorphism classes of absolutely
indecomposable representations of $Q$.
The refined Kac polynomials $A_{\lambda}(q)$ are labeled by partitions $\lambda,$ and satisfy $A_d(q) = \sum_{|\lambda|=d}A_\lambda(q).$
A special role will be played by $\lambda = 1^d = (1,1,...,1).$
In Proposition \ref{prop:kac-prop}, we show that when $Q$ is the quiver with one node and $h\geq 1$ arrows,
and if $A$ is the one-by-one matrix $(2-2h),$ then $A_{1^d}(1) = n^{(A)}_d.$


\begin{remark}
Many of the results which establish this equality were
performed by Kontsevich-Soibelman in \cite{KS}.
For the genus-one case studied by Aganagic and Vafa, the connection between DT invariants and open GW invariants in different framings  was observed also in \cite{LZ}.  
As for other Legendrians, also in genus-one,
wavefunctions for knot and link conormals were considered in \cite{AENV}. 
Finding quiver duals for knot conormals is known as the Knot-Quiver Correspondence
\cite{KRSS}.
The relationship \eqref{eq:framingduality}
suggests that the quiver invariants arise from an effective
quiver quantum-mechanical theory described by the capping data for
the noncompact threefolds we construct from Harvey-Lawson
components --- see, e.g.~\cite[Section 5.1.1]{CEHRV}.
Framing duality is thus
in the spirit as the knots-quiver correspondence of \cite{KRSS},
whose geometric and physical interpretations were proposed in \cite{EKL}.
It is however more general, in the following sense.
The Legendrian surfaces considered here are higher genus and not tori,
giving rise to all symmetric quivers and DT invariants depending on all $g$ variables. 
In contrast, framings of a fixed knot
are labeled by a single integer, corresponding to a one-parameter
set of quivers, with DT invariants determined by
specializing the $g$ variables to a one-dimensional slice ---
see \cite[Equation (4.2)]{KRSS}.
It would be interesting to pursue a geometric interpretation of
framing duality along the lines of \cite{EKL}.

\end{remark}

\begin{remark}
One wonders if the above relations extend to other cubic graphs and/or
nonsymmetric quivers.
\end{remark}

\subsection{Seminal Prior Works}


Very similar constructions were considered from related physical perspectives
in prior works.
In \cite{CEHRV} and \cite{DGGo} the authors consider an M5-brane on $S^3\times L,$ where $L$ is a Lagrangian
submanifold of a compactifying space.  (Those authors call this Lagrangian $M$.)
They describe the partition function of the effective 3d theory on $S^3$ as a quantum-mechanical state. 
The M-theory set-up expresses this partition function as an integer combination of dilogarithms.
The partition function can also be computed by reduction to $L$.
It is a general property of quantum field theory that 
the path integral on a manifold with boundary always defines
a state in the Hilbert space defined by the boundary.
In the present case, the boundary Hilbert space is a quantization of the space of
flat $\mathrm{U}(1)$ connections on the genus-$g$ Legendrian boundary surface
$S_\Gamma$ {(or a torsor over such --- see Section \ref{sec:cpgafm})}.
The wavefunction $\Psi$ should be understood as the wavefunction of this
quantum state.

On top of all this, many of the results of this paper have also appeared in important
previous works, to which we owe a debt of gratitude.
In \cite{CCV} and \cite{CEHRV}, the authors studied the
behavior of these wavefunctions under symplectic transformations, although
not via cluster theory and without relating
the results to Gromov-Witten invariants.
The papers
\cite{DGGo} and \cite{DGGu} overlap with the present paper,
as well as \cite{CEHRV}, in considering Lagrangian
double covers branched over tangles, and
studied the corresponding Lagrangian moduli space.
The paper \cite{KS} studied quiver representations
and preservation of integrality under changes of framings, providing
many of the key formulas that we use.
The idea of quantizing mirror curves goes back to \cite{ADKMV}
and has been integral to the spectral approach of \cite{GHM},
applications to knot polynomials in \cite{GS},
and difference equations for partition functions in \cite{T,NT}.
Finally, the relation of wavefunctions to open Gromov-Witten invariants appeared previously
in \cite{TZ},
as well as in \cite{ES}, though these references
did not explore the interaction with cluster theory.
(The paper \cite{Za} outlined the strategy employed in the present paper, but without details.)
Further citations are made in the text.

\subsection*{Acknowledgements}

We dedicate this paper to Steve Zelditch, our late colleague and friend.
A generous giant of a mathematician, Steve clarified several
analytical and representation-theoretic issues we confronted
in preparing this paper.
We are greatly endebted to David Treumann, who was involved in a
significant part of this collaboration.
It is a pleasure to thank Roger Casals, Melissa Liu, Lenny Ng, {Piotr Su\l kowski},
and Boris Tsygan for helpful conversations.
We thank Peng Zhou for asking about mutations very early in this project. 
L.S.~has been supported by NSF grant DMS-2200738.
E.Z.~has been supported by NSF grants DMS-1406024,
DMS-1708503 and DMS-2104087.

%% file: Cluster.tex
\section{Cluster Poisson Varieties and Quantizations}
\label{sec:cluster}
For the convenience of the reader, we briefly recall the needed background on cluster Poisson varieties and their quantizations. 
Within this paper, we focus on the  cluster Poisson varieties that are skew-symmetric and without frozen variables. A more general definition of cluster Poisson varieties can be found in \cite{FG2}.

\subsection{Cluster Poisson varieties}

\begin{definition} \label{seed.cluster}
A seed is a pair  ${\bf i}=(\{x_1, \ldots, x_n\} , W)$, where $\{x_1, \ldots, x_n\}$ is a collection of commuting algebraically independent variables, and
$
W=\sum_{i, j} a_{ij} x_i \frac{\partial }{\partial x_i} \wedge x_j \frac{\partial }{\partial x_j} 
$
 is a bi-vector encoded by an integer skew-symmetric matrix $A=(a_{ij})$. 
 Correspondingly, we get a quiver $Q_A$ such that its vertices are labelled by $1$ through $n$ and the number of arrows from $i$ to $j$ is $[a_{ij}]_+:= \max\{0, a_{ij}\}$.
\end{definition}

Let ${\bf i}$ be a seed.  Every  $k\in \{1,\ldots, n\}$ creates a new seed $\mu_k({\bf i})=(\{x_1',\ldots, x_n'\}, W)$ such that 
\[
x_i'= \left\{  \begin{array}{ll} 
      x_k^{-1} & \mbox{if } i=k, \\
      x_i(1+x_k^{-{\rm sgn}(a_{ik})})^{-a_{ik}} & \mbox{if } i\neq k. \\
   \end{array} \right.
\]
In terms of $\{x_i'\}$, the bi-vector $W$ can be presented as $\sum_{i, j} a_{ij}' x_i'\frac{\partial }{\partial x_i'}\wedge x_j'\frac{\partial }{\partial x_j'}$, where 
\[a_{ij}'=\left\{\begin{array}{ll}
-a_{ij} & \mbox{if $i=k$ or $j=k$;}\\
a_{ij}+\frac{|a_{ik}|a_{kj}+a_{ik}|a_{kj}|}{2}& \mbox{otherwise}.
\end{array}\right.\]
The process of obtaining the new seed  $\mu_k({\bf i})$ is called a {\it cluster mutation} in the direction $k$. The cluster mutation $\mu_k$ in the same direction is involutive: $\mu_k^2({\bf i})={\bf i}$. 

\smallskip 

Let $\sigma$ be a permutation of $\{1, \ldots, n\}$. It gives rise to a  seed $\sigma({\bf i})=(\{x_1',\ldots, x_n'\}, W)$ such that
\[
x_i'= x_{\sigma^{-1}(i)}, \hskip 14mm  i\in \{1, \ldots, n\}.
\]
A composition $\tau= \sigma\circ \mu_{i_1}\circ \cdots \circ \mu_{i_m}$ of cluster mutations and permutations taking a seed ${\bf i}$ to ${\bf i}'$ is called a cluster transformation.

\begin{definition}
Let $\cX$ be a rational variety over $\bC$ equipped with a rational bi-vector $W$. A cluster chart of $\cX$ is a birational map 
\[
\pi=(x_1, \ldots, x_n): ~ \cX \lra \bC^n
\]
such that ${\bf i}_\pi:=(\{x_1, \ldots, x_n\}, ~\pi_*(W))$ forms a seed. Two cluster charts  are  called equivalent if their corresponding seeds  are related by a cluster transformation. The equivalence class of a cluster chart $\pi$ is denoted by $|\pi |$.

Abusing notation\footnote{Within this paper, we only take into account the birational structure of $\cX$.}, a variety $\cX$ equipped with a pair $(|\pi|, W)$ is called a cluster Poisson variety.  

\end{definition}

Let $\mathbb{C}(\cX)$ be the field of rational functions on $\cX$. For a cluster chart $\pi'=\{x_1', \ldots, x_n'\}$, let 
$\mathcal{T}_{\pi'}=\mathbb{C}[x_1'^{\pm 1}, \ldots, x_n'^{\pm 1}] \subset \mathbb{C}(\cX)$ denote the ring of Laurent polynomials in $x_1',\ldots, x_n'$.
The  {\it cluster Poisson algebra} is the intersection
\begin{equation}
\label{upper.Poisson.algebra}
\mathbb{L}_{\cX}:=\bigcap_{\pi'\in |\pi|} \mathcal{T}_{\pi'}.
\end{equation} 
Note that the bivector $W$ induces a natural Poisson bracket on  $\mathbb{L}_{\cX}$:
\[
\{\cdot, \cdot\}: ~ \mathbb{L}_{\cX} \times \mathbb{L}_{\cX} \longrightarrow \mathbb{L}_{\cX}, \qquad \{f, g\}:= W(f,g).
\]

\smallskip 

Let $p$ be a birational automorphism of $\cX$. We say $p$ is a cluster automorphism  if 
\begin{itemize}
\item $p$ preserves the bi-vector: $p_\ast(W)=W$,
\item $p$ preserves the equivalence class of cluster charts: $\pi\circ p \in |\pi|$.
\end{itemize}
The set of cluster automorphisms forms a group. Denote it by $\mathcal{G}_\cX$ and call it the {\it cluster modular group} of $\cX$. The group $\mathcal{G}_\cX$ acts by Poisson automorphisms on the  algebra $\mathbb{L}_{\cX}$.

\subsection{Quantization}
\label{subsec:quant}
Let $\cX$ be a cluster Poisson variety.
Let $A=(a_{ij})$ be the $n\times n$ integer skewsymmetric matrix appearing in an initial seed defining $\cX$ as in Definition \ref{seed.cluster}. To $A$ is associated a triple $(\Lambda, \Pi, (\ast , \ast ))$, where $\Lambda$ is a rank  $n$ lattice, $\Pi=\{e_1, \ldots, e_n\}\subset \Lambda$ is a basis, and $(\ast , \ast)$ is a bilinear form on $\Lambda$ such that $( e_i, e_j ) =a_{ij}$. 
We also set
\[
\Lambda^+=\bigoplus_{i=1}^n \mathbb{Z}_{\geq 0} e_i, \hskip 14mm \Lambda^- = \bigoplus_{i=1}^n \mathbb{Z}_{\leq 0} e_i.
\]
Let $\mathbb{C}[q^{\pm1}]$ be the ring of Laurent polynomials in $q$.
Let $\mathcal{T}^q$ be the quantum torus algebra  over $\mathbb{C}[q^{\pm 1}]$ with the generators $X_v$ ($v\in \Lambda$), subject to the relations
\begin{equation}
\label{quantum.relations}
X_v X_w = q^{(v, w)} X_{v+w}.
\end{equation}
Denote by ${\bf Frac}(\mathcal{T}^q)$ the non commutative field of fractions of $\mathcal{T}^q$ (cf.~\cite[Appendix]{BZ}).
The positive cone $\Lambda^+$ determines a formal completion of the algebra $\mathcal{T}^q$. We will consider the group of formal power series with leading term 1
\[
\widehat{\mathcal{R}}_\Pi=\Big\{\sum_{v\in \Lambda^+} a_v(q) X_v ~\Big|~ a_0(q)=1, ~~a_v(q)\in \mathbb{C}(q)\Big\}.
\]

Now let us consider the mutations of the basis $\Pi=\{e_1,\ldots, e_n\}$. Let $\Pi^*=\{\alpha_1, \ldots, \alpha_n\} \subset \Lambda^\ast$ be the dual basis of $\Pi$. Let $k\in \{1,\ldots, n\}$.
For an $n$-tuple $S=\{v_1, \ldots, v_n\}$  of elements in $\Lambda$, the mutated  $\mu_k(S)=\{v_1',\ldots, v_n'\}$ consists of elements
\begin{equation}
\label{basis.mutatation}
v_i'=  \left\{\begin{array}{ll} 
      -v_k & \mbox{if } i=k, \\
       {\displaystyle v_i+ \sum_{l=1}^n \max\{0, (v_i, v_k) \alpha_l(v_k)\} \cdot {\rm sgn}(\alpha_l(v_k)) e_l } & \mbox{if } i\neq k.  \\
   \end{array}\right.
\end{equation}
\begin{remark} There is a slightly more general version of mutations, which we will consider in Section \ref{groupoid}.
\end{remark}

Let $(k_1, \ldots, k_m)$ be a sequence of indices in $\{1,\ldots, n\}$. Let us start with the set $\Pi=S$. Applying the mutations \eqref{basis.mutatation} recursively, we obtain a sequence of bases of $\Lambda$ 
\begin{equation}
\label{mutation.quiver}
\Pi=\Pi_1\stackrel{\mu_{k_1}}{\lra} \Pi_2\stackrel{\mu_{k_2}}{\lra}\ldots \stackrel{\mu_{k_m}}{\lra} \Pi_{m+1} =\Pi', \hskip 14mm \mbox{where }  \Pi_j=\{e_{1}^{(j)}, \ldots, e_n^{(j)}\}.
\end{equation}
A basis $\Pi'$ obtained  this way is said to be equivalent to $\Pi$. Let  $|\Pi|$ consist of bases equivalent to $\Pi$.

The elements $e_i^{(j)}$ in \eqref{mutation.quiver} are called  $c$-vectors by Fomin-Zelevinsky \cite{FZ4}.
The sign coherence of $c$-vectors asserts that each $e_{i}^{(j)}$ lies  either in $\Lambda^+$ or $\Lambda^-$  \cite{DWZ}.
Hence there is a unique sequence of signs $(\varepsilon_{1}, \ldots, \varepsilon_{m}) $ such that
\begin{equation}
\label{f.vectors}
f_j= \varepsilon_{j} e_{k_j}^{(j)} \in \Lambda^+, \hskip 14mm  j=1, \ldots, m.
\end{equation}
We define the formal power series 
\[
\Phi_{\Pi'}= \Phi(X_{f_1})^{\varepsilon_1}\Phi(X_{f_2})^{\varepsilon_2}\cdots \Phi(X_{f_m})^{\varepsilon_m} ~~~ \in ~~~ \widehat{\mathcal{R}}_\Pi.
\]
where we recall that
\[
\Phi(X)= \prod_{n=0}^\infty (1+q^{2n+1}X)^{-1}
\]
is the (compact) quantum dilogarithm function. The formal power series $\Phi(X)$ is a close relative of the infinite {\it q-Pochhammer} symbol 
\begin{align}
\label{Pochhammer.infty}
(x; q^2)_\infty:&=\prod_{n=0}^{\infty}(1-q^{2n}x)\\
&=1+\sum_{k=1}^{\infty}\frac{(-1)^k q^{k(k-1)}}{\prod_{i=1}^k (1-q^{2i})}x^k 
\nonumber\\
&= \exp\Big(\sum_{k=1}^\infty \frac{x^k}{k(q^{2k}-1)}\Big) \nonumber\\
&= \Phi(-q^{-1}x)^{-1}    \hskip 2cm \in \mathbb{Z}((q))[[x]]. \nonumber
\end{align}
The latter is  the unique formal power series starting from 1 and satisfying the  difference relation
\begin{equation}
\label{difference.rel}
(x; q^2)_\infty =(1-x)\cdot (q^2x; q^2)_\infty.
\end{equation}
For $m\in \mathbb{Z}$, we define the finite {\it q-Pochhammer symbol} by
\[
(x;q^2)_m:= \frac{(x;q^2)_\infty}{(q^{2m}x; q^2)_\infty} 
\]
We have the following fundamental result, which guarantees that the series $\Phi_{\Pi'}$ is a well-defined function of the set $\Pi'$:
\begin{theorem}[{\cite[Th.4.1]{K}}] The power series $\Phi_{\Pi'}$ only depends on the set $\Pi'$, not on the mutation sequences that take $\Pi$ to $\Pi'$. 
\end{theorem}
Associated with each $\Pi'=\{e_1',\ldots, e_n'\}\in |\Pi|$ is a quantum torus algebra $\mathcal{T}_{\Pi'}^q$ over $\mathbb{C}[q^{\pm 1}]$ with generators
\[
X_v'= {\rm Ad}_{\Phi_{\Pi'}} (X_v) ~ \in {\bf Frac}(\mathcal{T}^q), \hskip 14mm v\in \Lambda.
\]
The generators  $X_v'$ satisfy the relations \eqref{quantum.relations}. In particular, the variables $X_{e_1'}',\ldots, X_{e_n'}'$ are called {\it quantized cluster $\mathcal{X}$-variables.}
The pair $(\Pi', \mathcal{T}_{\Pi'}^q)$ is called a quantum cluster seed. The {\it quantum cluster algebra} is the intersection
\begin{equation}
\label{quantum.cluster.a}
\mathbb{L}^q_\cX=\bigcap_{\Pi'\in |\Pi|} \mathcal{T}_{\Pi'}^q \subset {\bf Frac}(\mathcal{T}^q).
\end{equation}
The quasiclassical limit $q\mapsto 1$ of \eqref{quantum.cluster.a} recovers the Poisson algebra \eqref{upper.Poisson.algebra}. 

The cluster modular group $\mathcal{G}_\cX$ acts on $\mathbb{L}^q_{\mathscr{X}}$ via {\it quantum cluster automorphism}, constructed as follows. Every element in $\mathcal{G}_\cX$ one-to-one corresponds to a linear automorphism $\tau$ of the lattice $\Lambda$
such that $\tau$ preserves the  bilinear form on $\Lambda$ and maps the initial basis set $\Pi$ to $\Pi':=\tau(\Pi) \in |\Pi|$. Each $\tau$ gives rise to an algebra isomorphism
\[
g_\tau: \mathcal{T}_{\Pi'}^q \stackrel{\sim}{\lra} \mathcal{T}_{\Pi}^q,\hskip 14mm X_v'\longmapsto X_{\tau^{-1}(v)}.
\]
The restriction of $g_\tau$ on $\mathbb{L}^q_\cX$ induces an algebra automorphism of $\mathbb{L}^q_\cX$, called a quantum cluster automorphism.

\subsection{Casimirs} The bilinear form $(\ast, \ast)$ on $\Lambda$ gives rise to a linear map $c$ from $\Lambda$ to its dual $\Lambda^\ast$
\[
\forall v\in \Lambda, \hskip 10mm c(v)(\ast)=(v, \ast).
\]
The  kernel of $c$ forms a sub-lattice $\Lambda_c$  of $\Lambda$. The quotient $\Lambda / \Lambda_c$ is a symplectic lattice.

If $v\in \Lambda_c$, then $X_v$ commutes with every generator $X_w$ by \eqref{quantum.relations}. For  every $\Pi'\in |\Pi|$, we have
\[
X_v'= {\rm Ad}_{\Phi_{\Pi'}}(X_v)= X_v. 
\]
Therefore $X_v$ $(v\in \Lambda_c)$ are contained in $\mathbb{L}_\cX^q$ and are called {\it Casimirs}.  
It is easy to see that the center $Z(\mathbb{L}_\cX^q)$ of $\mathbb{L}_\cX^q$ is the torus algebra  generated by Casimirs. 

\begin{definition}
Let ${\bf t}$ be  a homomorphism from $Z(\mathbb{L}_\cX^q)$ to $\mathbb{C}[q^{\pm1}]$.
 The quotient algebra  $\mathbb{L}_{\cX, {\bf t}}^q$ of $Z(\mathbb{L}_\cX^q)$ is obtained by modulo the relations 
\[
 X_v={\bf t}(X_v),\] 
 where $v$ goes through $\Lambda_c$.
\end{definition}

\subsection{Moduli space of \texorpdfstring{$G$}{G}-local systems}
\label{sec:FG-moduli-space}
Let $G$ be a split semisimple algebraic group over $\mathbb{Q}$ with trivial center. Let ${S}$ be an oriented compact topological surface with $n$ {\it punctures} $p_1, \ldots , p_n$  removed. Denote by $\chi(S)$ the Euler characteristic of ${S}$. We require that 
\[
n> \max\{0, \chi(S)\}
\]
so that $S$ admits a triangulation whose vertices are the punctures.
The Fock-Goncharov moduli space $\mathscr{X}_{G, S}$, introduced in \cite{FG1}, provides an important class of cluster Poisson varieties. Below we briefly recall the definition and several basic properties of $\mathscr{X}_{G, S}$ for later use.
\smallskip 

We start with a local model. The flag variety $\mathcal{B}$ parametrizes the Borel subgroups of $G$. Recall the {\it Grothendieck-Springer} resolution
\[
\widetilde{G}:=\left\{ (g, B) \in G\times \mathcal{B} ~|~ g\in B\right\}.
\]
The projection from $\widetilde{G}$ to $\mathcal{B}$ makes $\widetilde{G}$ a smooth $B$-bundle over $\mathcal{B}$. Let ${H}$ be the Cartan subgroup of $G$. For each Borel subgroup ${B}\in \mathcal{B}$, there is a canonical group homomorphism 
\begin{equation}
\label{pr.h}
\pi_{B}:~ {B}\longrightarrow {B}/[{B}, {B}] \stackrel{\sim}{\longrightarrow} {H}.
\end{equation}
Consequently, we get a regular map
\[
\pi:~\widetilde{G}\longrightarrow H, \qquad (g, B) \longmapsto \pi_B(g).
\]
The variety $\widetilde{G}$ carries a Poisson structure such that $\pi$ is a symplectic fiberation. 
For example, see \cite{EL} for more details on the Poisson geometry of $\widetilde{G}$. 
An element $g\in G$ is unipotent if and only if $\pi_B(g)=1$. The subvariety
\[
\widetilde{\mathcal{N}}:=\pi^{-1}(1) \subset \widetilde{G}
\]
is the usual Springer resolution of the unipotent cone $\mathcal{N}\subset G$. Note that $\widetilde{\mathcal{N}}$ is naturally isomorphic to the cotangent bundle $T^*\mathcal{B}$. Therefore it admits a symplectic structure, although we caution the reader that this is
not the same as the one determined by the cluster structure associated with the model of once-punctured disk. Its zero section consists of elements $(1,B)$ for all $B\in \mathcal{B},$ and is a Lagrangian subvariety of $\widetilde{\mathcal{N}}$.

\vskip 2mm

Now we generalize the above construction to the moduli space of $G$-local systems. 
\begin{definition} 
 A framed $G$-local system over $S$ consists of the data $({\mathcal L}, \{B_1, \ldots , B_n\})$ where 
\begin{itemize}
\item ${\mathcal L}\in {Hom}(\pi_1(S), G)$ is a $G$-local system over $S$;
\item $B_i$ is a flat section of the associated bundle $\mathcal{L}\times_G \mathcal{B}$ over the loop around the puncture $p_i$.
\end{itemize}
The moduli space $\cX_{G, S}$ consists of the framed $G$-local systems modulo the conjugation of $G$.   
\end{definition}

\begin{theorem} 
\label{cluster.x}
The space $\cX_{G, S}$ is a cluster Poisson variety. The mapping class group of $S$ acts on $\cX_{G, S}$ via cluster Poisson transformations.
\end{theorem}
\begin{remark} The cluster Poisson structure on $\cX_{G, S}$ has been constructed by Fock and Goncharov \cite[\S 9]{FG1} for $G={PGL}_{r+1}$, by Le \cite{Le} for $G$ being a classical group, and finally by Goncharov and
Shen \cite{GS2} for an arbitrary semisimple group. Theorem 1 of \cite{S} further shows that the ring of regular functions $\mathcal{O}(\cX_{G, S})$ is a cluster Poisson algebra and therefore admits a quantization.
\end{remark}

\begin{example}
Let $G=PGL_2$ and let $\mathcal{T}$ be an ideal triangulation of $S$, i.e., a triangulation whose vertices are the punctures. For simplicity, we shall avoid self-folded triangles. We place a vertex at the center of every edge in $\mathcal{T}$. Within each triangle in $\mathcal{T}$, we add three arrows in the counter-clockwise orientation, as shown in Figure \ref{figurecl01}. In this way, we obtain a quiver $Q_{\mathcal{T}}$.
\begin{figure}[ht]
 \begin{tikzpicture}[scale=1.5]
 \draw (0,1) -- (1,0) -- (0, -1) -- (0, 1) --(-1,0)--(0,-1);
  \node[red]  at (0,1.2) {$a$};
 \node[red]  at (-1.2,0) {$d$};
 \node[red]  at (0,-1.2) {$c$};
 \node[red]  at (1.2,0) {$b$};
  \node[blue, label=above: $x$] (x) at (0,0) {$\bullet$};
 \node[blue, label=right: $t$] (t) at (0.5,0.5) {$\bullet$};
 \node[blue, label=right: $u$] (u) at (0.5,-0.5) {$\bullet$};
 \node[blue, label=left: $w$] (w) at (-0.5,0.5) {$\bullet$};
  \node[blue, label=left: $v$] (v) at (-0.5,-0.5) {$\bullet$};
  \draw[blue, thick, -latex] (t) edge (x) (x) edge (u) (u)edge (t) (x) edge (w) (w) edge (v) (v) edge (x);
  \begin{scope}[shift={(4,0)}]
   \draw (-1,0) -- (1,0) -- (0, 1) -- (-1, 0) --(0,-1)--(1,0);
    \node[red]  at (0,1.2) {$a$};
 \node[red]  at (-1.2,0) {$d$};
 \node[red]  at (0,-1.2) {$c$};
 \node[red]  at (1.2,0) {$b$};
  \node[blue, label=above: $x'$] (x') at (0,0) {$\bullet$};
 \node[blue, label=right: $t'$] (t') at (0.5,0.5) {$\bullet$};
 \node[blue, label=right: $u'$] (u') at (0.5,-0.5) {$\bullet$};
 \node[blue, label=left: $w'$] (w') at (-0.5,0.5) {$\bullet$};
  \node[blue, label=left: $v'$] (v') at (-0.5,-0.5) {$\bullet$};
  \draw[blue, thick, -latex] (t') edge (w') (w') edge (x') (x')edge (v') (v') edge (u') (u') edge (x') (x') edge (t');
  \end{scope}
  \end{tikzpicture}
  \caption{A cluster structure associated with $\mathscr{X}_{PGL_2, S}$}
  \label{figurecl01}
\end{figure}
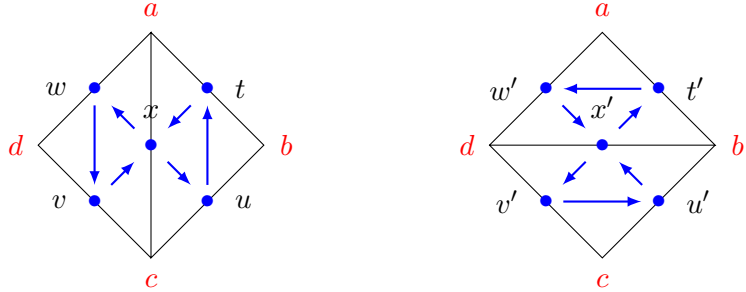
Note that $\mathcal{B}_{PGL_2}=\mathbb{P}^1$. Each framed local system  in $\mathscr{X}_{PGL_2, S}$ assigns a quadruple $a,b,c,d\in \mathbb{P}^1$ to the vertices of each quadrilateral in $\mathcal{T}$. We define the cluster variable placed on the diagonal of the quadrilateral to be the cross ratio
\[
x=-\frac{(a-b)(c-d)}{(b-c)(d-a)}.
\]
In this way, we obtain a cluster seed ${\bf i}_{\mathcal{T}}$ for $\mathscr{X}_{PGL_2, S}$. 

As  in Figure \ref{figurecl01}, a flip of each edge gives rise to a cluster mutation, whose new variables become
\[
x'=x^{-1}, \qquad t'=t(1+x), \qquad u'=u(1+x^{-1})^{-1},
\qquad 
v'=v(1+x),\qquad w'=w(1+x^{-1})^{-1},
\]
and the rest of the variables are invariant. 
\end{example}

For general $G$, following \eqref{pr.h}, the flat section $B_i$ chosen for each puncture $p_i$ gives rise to a map from ${\cX}_{G,S}$ to the Cartan subgroup $H$. 
Therefore we get a map
\begin{equation}
\label{eq:GS-pi}
\pi=(\pi_1, \ldots, \pi_n): {\cX}_{G,S} \lra {H}^n.
\end{equation}
By Theorem 2.10 of \cite{GS2}, the fibers of $\pi$ are symplectic varieties. In particular, for each simple positive root $\alpha$ of $G$, the regular function $\alpha\circ \pi_i$ is a Casimir of $\cX_{G,S}$.
Let us set 
\[{\cX}_{G, S}^{\rm un}:= \pi^{-1}(1).
\]
We have
\[
\dim{\cX}_{G, S}^{\rm un} = \dim{\cX}_{G, S}-n\dim H= 2n \dim \mathcal{B} -\chi(S)\dim G.  
\]

\subsection{Example: the sphere cases.}
\label{sec:sphere-cases}
Within this subsection, we assume that $S$ is a sphere with $n$ punctures. As illustrated by Figure \ref{figurecl}, we have
\begin{equation}
\label{acn123}
\cX_{G, S}^{\rm}=\left\{\left((u_1, B_1), (u_2,B_2), \ldots, (u_n, B_n)\right)~\middle|~ (u_i,B_i) \in \widetilde{\mathcal{N}}, ~\prod_{i=1}^n u_i=1\right\} {\Big \slash} G
\end{equation}
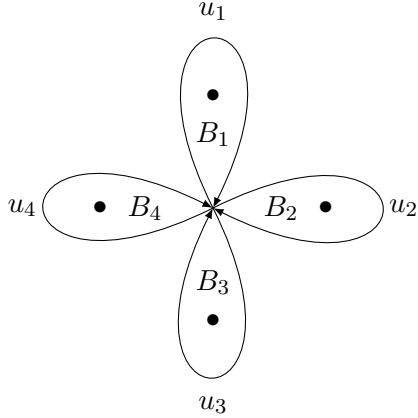
\begin{figure}[ht]
 \begin{tikzpicture}[scale=1.5]
 \node[label=below: ${B_1}$] at (0,1) {$\bullet$};
 \node[label=left: ${B_2}$] at (1,0) {$\bullet$};
 \node[label=above: ${B_3}$] at (0,-1) {$\bullet$};
 \node[label=right: ${B_4}$] at (-1,0) {$\bullet$};
  \node[label=below: ${u_1}$] at (0,2) {};
 \node[label=left: ${u_2}$] at (2,0) {};
 \node[label=above: ${u_3}$] at (0,-2) {};
 \node[label=right: ${u_4}$] at (-2,0) {};
\draw[-latex] (0,0) .. controls (-1,2) and (1,2) .. (0,0);
\draw[ -latex] (0,0) .. controls (2,1) and (2,-1) .. (0,0);
\draw[-latex] (0,0) .. controls (1,-2) and (-1,-2) .. (0,0);
\draw[ -latex] (0,0) .. controls (-2,-1) and (-2,1) .. (0,0);
  \end{tikzpicture}
  \caption{A framed local system on a sphere with 4 punctures.}
  \label{figurecl}
  \end{figure}

Let $G_{sc}\rightarrow G$ the simply connected covering of $G$. The center $Z(G_{sc})$ coincides with the kernel of the covering map.  
Let $d$ be the determinant of the Cartan matrix of $G$, as given in the following table
\begin{center}
\begin{tabular}{|c| c| c| c| c| c |c |c|} 
 \hline
 & $A_r$ & $B_r$ & $C_r$ & $D_r$ & $E_r$~($r=6,7,8$) & $F_4$ & $G_2$  \\ 
 \hline
 $d$ & $r+1$ & 2 & 2 & 4 & $9-r$ &1 &1  \\ 
 \hline
\end{tabular}.
\end{center}
It is known that the order of $Z(G_{sc})$ is $d$.

\begin{proposition}  \label{prop:components}
Let $S$ be a sphere with at least three punctures. The space ${\cX}_{G, S}^{\rm un}$ has $d$ many top dimensional irreducible components. 
\end{proposition}
\begin{proof}
Every unipotent element $u\in G_{sc}$ has a unique lift to a unipotent element $\tilde{u}\in G_{sc}$. Then the product condition in \eqref{acn123} becomes
\[
\prod_{i=1}^n \tilde{u}_i \in Z(G). 
\]
Accordingly, we obtain a decomposition
\[
\cX_{G,S}^{un} =\bigsqcup_{c\in Z(G)} \cX_{G,S}^{un}(c),
\]
where $\cX_{G,S}^{un}(c)$ consists of the points such that $\prod_{i=1}^n \tilde{u}_i=c$.

Now we show that every $\cX_{G,S}^{un}(c)$ contains a unique top dimensional irreducible component. Let $D$ be a disk with $n-2$ punctures and $2$ marked points on its boundary. Following \cite[Definition 2.4]{FG1}, the moduli space $\mathscr{A}_{G_{sc}, D}$ parametrizes the decorated twisted unipotent $G$-local systems on $D$. Each boundary interval of $D$ corresponds to an invariant in the Cartan subgroup of $G_{sc}$, denoted by $h$ and $h'$ respectively as in the following figure. 
\begin{figure}[ht]
\begin{tikzpicture}
\draw (0,0) circle (2cm);
\node[label=left: $1$] at (0,2) {$\bullet$};
\node[label=left: $2$] at (0,1) {$\bullet$};
\node[label=left: $n-1$] at (0,-1) {$\bullet$};
\node at (0,0) {$\vdots$};
\node[label=left: $n$] at (0,-2) {$\bullet$};
\draw[thick, blue, -latex] (60:2) .. controls (0,1.5) .. (120:2);
\draw[thick, blue, -latex] (150:2) .. controls (1.5,2) and (1.5,0) .. (150:2);
\draw[thick, blue, -latex] (210:2) .. controls (1.5,0) and (1.5,-2) .. (210:2);
\draw[thick, blue, -latex] (240:2) .. controls (0,-1.5) .. (300:2);
\node[red] at (-2.2,0) {$h'$};
\node[red] at (2.2,0) {$h$};
\node[blue] at (70:1.6) {$u_1$};
\node[blue] at (1.2,1) {$u_2$};
\node[blue] at (1.2,-1) {$u_{n-1}$};
\node[blue] at (-70:1.6) {$u_n$};
\end{tikzpicture}
\end{figure}

As constructed in \cite{GS2}, the space $\mathscr{A}_{G_{sc}, D}$ carries a cluster $K_2$ structure, with $2{\rm rk}(G)$ many frozen variables, given by $\omega_i(h)$ and $\omega_i(h')$ respectively, where $\omega_i$ are the fundamental weights  of $G_{sc}$.  
We impose an extra condition that $h'=1$ and $h\in Z(G_{sc})$, obtaining a subspace $\mathscr{A}_{G_{sc}, D}' \subset \mathscr{A}_{G_{sc}, D}$. Depending on the value of $h$, we get a decomposition
\[
\mathscr{A}_{G_{sc}, D}'=\bigsqcup_{h\in Z(G_{sc})} \mathscr{A}_{G_{sc}, D}'(h).
\]
Here every component $\mathscr{A}_{G_{sc}, D}'(h)$ is rational, with the usual cluster coordinates for the mutable ones, and a specialization on the frozen ones. 

Since $h'=1$ and $h\in Z(G_{sc})$, when passing from $G_{sc}$ to $G$, one may identify the pinnings given by the two boundary intervals, obtaining a map 
\[
\pi:\mathscr{A}_{G_{sc}, D}' \longrightarrow \mathscr{X}_{G, S}^{un}.
\]
More precisely, recall the central element $s_{G_{sc}}$ as in Corollary 2.1 of \cite{FG1}. By comparing the geometric meanings of both spaces, we see that $\pi$ maps $\mathscr{A}_{G_{sc},D}'(s_{G_{sc}}^n\cdot c)$ to $\mathscr{X}^{\rm un}_{G, S}(c)$. 

Now we fix a  simple path $\gamma$ on the sphere $S$ connecting the puncture 1 and n. Given a generic point in $\mathscr{X}^{\rm un}_{G, S}(c)$, we may choose a decoration for each of the flags $B_1,...B_{n-1}$. 
Let us cut along the path $\gamma$, obtaining the disk $D$. Finally, we choose a decoration for $B_n$ such that $h'=1$. 
In this way, we obtain a lift of the generic point in $\mathscr{X}^{\rm un}_{G, S}(c)$ to  $\mathscr{A}_{G_{sc},D}'(s_{G_{sc}}^n\cdot c)$. Through the process, we see that the map $\pi$ is dominant, with a fiber isomorphic to ${H}_{sc}^{n-1}$ for every generic point in $\mathscr{X}^{\rm un}_{G, S}(c)$.  As a consequence, we get the desired dimension 
\[
\dim \cX_{G,S}^{un}(c) = \dim \mathscr{A}_{G_{sc}, D}' - \dim H_{sc}^{n-1}= 2n \dim \mathcal{B} - 2 \dim \G. 
\] \qedhere
\end{proof}

Now let $i$ be a reflection of $S$ that fixes the punctures. For example, if $S$ is a sphere, then one can put all the punctures on the equator, and $i$ exchanges the two hemispheres. Note that $i$ changes  the orientation of $S$.
Therefore $i$ induces an anti-Poisson involution of ${\mathcal X}_{G,S}$. Let $s$ be the inverse map of $H^n$ which takes $(h_1, ..., h_n)$ to $(h_1^{-1}, \ldots, h_{n}^{-1})$. By definition, the following maps commute

\begin{center}
\begin{tikzpicture}
\node (a) at (0,0) {${\cX}_{G,S}$};
\node (b) at (2,0) {${H}^n$};
\node [label=above: $\pi$] at (1, -0.2) {};
\node [label=above: $\pi$] at (1, -1.7) {};
\node [label=left: $i$] at (0.2, -0.75) {};
\node [label=right: $s$] at (1.8, -0.75) {};
\node (c) at (0, -1.5) {${\cX}_{G,S}$};
\node (d) at (2,-1.5) {${H}^n$};
\draw [->] (a) --(b);
\draw [->] (c) -- (d);
\draw [->] (a) -- (c);
\draw [->] (b) --(d);
\end{tikzpicture}
\end{center}
Therefore $i$ maps ${\cX}_{G, S}^{\rm un}$ to ${\cX}_{G, S}^{\rm un}$.

Taking all the fixed points of the map $i$, we get a subvariety ${\cM}_i$ of ${\cX}_{G, S}^{\rm un}$.

\begin{theorem} 
\label{thm:general.lagrangian}
${\cM}_i$ is a Lagrangian subvariety of ${\cX}_{G, S}^{\rm un}$. 
\end{theorem}

\begin{proof} Let $w$ be the symplectic form on ${\cX}_{G, S}^{\rm un}$. Note that $i^*(w)=-w$. Since $i$ is the identity map on $\cM_i$, the restriction of $w$ to $\cM_i$ is trivial. It remains to check the dimension: $\dim \cM_i= \frac{1}{2} \dim {\cX}_{G, S}^{\rm un}$.
\end{proof}


\input{Globalization.tex}

%% file: Globalization.tex
\begin{example} \label{rmk:special.case} 
Let $G=PGL_2$.
The following triangulations show an example of involution for a sphere with 4 punctures.

\begin{center}
\begin{tikzpicture}[scale=0.9]
\node at (0,0) {$\bullet$};
\node at (2,0) {$\bullet$};
\node at (0,2) {$\bullet$};
\node at (2,2) {$\bullet$};
\draw (0,0)--(2,0)--(2,2)--(0,2)--cycle;
\draw [dashed] (0,0) --(2,2);
\draw (2,0)--(0,2);
\node at (-0.2, 1) {$w$};
\node at (1, 2.2) {$x$};
\node at (2.2, 1) {$y$};
\node at (1, -0.2) {$z$};
\node at (1.5, 1.3) {$u$};
\node at (0.5, 1.3) {$v$};
\node at(4.25, 1.2) {${\rm reflection}$};
\draw[thick, -latex] (3,1)--(5.5,1);
\end{tikzpicture}
\begin{tikzpicture}[scale=0.9]
\node at (0,0) {$\bullet$};
\node at (2,0) {$\bullet$};
\node at (0,2) {$\bullet$};
\node at (2,2) {$\bullet$};
\draw (0,0)--(2,0)--(2,2)--(0,2)--cycle;
\draw (0,0) --(2,2);
\draw[dashed] (2,0)--(0,2);
\node at (-0.2, 1) {$w^{-1}$};
\node at (1, 2.2) {$x^{-1}$};
\node at (2.2, 1) {$y^{-1}$};
\node at (1, -0.2) {$z^{-1}$};
\node at (1.5, 1.3) {$u^{-1}$};
\node at (0.5, 1.3) {$v^{-1}$};
\draw[thick, -latex] (3.5,1)--(7,1);
\node at(5.2, 1.2)  {cluster~ mutations};
\end{tikzpicture}
\begin{tikzpicture}[scale=0.9]
\node at (0,0) {$\bullet$};
\node at (2,0) {$\bullet$};
\node at (0,2) {$\bullet$};
\node at (2,2) {$\bullet$};
\draw (0,0)--(2,0)--(2,2)--(0,2)--cycle;
\draw[dashed] (0,0) --(2,2);
\draw (2,0)--(0,2);
\node at (-0.2, 1) {$\tau(w)$};
\node at (1, 2.2) {$\tau(x)$};
\node at (2.2, 1) {$\tau(y)$};
\node at (1, -0.2) {$\tau(z)$};
\node at (1.5, 1.3) {$v$};
\node at (0.5, 1.3) {$u$};
\end{tikzpicture}
\end{center}
Here
\[
\tau(u)=v;  \hskip 14mm \tau(v)=u.
\]
\[
\tau(w)=w^{-1}(1+v^{-1})(1+u^{-1}); \hskip 14mm \tau(y)=y^{-1}(1+v^{-1})(1+u^{-1}); 
\]
\[
\tau(x)=x^{-1}(1+v)^{-1}(1+u)^{-1}; \hskip 14mm \tau(z)=z^{-1}(1+v)^{-1}(1+u)^{-1}; 
\]
Note that the mapping class group of punctured sphere acts on $\mathcal{P}$ by symplectomorphisms. The mapping class group preserves $\cM$, but it interchanges the other components of $\cM_i$.

In general, $\cM$ in Theorem \ref{thm:comparison} is a connected component of $\cM_i$. Therefore Theorem \ref{thm:comparison} is a special case of Theorem \ref{thm:general.lagrangian}.

\end{example}

%% file: groupoid.tex
\section{Groupoids of polarized and framed seeds}
\label{groupoid}
In this section we define the groupoid of framed seeds, an enhancement of the standard cluster modular groupoid that we shall use to describe concrete models for representations of the corresponding cluster variety.  

\subsection{Polarizations and framings for seeds}

\label{sec:polarizationsandframings}

Suppose the rank of the skew-form $(\ast, \ast)$ associated to the seed $\mathbf{i}$ is $2g$, and write  $\Lambda_c\subset\Lambda$ for its kernel. In what follows, we will write $\underline\Lambda :=\Lambda/\Lambda_c$ for the corresponding rank-$2g$ symplectic lattice, which fits into the short exact sequence
\begin{align}
    \label{eq:SES-center}
    0\rightarrow \Lambda_c\rightarrow\Lambda\rightarrow\underline\Lambda\rightarrow0.
\end{align}

A \emph{polarization} for $\mathbf{i}$ is the choice of an isotropic sublattice $K\subset \underline{\Lambda}$ of maximal rank $g$, such that the skew form induces a short exact sequence of lattices
\begin{align}
\label{eq:polarization-SES}
0\rightarrow K\rightarrow \underline\Lambda\rightarrow K^\vee\rightarrow 0.
\end{align}

We consider two polarized seeds $(\mathbf{i},K)$ and $(\mathbf{i}',K')$ to be equivalent if the canonical map  $\underline\Lambda_{\mathbf{i}}\rightarrow \underline\Lambda_{\mathbf{i}'}$ is an isometry which sends $K$ to $K'$. If $(\mathbf{i},K)$ is a polarized seed and $\mathbf{i}'$ is a seed related to $\mathbf{i}$ by a signed mutation or permutation, then the induced isomorphism  of symplectic lattices $\iota:\underline{\Lambda}_{\mathbf{i}}\simeq\underline{\Lambda}_{\mathbf{i}'}$  determines a polarization $K'=\iota(K)$ for $\mathbf{i}'$. 

Our reason for introducing the additional data of polarizations is that they define representations of the symplectic torus $\mathcal{T}_{\underline\Lambda}^q$ associated to the seed $\mathbf{i}$. Indeed, a polarization $K$  for $\mathbf{i}$ determines a commutative subalgebra $\mathcal{T}_{K}^q\subset\mathcal{T}_{\underline\Lambda}^q$. The subalgebra $\mathcal{T}_{K}$ is identified with the coordinate ring of a split algebraic torus of rank $g$, and let us write $\mathbf{1}_K$ for its 1-dimensional representation given by evaluation at the identity element. From the latter we may construct an induced representation of $\mathcal{T}_{\underline\Lambda}^q$:
\begin{align*}
\mathcal{V}_{K} : &= \mathrm{Ind}_{\mathcal{T}_{K}^q}^{\mathcal{T}_{\underline\Lambda}^q}\left(\mathbf{1}_K\right)\\
&=\mathcal{T}_{\underline\Lambda}^q\otimes_{\mathcal{T}_{K}^q}\mathbf{1}_K.
\end{align*}

The representation $\mathcal{V}_{K}$ is a $\mathbb{Z}$-module of infinite rank. In order to give a concrete model for it, it is necessary to equip the polarized seed $(\mathbf{i},K)$ with another piece of additional data, which we now describe.

\begin{definition}
\label{def:framedseed}
A \emph{framing} $\mathbf{f}$ for a polarized seed $(\mathbf{i},K)$ is the following data:
\begin{enumerate}
    \item a basis $\{a_i\}$ for $K^\vee$;
    \item  \label{dat:polsplit} a splitting $s:K^\vee\rightarrow \underline\Lambda$ of the short exact sequence~\eqref{eq:polarization-SES}, such that the image of $K^\vee$ in $\underline\Lambda$ is isotropic; and
    \item \label{dat:cc} a group homomorphism
$$
\widetilde{\mathbf{t}} : \underline{\Lambda} \rightarrow \mathbb Z.
$$
\end{enumerate}
\end{definition}
Let us now reformulate the notion of a framing for a seed in concrete terms. Consider the standard quantum torus
\[
\mathcal{D}_{2g}:= \mathbb{Z}[q^{\pm1}]\langle U_1^{\pm},\ldots, U_n^\pm, V_1^\pm, \ldots, V_n^\pm\rangle 
\]
with the relations
\[
U_i U_j = U_j U_i , \hskip 1cm V_i V_j = V_j V_i, \hskip 1cm V_i U_j = q^{2\delta_{ij}} U_j V_i
\]
The choice of a framing $\mathbf{f}$ for a polarized seed $(\mathbf{i},K)$  determines an isomorphism
\begin{equation}
\label{eq:iso-qtori}
 \underline{\iota}_{\mathbf{f}}: \mathcal{T}_{\underline\Lambda}^q \longrightarrow \mathcal{D}_{2g}
\end{equation}
which is uniquely characterized by the requirement that the element $X_{s(a_i)}$ of $\mathcal{T}_{\underline\Lambda}^q$ is mapped to the generator $U_i$ of $\mathcal{D}_{2g}$. The generators $V_i$ then correspond under the inverse isomorphism to elements of the basis $\{b_i\}$ of $K$ dual to the basis $\{a_i\}$ for $K^\vee$. Additionally, the data ~\eqref{dat:cc} of the homomorphism $\widetilde{\mathbf{t}}$ in the definition of a framing determines a surjection of quantum tori 
$$
 \mathcal{T}_{\Lambda}^q \longrightarrow \mathcal{T}_{\underline\Lambda}^q, \qquad X_{\lambda}\mapsto (-q)^{\widetilde{\mathbf{t}}(\lambda)}X_{\lambda+\Lambda_c}
$$
which factors through the central quotient of $\mathcal{T}_{\Lambda}^q$ by the double sided ideal $\langle X_z-(-q)^{\widetilde{\mathbf{t}}(z)} | z\in\Lambda_c\rangle$.
Putting everything together, we see that a framing $\mathbf{f}$ gives rise to a surjection of quantum tori
$$
\iota_\mathbf{f} :  \mathcal{T}_{\Lambda}^q\longrightarrow \mathcal{D}_{2g},    
$$
and that all the data of the framing and polarization can be uniquely recovered from that of the surjection $\iota_{\mathbf{f}}$.

Now let  $\mathcal{R}:= \mathbb{Z}[q^{\pm}][X^{\pm}_1, \ldots, X^{\pm}_g]$ be the ring of Laurent  polynomials in $g$ variables. Then there is a representation of $\mathcal{D}_{2g}$  on $\mathcal{R}$ such that
\begin{equation}
\label{eq:rep-def}
\forall F\in \mathcal{R}, \hskip 1cm U_i \cdot F = X_iF, \hskip 1cm V_i \cdot F = F(X_1, \ldots, q^2 X_i, \ldots, X_n),
\end{equation}
and we obtain an isomorphism of $\mathcal{T}_{\underline\Lambda}^q $-modules
$$
\mathcal{V}_K\simeq \underline{\iota}_{\mathbf{f}}^*\mathcal{R},
$$
thus providing the promised model for the induced representation $\mathcal{V}_K$.

A \emph{framed seed} $\underline{\mathbf{i}} $ is the data $(\mathbf{i},K,\mathbf{f})$ of a seed $\mathbf{i}$ together with a polarization $K$ and framing $\mathbf{f}$. We consider two framed seeds to be equivalent if the isomorphism of quantum tori $\mathcal{T}^q_{\Lambda_{\mathbf{i}}}\rightarrow \mathcal{T}^q_{\Lambda_{\mathbf{i}'}}$ induced by canonical map of lattices $\underline\Lambda_{\mathbf{i}}\rightarrow \underline\Lambda_{\mathbf{i}'}$ fits into a commutative diagram
$$
\begin{tikzcd}
\mathcal{T}^q_{\Lambda_{\mathbf{i}}}\arrow{rd}{\iota_{\mathbf{f}}}\arrow{dd}{a} &  \\
 &\mathcal{D}_{2g}\\
\mathcal{T}^q_{\Lambda_{\mathbf{i}'}}\arrow{ru}{\iota_{\mathbf{f}'}}
\end{tikzcd}
$$
\subsection{Operations on framed seeds}
Suppose that seeds $\mathbf{i},\mathbf{i}'$ are related by a signed mutation in direction $k$, so that we have an isometry of lattices $\nu_k^\pm:\Lambda_{\mathbf{i}'}\rightarrow\Lambda_{\mathbf i}$. If $K,K'$  and $\mathbf{f},\mathbf{f}'$ are polarizations and framings for  $\mathbf{i},\mathbf{i}'$, we say that the framed seeds $(\mathbf{i},K,\mathbf{f})$ and $(\mathbf{i}',K',\mathbf{f}')$ are related by the signed mutation in direction $k$ if $K =(\nu_{k}^\pm)^*(K')$, and similarly all pieces of framing data for $\mathbf{f}$ in Definition~\ref{def:framedseed} are identified with those for $\mathbf{f}'$ under the lattice isomorphism $\nu_k^\pm$. In particular, for any pair of framed seeds related by a signed mutation, there is a unique monomial map $\underline\nu^\pm:\underline\Lambda'\rightarrow \underline\Lambda$ such that the following diagram commutes:
$$
\begin{tikzcd}
\mathcal{T}^q_{\Lambda_{\mathbf{i}}}\arrow{r}\ar{d}{\nu^\pm} &\mathcal{T}^q_{\underline\Lambda_{\mathbf{i}}}\ar{d}{\underline\nu^\pm} \\
\mathcal{T}^q_{\Lambda_{\mathbf{i}'}}\arrow{r} &\mathcal{T}^q_{\underline\Lambda_{\mathbf{i}'}} 
\end{tikzcd}
$$
Recall that a framed seed $\underline{\mathbf{i}} $ gives rise to a symplectic basis $\{s(a_i),b_i\}$ for $\underline\Lambda$, where we again write $b_i$ for the elements of the basis for $K$ dual to the basis $\{a_i\}$ for $K^\vee$. 
We say that two framed seeds are related by a \emph{framing change morphism} if all pieces of the framing data are identical except for the datum~\eqref{dat:polsplit} given by the splitting $s$ of $\underline\Lambda$. The space of framing change morphisms based at a given framed seed is naturally identified with the space of $g\times g$ symmetric integer matrices $\Omega=(\omega_{ij})$, where the new splitting $s'$ is related to the original by
$$
s'(a_i) = s(a_i) + \sum_{j=1}^g\omega_{ij}b_j, \quad i=1,\ldots,g.
$$
\begin{remark}
\label{rmk:dual-basis-transformation}
We recall that if $a'_i = \sum C_{ij}a_j$ is another basis of $K^\vee$, then the corresponding dual basis is given by $b'_i = \sum (C^{-1})_{ji}a_j$. Hence the symmetric matrix $\Omega$ transforms under such a change of basis $C$ as
$$
\Omega \longmapsto C\Omega C^T.
$$
\end{remark}
Given a vector $\mathbf{d}=(d_1,\ldots, d_g)\in\mathbb{Z}^g$, consider the algebra automorphism $\sigma_\mathbf{d}$ of $\mathcal{D}_{2g}$
defined by
\begin{align}
\label{eq:sigma-transf}
    \sigma_\mathbf{d}(U_i) = (-q)^{d_i}U_i,\qquad \sigma_\mathbf{d}(V_i)=V_i.
\end{align}

We say that two framed seeds with identical underlying lattice $\Lambda$ are related by a \emph{coordinate rescaling} if the surjections $
\iota_\mathbf{f},\iota_\mathbf{f}' : \mathcal{T}_{\Lambda}^q\longrightarrow \mathcal{D}_{2g},
$ are related by $\iota_\mathbf{f}' = \sigma_{\mathbf{d}}\circ\iota_\mathbf{f}$ for some $\mathbf{d}\in\mathbb{Z}^{2g}$.

The \emph{framed seed groupoid} is a category whose objects are equivalence classes of framed seeds. The arrows are generated by those of four elementary kinds: signed mutations, permutations,  framing change morphisms, and coordinate rescalings. Each arrow $a:(\mathbf{i},K,\mathbf{f})\rightarrow (\mathbf{i}',K',\mathbf{f}')$ induces a birational automorphism of $\mathcal{D}_{2g}$: those corresponding to permutations, changes of framing, and coordinate rescalings induce the natural biregular automorphisms, and a signed mutation in direction $k$ induces a birational automorphism via the monomial isomorphism $\underline\nu_k^\pm$ and conjugation by $\Phi(\iota_{\mathbf{f}}(X_{\pm e_k}))^{\pm}$. We put a relation on the arrows in the framed seeds groupoid by identifying arrows with the same source and target which induce identical birational automorphisms of $\mathcal{D}_{2g}$.

\subsection{Framed seeds and representations}
\label{sec:sectiontitle}

Suppose that $\underline{\mathbf{i}}$ is a framed seed, and recall the corresponding representation 
$$
\mathcal{V}_{\underline{\mathbf{i}}}\simeq \mathbb{Z}[q^\pm][X_1^\pm,\ldots, X_g^\pm]
$$
 of the quantum torus $\mathcal{T}^q_{\underline\Lambda}$.
The embedding of the Laurent series ring into the ring 
$$
\mathcal{K}:=\mathbb{Z}((q))((X_1,\ldots, X_g))
$$  
of formal Laurent series also gives rise to a representation of $\mathcal{T}_{\underline\Lambda}^q$ which we denote by $\widehat{\mathcal V}_{\underline{\mathbf i}}$. 

For the purposes of constructing wavefunctions, it will be necessary to consider the action of a somewhat larger algebra on the representation $\widehat{\mathcal V}_{\underline{\mathbf i}}$. Write $\widehat{\mathcal{D}}_{2g}$ for the `complete quantum torus' associated to $\mathcal{D}_{2g}$, which may be regarded as the ring of non-commutative formal Laurent series in $U_i,V_i$. Inside $\widehat{\mathcal{D}}_{2g}$, consider the subalgebra
$$
\mathcal{A}_{2g} := \mathbb{Z}((q))((U_1,\ldots, U_g))\langle V_1^{\pm1},\ldots V^{\pm 1}_g\rangle
$$
consisting of formal Laurent series in the $U_i$ whose coefficients are Laurent polynomials in the $V_i$. Unlike in the case of $\widehat{\mathcal{D}}_{2g}$, there is a well-defined action of the algebra $\mathcal{A}_{2g}$ on $\widehat{\mathcal V}_{\underline{\mathbf i}}$. Indeed, under~\eqref{eq:rep-def} each $V_i$ acts on the `vacuum vector' $1\in \widehat{\mathcal V}_{\underline{\mathbf i}}$ by $V_i\cdot1 = 1$, and so the action of a arbitrary Laurent polynomial in the $V_i$, being a finite $\mathbb{Z}((q))$-linear combination of such, is also well-defined. 

Recall that the space of change of framing morphisms based at a given framed seed can be identified with the additive group $\mathfrak{p}_g$ of $g\times g$ symmetric matrices $\Omega=(\omega_{ij})$ with $\omega_{ij}\in \mathbb{Z}$. Its group algebra $\mathbb{Z}\mathfrak{p}_g$ is generated by symbols $T_\Omega, \Omega\in\mathfrak{p}_g$ satisfying $T_\Omega T_{\Omega'} = T_{\Omega+\Omega'}$. The group $\mathfrak{p}_g$ acts on $\mathcal{A}_{2g}$ by automorphisms called \emph{changes of framing}:
\begin{equation}
\label{frame.change.matrix}
T_\Omega: \mathcal{A}_{2g} \stackrel{\sim}{\lra} \mathcal{A}_{2g},  \hskip 1cm V_{j}\longmapsto V_j, \hskip 7mm U_j \longmapsto q^{\omega_{jj}} U_j \prod_{k=1}^g V_k^{\omega_{jk}},
\end{equation}
and we may form the semi-direct product algebra
$$
\widehat{\mathcal A}_{2g} = \mathcal{A}_{2g}\otimes_{\mathbb{Z}}\mathbb{Z}\mathfrak{p}_g.
$$
Given $U_{\mathbf w} = \prod_{j}U_j^{w_j}$, it follows from~\eqref{frame.change.matrix} that we have
$$
T_\Omega(U_{\mathbf w}) = q^{{\bf w}^t \Omega {\bf w}}U_{\mathbf w} V_{\Omega\mathbf w}.
$$

\begin{remark}
\label{rmk:framing-BCH}
The reader may find the following interpretation of the framing shift automorphisms useful. Consider the topological Heisenberg algebra $\mathcal{H}_g$ over $\mathbb{C}[[\hbar]]$ generated by $\{u_j,v_j\}$ subject to the relations
$$
[u_j,v_k] = \frac{\delta_{j,k}}{2\pi i },
$$
and set $q = e^{\pi i \hbar^2}$. The algebra $\mathcal{A}_{2g}$ embeds into this Heisenberg algebra via $U_k\mapsto e^{2\pi\hbar u_k},~V_k\mapsto e^{2\pi\hbar v_k}$. Now given a $g\times g$ symmetric matrix $\Omega\in\mathfrak{p}_g$, consider the associated quadratic form
$$
Q(\mathbf v) = \sum_{j,k=1}^g\omega_{jk}v_jv_k,
$$
and write $e^{-\pi i Q(\mathbf v)}$ for the corresponding element of the group algebra $\mathbb{Z}\mathfrak{p}_g$. Note that the $e^{-\pi i Q(\mathbf v)}$ are not elements of the Heisenberg algebra $\mathcal{H}_g$, but one can nonetheless formally compute the result of conjugating the generators of $\mathcal{H}_g$ by them using the Baker-Campbell-Hausdorff formula:
\begin{align*}
\mathrm{Ad}_{e^{-\pi i Q(\mathbf v)}}(u_j) &= u_j - \pi i [Q(\mathbf v),u_j]\\
&= u_j + \sum_k \omega_{jk}v_k,
\end{align*}
so that 
\begin{align*}
\mathrm{Ad}_{e^{-\pi i Q(\mathbf v)}}(U_j) &= \mathrm{Ad}_{e^{-\pi i Q(\mathbf v)}}(e^{2\pi \hbar u_j})\\
&= e^{2\pi \hbar(u_j+ \sum_k \omega_{jk}v_k)}\\
&=q^{\omega_{jj}}e^{2\pi \hbar u_j}e^{2\pi\hbar \sum_k \omega_{jk}v_k}\\
&= q^{\omega_{jj}} U_j \prod_{k=1}^g V_k^{\omega_{jk}},
\end{align*}
recovering~\eqref{frame.change.matrix}.

\end{remark}

%

The extended algebra $\widehat{\mathcal{A}}_{2g}$ also acts in the representation $\widehat{\mathcal V}_{\underline{\mathbf i}}
\simeq\mathcal{K}$: given $F= \sum_{{\bf w}} C_{\bf w}(q) X^{\bf w} \in \mathcal{K}$, we define
\begin{align}
\label{eq:frame-change-action}
T_\Omega\cdot F :=\sum_{{\bf w}} q^{{\bf w}^t \Omega {\bf w}} C_{\bf w}(q) X^{\bf w}.
\end{align}
That~\eqref{eq:frame-change-action} indeed defines a representation of the extended algebra $\mathcal{A}_{2g}$ follows easily from the considerations of Remark~\eqref{rmk:framing-BCH}, or can be readily verified directly. 
{Finally, let us remark that the coordinate-rescaling operators $\sigma_{\mathbf d}$ defined in~\eqref{eq:sigma-transf} also act naturally in the representation $\widehat{\mathcal V}_{\underline{\mathbf i}}$ via
\begin{align}
\label{eq:sigma-transf-action}
\sigma_{\bf d}\cdot F :=\sum_{{\bf w}} (-q)^{\bf d\cdot \bf w} C_{\bf w}(q) X^{\bf w}.
\end{align}
}

\subsection{{Admissible and primitive mutations}} 
\label{sec:groupoid-rep}
Suppose that $\underline{\mathbf{i}}$ is a framed seed, and $e_k$ is an element of the basis $\Pi$ for $\Lambda$ associated to the underlying seed $\mathbf{i}$. Recall that the data of the framing $\mathbf{f}$ allows us to associate to $\pm e_k$ a monomial $\iota_\mathbf{f}(X_k^\pm)\in\mathcal{D}_{2g}$ of the form
$$
\iota_\mathbf{f}(X_k^\pm) = (-q)^r\exp\left(2\pi\hbar{\sum_{j=1}^g{m_j}u_j + n_jv_j}\right), \qquad m_j,n_j,r\in\mathbb{Z},
$$
where we adopt the notations of Remark~\ref{rmk:framing-BCH}. 
We say that a mutation of the framed seed $\underline{\mathbf{i}}$ in direction $e_k$ with sign $\epsilon$ is \emph{admissible} if in the monomial $\iota_\mathbf{f}(X_k^\epsilon)$ we have $m_j\geq0$ for all $j=1,\ldots, g$, and in addition there is at least one $j$ for which $m_j\neq0$.
Let us make a few simple remarks about this definition. 
\begin{remark}
If two framed seeds $\underline{\mathbf{i}},\underline{\mathbf{i}}'$ are related by a change of framing, then evidently a signed mutation is admissible with respect to $\underline{\mathbf{i}}$ if and only if it is admissible with respect to $\underline{\mathbf{i}}'$.
\end{remark}
\begin{remark}
Let $a$ be an admissible mutation of framed seed $\underline{\mathbf{i}}$ in direction $k$ with sign $\epsilon$, and let $\underline{\mathbf{i}}' = a(\underline{\mathbf{i}})$ be the resulting framed seed.
Then the mutation of $\underline{\mathbf{i}}'$ in direction $k$ with sign $-\epsilon$, which is the inverse of $a$ in the framed seed groupoid, is also an admissible mutation.
\end{remark}
It follows from these remarks there is a sub-groupoid of the framed seeds groupoid whose morphisms are generated by framing shifts and admissible mutations.

Our reason for introducing the notion of admissibility of mutations is the following: a mutation of a framed seed in direction $e_k$ with sign $\epsilon$ is admissible (if and) only if  the  quantum dilogarithm formal power series $\Phi\left(\iota_\mathbf{f}(X_k^\epsilon)\right)^\epsilon$ is an element of the algebra $\mathcal{A}_{2g}$.

Suppose that $\vec{a} = (a_1,\ldots, a_l)$ is a morphism in the framed seed groupoid, i.e. a sequence of mutations, framing shifts and coordinate rescalings. Let us say that such a morphism is admissible if each signed mutation in the corresponding sequence is. We define the groupoid $\mathbb{G}_{\mathrm{ad}}$ to be the subcategory of $\mathbb{G}$ whose morphisms are the admissible ones.

To each admissible morphism we may associate an invertible element $\Phi_{\vec{a}}$ of the extended algebra $\widehat{\mathcal{A}}_{2g}$. This element $\Phi_{\vec a}$ determines a birational automorphism of $\mathcal{D}_{2g}$ (by conjugation), along with an automorphism of $\mathcal{K}$ (via the representation~\eqref{eq:rep-def},~\eqref{eq:frame-change-action}.)

\begin{lemma}
\label{lem:repwelldef}
Suppose that two chains of $\vec{a}_1,~\vec{a}_2$ of admissible mutations and framing shifts induce the same birational automorphism of $\mathcal{D}_{2g}$. Then $\Phi_{\vec{a}_1} = \Phi_{\vec{a}_2}$.
\end{lemma}
\begin{proof}
The Lemma is proved by the following standard argument, {\it cf.}~\cite{KN}. If the $\Phi_{\vec{a}_i}$ induce the same birational automorphism of $\mathcal{D}_{2g}$, then the element $\Phi_{\vec{a}_1}^{-1}\Phi_{\vec{a}_2}\in\widehat{\mathcal{A}}_{2g}$ commutes with all generators $U_i,V_i$. An easy calculation shows that this implies that $\Phi_{\vec{a}_1}^{-1}\Phi_{\vec{a}_2}$ must be an element of the ground ring $\mathbb{Z}((q))$. But since each quantum dilogarithm corresponding to an admissible mutation is a formal power series in $U_i$ starting from 1, we see that $\Phi_{\vec{a}_1}^{-1}\Phi_{\vec{a}_2}=1$, and the Lemma is proved.
\end{proof}

For the purposes of understanding the integrality properties of wavefunctions, we introduce the following strengthening of the notion of admissible mutations.  Let us say that an admissible mutation in direction $e_k$ is \emph{primitive} if in the monomial
$$
\iota_\mathbf{f}(X_k) = (-q)^r\exp\left(2\pi\hbar{\sum_{j=1}^g{m_j}u_j + n_jv_j}\right), \qquad m_j,n_j,r\in\mathbb{Z}
$$
the vector 
\begin{align}
\label{eq:exp-vector}
\mathbf{m}= (m_1,\ldots,m_g) 
\end{align}
is a primitive vector in $\mathbb{Z}^g$.

%% file: cubic-graph.tex
\section{The Chromatic Lagrangian}
Fix $\G := \PGL_2$ in this section.
We begin by reviewing the constructions and results of~\cite{TZ}.

\subsection{Cubic Planar Graphs and Fukaya Moduli}
\label{sec:cpgafm}

Let $\Gamma \subset S^2$ be a cubic planar graph.  There is an integer $g$ such that $\Gamma$ has ${\sf v} = 2g+2$ vertices, ${\sf e} = 3g+3$ edges, and ${\sf f} = g+3$ faces.  As in \cite{TZ}, one may associate the following objects to $\Gamma$.

\begin{enumerate}
\item A Legendrian surface $\leg_{\Gamma} \subset T^{\infty} \bR^3 \subset S^5$ of genus $g$ \cite[Def. 2.1]{TZ}.
The surface $\leg_\Gamma$ is a branched double cover of $S^2,$ branched over the vertices of $\Gamma.$
It is defined by its front projection, which is taken to be a two-sheeted cover of $S^2$ with crossing
locus over the edges of $\Gamma$ and looking like the following near vertices:
\begin{figure}[H]
\includegraphics[scale = .25]{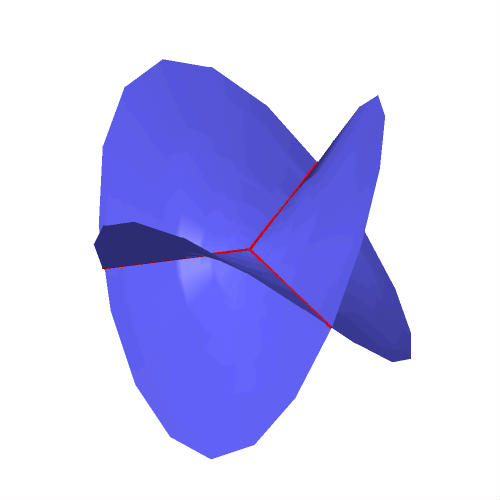}
\caption{The front projection of $\leg_\Gamma$ near a vertex.}
\label{fig:legprojs}
\end{figure}
\item A period domain $\cP_{\Gamma} \cong H^1(\leg_{\Gamma},\bC^*)$, which is an algebraic torus equipped with an algebraic symplectic form coming from the intersection pairing on $H_1(S_{\Gamma})$ \cite[\S 4.6]{TZ}.  
More precisely, let ${\bf B}$ be the set of branch points of $S_\Gamma$, corresponding to the vertices of $\Gamma$. The period domain $\cP_\Gamma$ is the moduli space that parametrizes flat line bundles over $S_\Gamma \backslash {\bf B}$ such that the monodromies surrounding the branch points are $-1.$  Note that $H^1(\leg_\Gamma, \bC^*)$ can be identified with the moduli space of flat line bundles over $S_\Gamma$. It acts on $\cP_\Gamma$ by taking the tensor product of corresponding line bundles, and this action equips $\cP_\Gamma$ with the structure of an $H^1(\leg_\Gamma, \bC^*)$ torsor.

\item \label{item:chromatic} A moduli space $\cM_{\Gamma}$ of microlocal-rank-one constructible sheaves on $\bR^3$, whose singular support lies in $S_{\Gamma}$  \cite[\S 4.3]{TZ}.  More concretely, $\cM_{\Gamma}$ is the space of $\PGL_2$-equivalence classes of $\bP^1$-colorings of the faces of $\Gamma$.
\item A Lagrangian microlocal monodromy map $\cM_{\Gamma} \to \cP_{\Gamma}$ \cite[\S 4.7]{TZ}.
It can be described as follows.  Every edge $e$ of $\Gamma$ connects branch points and therefore defines an element of $H_1(S_\Gamma).$ 
It gives rise to a character $x_e: \cP_\Gamma \to \bC^\ast$ by the canonical pairing between $H_1$ and $H^1$.  The sum of edges surrounding a face $f$ is a trivial cycle in $H_1$, so
$\prod_{e\in \partial f}x_e = 1.$  The map $\cM_\Gamma \to \cP_\Gamma$ is defined by setting $x_e$ to be the cross ratio
\begin{equation}
\label{eq:cross-ratio}
x_e = -\frac{a - b}{b-c}\cdot \frac{c-d}{d-a}
\end{equation}
where $a,b,c,d\in \bP^1 = \bC\cup \{\infty\}$ are the colors of faces surrounding an edge $e$ in the following pattern:
\begin{center}
\begin{tikzpicture}
\node at (0,.65) {$a$};
\node at (1.75,0) {$d$};
\node at (0,-.65) {$c$};
\node at (-1.75,0) {$b$};
\draw[thick] (-2,1)--(-1,0)--(1,0)--(2,1);
\draw[thick] (-2,-1)--(-1,0);
\draw[thick] (2,-1)--(1,0);
\end{tikzpicture}
\end{center}
One easily verifies the relations $\prod_{e \in \partial{f}} x_e = 1$.

We exhibit defining equations for $\cM_\Gamma.$
The characters $x_e$ generate the coordinate ring of $\cP_{\Gamma}$, obeying the relation 
\begin{equation}
\label{eq:cond-global}
\prod_e x_e = (-1)^{g+1}
\end{equation}
and further the equation
\begin{equation}
\label{eq:cond1}
x_{e_1} x_{e_2} \cdots x_{e_n}  =  1.
\end{equation}
whenever $e_1,\ldots,e_n$ label the edges of a face of $\Gamma$.  
In these coordinates, the map $\cM_\Gamma \to \cP_\Gamma$ is given parametrically by the cross ratio \eqref{eq:cross-ratio}.  But it is also given by equations, as a complete intersection, in the following way.  
Let $F$ be the set of faces of $\Gamma$.  If $e_1,\ldots,e_n$ are the edges around a face $f\in F$ taken counterclockwise, then the expression
\begin{equation}
\label{eq:face-polynomial}
V_f := 1 + x_{e_1} + x_{e_1} x_{e_2} + \cdots + x_{e_1}\cdots x_{e_{n-1}}
\end{equation}
is independent of which edge is called $e_1$.  $\cM_{\Gamma}$ is cut out by the equations $V_f = 0, f\in F$.
\newline
\end{enumerate}

Now let $\hat{\Gamma}$ denote the dual planar graph, with vertex set $V(\hat{\Gamma})$.  Since $\Gamma$ is cubic, $\hat{\Gamma}$ is a triangulation of $S^2$, and we regard its $g+3$ vertices as punctures on the sphere in the sense of Section~\ref{sec:FG-moduli-space}.  Now let $S$ be a sphere with $g+3$ punctures, and recall the corresponding moduli space $\cX_{\G,S}$ of decorated ${\rm PGL}_2$ local systems on $S$.


\begin{theorem}
\label{thm:comparison}
Let $\cP$ be the  symplectic subvariety of the cluster Poisson moduli space
${\cX_{\G,S}}$ cut out by equations~\eqref{eq:cond1}. There is a canonical algebraic Lagrangian subvariety $\cM \subset \cP$ with the following property: for every cubic planar graph $\Gamma$ with $2g+2$ vertices, there is a cluster chart $\cP_\Gamma \subset \cP$ such that the embedding $\cM_{\Gamma} \to \cP_{\Gamma}$ is isomorphic to $\cM \cap \cP_\Gamma \to \cP_\Gamma$.
\end{theorem}
\begin{proof}  
The  subvariety $\cM$ is given by the subvariety of decorated local systems whose underlying local system is trivial. We show that the intersection of $\cM$ with $\cP_\Gamma$ coincides with $\cM_\Gamma$ using the prescription for constructing a decorated local system corresponding to a point in a cluster torus described in Section 9.10 of ~\cite{FG1}. 

Let $T_{\Gamma}=(\mathbb{C}^{\ast})^{3g+3}$ be the torus parametrized by the edges of the triangulation $\widehat{\Gamma}$. As in Corollary 9.1 of {\it loc.cit.}, there is an open torus embedding
\[
\psi_{\Gamma}: (\mathbb{C}^{\ast})^{3g+3} \longrightarrow \mathscr{X}_{{\rm PGL}_2, S}.
\] 
The image of $\psi_{\Gamma}$ together with the variables parametrized by the edges of $\widehat{\Gamma}$ give rise to a cluster seed of $\mathscr{X}_{{\rm PGL}_2, S}$. 

In more detail, let $a,b,c,d \in \mathbb{P}^1$ be as in \eqref{eq:cross-ratio}. Without loss of generality, assume that $(a, b, c)=(\infty,0,1)$. The cross ratio $r(a,b,c,d)=x_e$ implies that $d=1+x_e$. The M\"obius transformation taking the triple $(a,b,c)$ to $(a,c,d)$ corresponds to 
\[
g(x_e):=\begin{pmatrix}
    x_e & 1 \\
    0 & 1 \\
\end{pmatrix}\in PGL_2.
\]
The $g+3$ punctures on the sphere correspond to the vertices of $\widehat{\Gamma}$, or equivalently, to the faces of $\Gamma$.
Let $x_{e_1}, \ldots, x_{e_n}$ be the variables associated with the edges of $\widehat{\Gamma}$ surrounding a vertex $p$. 
Under the map $\psi_\Gamma$, the holonomy surrounding the vertex $p$ is
\begin{equation}
\label{p. holonomy}
M(p):=g(x_{e_1})g(x_{e_2})\cdots g(x_{e_n})= \begin{pmatrix}
    x_{e_1}x_{e_2}\cdots x_{e_n} & V_f \\
    0 & 1 \\
\end{pmatrix},
\end{equation}
where $V_f$ is given by the expression \eqref{eq:face-polynomial}.
\begin{center}
\begin{tikzpicture}
\draw (0,0)--(-90:2.2);
\draw (0,0)--(-30:2.2);
\draw (0,0)--(30:2.2);
\draw (0,0)--(90:2.2);
\draw (0,0)--(150:2.2);
\draw (0,0)--(210:2.2);
\node at (175:1) {$\vdots$}; 
\node[blue] at (-90:1.5) {$g(x_{e_1})$};
\node[blue] at (-30:1.5) {$g(x_{e_2})$};
\node[blue] at (30:1.5) {$g(x_{e_3})$};
\node[blue] at (90:1.5) {$g(x_{e_4})$};
\node[blue] at (150:1.5) {$g(x_{e_5})$};
\node at (175:1) {$\vdots$}; 
\node[blue] at (-150:1.5) {$g(x_{e_n})$};
\draw [-latex, blue] (-120:1.2)--(-60:1.2);
\draw [-latex, blue] (-60:1.2)--(0:1.2);
\draw [-latex, blue] (0:1.2)--(60:1.2);
\draw [-latex, blue] (60:1.2)--(120:1.2);
\draw [-latex, blue] (120:1.2)--(170:1.2);
\draw [-latex, blue] (190:1.2)--(-120:1.2);
\end{tikzpicture}
\end{center}

By the expression \eqref{p. holonomy}, the conditions~\eqref{eq:cond1} generate the defining ideal of the intersection of the unipotent subvariety ${\cX}_{G, S}^{\rm un} $ with the image of $\psi_{\Gamma}$. 
 The triviality of the holonomy of the underlying unipotent local system around such a puncture is equivalent to the vanishing of the corresponding expression~\eqref{eq:face-polynomial}.

It remains to prove that $\cM$ is Lagrangian. But since any decorated local system with trivial underlying local system is a fixed point of the involution $i$ described in Section~\ref{sec:sphere-cases}, the Lagrangianicity follows from Theorem \ref{thm:general.lagrangian} --- see Remark \ref{rmk:special.case}.
\end{proof}

\begin{remark}
Because of item \eqref{item:chromatic} above, we refer to $\cM$ as the \emph{chromatic Lagrangian}.  The definition of $\cM$ from a dual perspective, along with its Lagrangian nature, were initially established
in the earlier work of Dimofte-Gabella-Goncharov \cite{DGGo}.
In this paper, our equations \eqref{eq:cond1} provide an explicit description of $\cM$ for every cubic planar graph $\Gamma$. This presentation is novel and crucial for the quantization of $\cM$, as discussed in Section \ref{sec.quantizing.la}. Meanwhile, our proof of the Lagrangian property of $\cM$ via the involution $i$ in Section \ref{sec:sphere-cases} is new and relatively simple, which can be generalized to the setting of any semisimple group.
\end{remark}

\subsection{Mutation and quantization}
\label{sec:quantization}

We define $\mathcal{T}_\Gamma := \Lambda_\Gamma \otimes_\bZ \bC^*,$ a Poisson torus.  It has a canonical
quantization $\mathcal{T}^q_\Gamma$, generated by coordinates $X_v,\; v\in \Lambda,$ with relations
\begin{equation}
\label{eq:quantum-torus-relation}
X_v X_w = q^{(v,w)}X_{v+w}.
\end{equation}

Let $\underline{\mathbf{i}}$ be a framed seed with underlying cubic graph $\Gamma$, and let $\Gamma'$ be the graph obtained from $\Gamma$ by flipping a single edge $e_0$. Then the positive and negative lattive mutation maps $\nu^\pm: \Lambda_{\Gamma'} \rightarrow \Lambda_{\Gamma}$ deliver isometries of edge lattices $\Lambda_{\Gamma'} \cong \Lambda_{\Gamma}$, and so define framed seeds $\nu_0^\pm(\underline{\mathbf{i}})$. The corresponding isometries of lattices are illustrated below:
\begin{equation}
\label{eq:torusmutation}
\begin{tikzpicture}
\draw (-1.5,-1)--(-1,0)--(-1.5,1);
\draw (-1,0)--(0,0);
\draw (.5,-1)--(0,0)--(.5,1);
\node at (-.5,1.6) {\small$\Gamma$};
\node at (-.5,.3) {\small$X_{e_0}$};
\node at (-1.6,.6) {\small$X_{e_1}$};
\node at (.7,.6) {\small$X_{e_4}$};
\node at (-1.7,-.6) {\small$X_{e_2}$};
\node at (.7,-.6){\small$X_{e_3}$};
\end{tikzpicture}
\qquad
\begin{tikzpicture}
\node at (0,0) {$\xymatrix{{}\\ {}&&\ar[ll]^{\nu^+_0}{}\\{}}$};
\end{tikzpicture}
\qquad
\begin{tikzpicture}

\node at (.15,1.6) {$ \Gamma'$};

\draw (-1,1)--(0,.5)--(1,1);
\draw (0,.5)--(0,-.5);
\draw (-1,-1)--(0,-.5)--(1,-1);
\node at (.55,0) {\small$ X_{-e_0}$};
\node at (-.8,.6) {\small$X_{e_1}$};
\node at (1.2,.6) {\small$X_{e_4+e_0}$};
\node at (-1.15,-.6) {\small$X_{e_2+e_0}$};
\node at (.8,-.6){\small$X_{e_3}$};
\end{tikzpicture}
\end{equation}
\begin{center}
\begin{tikzpicture}
\draw (-1.5,-1)--(-1,0)--(-1.5,1);
\draw (-1,0)--(0,0);
\draw (.5,-1)--(0,0)--(.5,1);
\node at (-.45,.25) {\small$X_{e_0}$};
\node at (-1.6,.6) {\small$X_{e_1}$};
\node at (.7,.6) {\small$X_{e_4}$};
\node at (-1.9,-.6) {\small$X_{e_2}$};
\node at (.7,-.6){\small$X_{e_3}$};
\end{tikzpicture}
\qquad
\begin{tikzpicture}
\node at (0,0) {$\xymatrix{{}\\ {}&&\ar[ll]^{\nu^-_0}{}\\{}}$};
\end{tikzpicture}
\qquad
\begin{tikzpicture}
\draw (-1,1)--(0,.5)--(1,1);
\draw (0,.5)--(0,-.5);
\draw (-1,-1)--(0,-.5)--(1,-1);
\node at (.55,0) {\small$X_{-e_0}$};
\node at (-1.05,.6) {\small$X_{e_1+e_0}$};
\node at (.9,.6) {\small$X_{e_4}$};
\node at (-.9,-.6) {\small$X_{e_2}$};
\node at (1,-.6){\small$X_{e_3+e_0}$};
\end{tikzpicture}
\end{center}
Also associated to each flip of triangulation is a cluster transformation, i.e.  a birational map of tori $\mathcal{T}_{\Gamma'} \dashrightarrow \mathcal{T}_{\Gamma}.$  As explained in Section \ref{sec:cluster}, these maps admit quantizations $\mathcal{T}^q_\Gamma \dashrightarrow
\mathcal{T}^q_{\Gamma'}$, which in our case take the form

\begin{center}
\begin{equation}
\label{fig:q-mut}
\begin{tikzpicture}

\draw (-1.5,-1)--(-1,0)--(-1.5,1);
\draw (-1,0)--(0,0);
\draw (.5,-1)--(0,0)--(.5,1);
\node at (-.45,.25) {\small$X_{e_0}$};
\node at (-1.6,.6) {\small$X_{e_1}$};
\node at (.7,.6) {\small$X_{e_4}$};
\node at (-1.9,-.6) {\small$X_{e_2}$};
\node at (.7,-.6){\small$X_{e_3}$};
\end{tikzpicture}
\qquad
\begin{tikzpicture}
\node at (0,0) {$\xymatrix{{}\\ {}&&\ar[ll]^{\mu_0}{}\\{}}$};
\end{tikzpicture}
\quad 
\begin{tikzpicture}
\draw (-2,1)--(-1,.5)--(0,1);
\draw (-1,.5)--(-1,-.5);
\draw (-2,-1)--(-1,-.5)--(0,-1);
\node at (-.9,1.6) {$\Gamma'$};
\node at (-.4,0) {\small$X_{-e_0}$};
\node at (-2.8,.6) {\small$X_{e_1}(1+qX_{e_0})$};
\node at (.9,.6) {\small$X_{e_4}(1+qX_{-e_0})^{-1}$};
\node at (-2.7,-.6) {\small$X_{e_2}(1+qX_{-e_0})^{-1}$};
\node at (.75,-.65){\small$X_{e_3}(1+qX_{e_0})$};
\end{tikzpicture}
\end{equation}
\qquad
\end{center}


The map $\mu$ can be factored in one of two ways, corresponding to the choice of sign in the lattice isomorphism $\nu^\pm$. Indeed, one easily verifies that the quantum cluster transformation $\mu_k$ corresponding to the flip at edge $k$ can be written as
\begin{align*}
\mu_k &= \mathrm{Ad}_{\Phi(X_{e_k})}\circ\nu^+_k\\
&= \mathrm{Ad}_{\Phi(X_{-e_k})^{-1}}\circ\nu^-_k.
\end{align*}

Now consider a morphism in the framed seed groupoid represented by a sequence of $n$ signed edge mutations $a({\mathbf{k}}):\underline{\mathbf{i}}\rightarrow\underline{\mathbf{i}}'$:
$$
\underline{\mathbf{i}}= \underline{\mathbf{i}}_0 \to^{k_1} \underline{\mathbf{i}}_1 \to^{k_2} \cdots \to^{k_n} \underline{\mathbf{i}}_n = \underline{\mathbf{i}}',
$$ 
where the $j$th mutation takes place at edge ${k_j}$ and has sign $\epsilon_j$. It gives rise to an isomorphism of quantum tori $\nu_{\mathbf{k}}:\mathcal{T}^q_{\underline{\mathbf{i}}'}\rightarrow\mathcal{T}^q_{\underline{\mathbf{i}}}$ given by
$$
\nu_{\mathbf{k}} = \nu^{\epsilon_{1}}_{k_1}\circ\cdots\circ\nu^{\epsilon_n}_{k_n}.
$$
Moreover, if we write $M_{j}$ for the image in $\mathcal{T}_{\underline{\mathbf{i}}}$ of the quantum torus element $X_{e_{k_j}}^{\epsilon_j}\in\mathcal{T}_{\underline{\mathbf{i}}_{j-1}}$ under the isomorphism
$$
\nu^{\epsilon_{1}}_{k_{1}}\circ\cdots\circ\nu^{\epsilon_{j-1}}_{k_{j-1}}~:~\mathcal{T}^q_{\underline{\mathbf{i}}_{j-1}}\rightarrow\mathcal{T}^q_{\underline{\mathbf{i}}_{0}},
$$ 
then we have
\begin{align*}
\mu_{\mathbf{k}}&:=\mu_{k_1}\circ\cdots\circ\mu_{k_n} \\
&= \mathrm{Ad}_{\Phi(M_1)^{\epsilon_1}}\circ\cdots\circ\mathrm{Ad}_{\Phi(M_n)^{\epsilon_n}}\circ\nu_{\mathbf{k}}.
\end{align*}

Per Equation \ref{eq:iso-qtori}, such a sequence of mutations of framed seeds gives rise to a birational automorphism $\mu^{\mathcal{D}}_{\mathbf{k}}:=\underline{\iota}_{\mathbf{i}}\circ\mu_{\mathbf{k}}\circ\underline{\iota}_{\mathbf{i'}}^{-1}$ of $\mathcal{D}_{2g}$,
which evidently factors as 
$$
\mu^{\mathcal{D}}_{\mathbf{k}} = \mathrm{Ad}_{\Phi(\underline{\iota}_1(M_1))^{\epsilon_1}}\circ\cdots\circ\mathrm{Ad}_{\Phi(\underline{\iota}_n(M_n))^{\epsilon_n}}\circ\underline{\iota}^*\nu_{\mathbf{k}},
$$
where we have set 
$$
\underline{\iota}^*\nu_{\mathbf{k}} := \underline{\iota}_0\circ\nu_{\mathbf{k}}\circ\underline{\iota}_n^{-1}.
$$
The reader may find it convenient to visualize the automorphism $\underline{\iota}^*\nu_{\mathbf{k}}$ as follows. Recall that the data of a framing for a seed gives rise to a decoration of the edges of its cubic graph by monomials $\mathcal{D}_{2g}$. Then the automorphism $\underline{\iota}^*\nu_{\mathbf{k}}$ is characterized by the property that it maps the monomial sitting on edge $e$ of $\Gamma'$ in framed seed $\underline{\mathbf{i}}_n$ to the monomial sitting on the corresponding edge of $\Gamma$ in framed seed $\underline{\mathbf{i}}_0$.

Now let us suppose that each signed mutation in the sequence $\mathbf{k}$ is admissible, so that under the framing isomorphism $\underline{\iota}_j$ from $\underline{\mathbf{i}}_j$ the monomial $M_j$ is mapped to an element of the algebra $\mathcal{A}_{2g}$. Then we may form the product
$$
\Phi_{a(\mathbf{k})} := \Phi(\underline{\iota}_n(M_n))^{-{\epsilon_n}}\circ\cdots\circ\Phi(\underline{\iota}_1(M_1))^{-{\epsilon_1}} \in\mathcal{A}_{2g}.
$$

Recall the representation $\mathcal{K}\simeq\mathbb{Z}((q))((X_1,\ldots, X_g))$ of the algebra $\mathcal{A}_{2g}$. The action of $\Phi_{a(\mathbf{k})}$ defines an automorphism 
$$
a(\mathbf{k}) ~:~ \mathcal{K}\rightarrow \mathcal{K} , \qquad f\longmapsto \Phi_{a(\mathbf{k})}\cdot f,
$$
and for all $A\in\mathcal{D}_{2g}$, we have the following identity of operators on $\mathcal{K}$:
\begin{align}
\label{eq:admissible-mut-intertwining}
\Phi_{a(\mathbf{k})} \circ \mu^{\mathcal{D}}_{\mathbf{k}}(A)  = \underline{\iota}^*\nu_{\mathbf{k}}(A)\circ \Phi_{a(\mathbf{k})}.
\end{align}
In particular, if {$u\in\mathcal{T}^q_{\underline{\mathbf{i}}_0}$ and $u'\in \mathcal{T}^q_{\underline{\mathbf{i}}_n}$} are related by $u = \mu_{\mathbf{k}}(u')$, then we have
\begin{align}
\label{eq:admissible-mut-intertwining-example}
\red{\Phi_a\circ \underline{\iota}_0(u)= \underline{\iota}_n(u') \circ \Phi_a}
\end{align}
as operators on $\mathcal{K}$.

The torus $\mathcal{T}_\Gamma$ associated to a cubic graph $\Gamma$, or its quantization $\cT^q_\Gamma$, is the cluster chart $\cP_\Gamma$ of $\cP,$ described in Sections \ref{sec:cpgafm}. 
In the next section we show that
the global Lagrangian submanifold $\cM \subset \cP$ is compatible with this chart-wise quantization.

\subsection{Quantizing the Chromatic Lagrangian}
\label{sec.quantizing.la}
%
%
We begin by discussing the quantization of the relevant connected component of the moduli space of framed $PGL_2$ local systems with unipotent monodromy on the punctured sphere. Fix a cubic graph $\Gamma$ of genus $g$, and as in the previous section let $\mathcal{T}^q_{\Gamma}$ be the associated quantum torus. Suppose that $e_1,\ldots, e_n$ are the edges around a face $f$ of $\Gamma$, listed in counterclockwise cyclic order around the face; note that this means that each $e_{i+1}$ precedes $e_i$ in the counterclockwise order with respect to their common vertex, so that we have
$$
X_{e_i}X_{e_{i+1}}= q^{-2}X_{e_{i+1}}X_{e_i}.
$$
Then the relation \eqref{eq:cond1}, which imposes unipotency of the monodromy around the puncture dual to the face $f$, is quantized as
\begin{equation}
\label{quantum.face.0}
X_{e_1+\ldots+e_n} = q^{-2}.
\end{equation}
Note that the relation~\ref{quantum.face.0} can be equivalently formulated as $X_{e_1}\cdots X_{e_n} = q^{-n}.$  In order to pick out the required component, let $s=\sum_{e_i\in E}{e_i} \in \widetilde{\Lambda}_\Gamma$ be the sum of the edges. We then further impose the relation that
\begin{equation}
\label{quantum.global}
X_s= (-q)^{g+3}.
\end{equation}
After quotienting by these relations, we obtain a symplectic quantum torus algebra $\mathcal{T}^q_{\underline\Gamma}$. 

We now proceed to the quantization of the additive face relations that are equivalent to the triviality of the underlying unipotent local system at a point of $\mathcal{P}$.
To this end,  set
\begin{align}
\label{quantum.face.2h}
R_{f} &= {q^{-1}}+X_{e_1} + X_{e_1+e_2}+\ldots + X_{e_1+e_2+\cdots e_{n-1}}\\
\nonumber &= q^{-1} + X_{e_1} + qX_{e_1}X_{e_2} + q^2X_{e_1}X_{e_2}X_{e_3}+\cdots + q^{n-2}X_{e_1}X_{e_2}\cdots X_{e_{n-1}}.
\end{align}
\begin{remark}
\label{rmk:rebase}
 It follows from the multiplicative face relation~\eqref{quantum.face.0} that multiplying ~\eqref{quantum.face.2h} by $qX_{e_n}$ yields
$$
X_{e_n}+X_{e_n+e_1} + X_{e_n+e_1+e_2}+\ldots + q^{-1},
$$
so we see that the ideal in the quantum torus $\mathcal{T}^q_{\Gamma}$ generated by $R_f$ is independent of our arbitrary linearization of the cyclic order on the edges around the face $f$ implicit in~\eqref{quantum.face.2h}.
\end{remark}
Let $\mathcal{I}_{\Gamma}$ be the left ideal in $\mathcal{T}_{\Gamma}^q$ generated by all~\eqref{quantum.face.0} along the global relation~\eqref{quantum.global} and the relations $R_f$ for all faces $f$. As the quantization of a Lagrangian subvariety, the D-module $\mathcal{V}_{\Gamma}:=\mathcal{D}_{\Gamma}/ \mathcal{I}_{\Gamma}$ is holonomic.

Now suppose that two regular cubic graphs $\Gamma$ and $\Gamma'$ are related by mutation at  edge $e_0$. Let us write $\mathcal{T}_{\Gamma,\Gamma'}$ for the localization of the quantum torus $\mathcal{T}_{\Gamma}$ at the Ore set $\{\prod_{k}(1+q^{2k+1}X'_{e_0})^{n_k} : k\in\mathbb{Z},n_k\geq0\}$, and write $\mathcal{T}_{\Gamma',\Gamma}$ for the analogous localization of $\mathcal{T}_{\Gamma}'$. Then the quantum mutation map $\mu_0$ in~\eqref{fig:q-mut} defines an isomorphism
$ \mu_0~:~\mathcal{T}_{\Gamma',\Gamma}\rightarrow \mathcal{T}_{\Gamma,\Gamma'}$. Let us write $\mathcal{I}_{\Gamma,\Gamma'}$ for the ideal in $\mathcal{T}_{\Gamma',\Gamma}$ generated by the quantized chromatic ideal $\mathcal{I}_\Gamma$, and $\mathcal{I}_{\Gamma',\Gamma}$ for the ideal in $\mathcal{T}_{\Gamma,\Gamma}'$ generated by $\mathcal{I}_{\Gamma'}$.

\begin{theorem}
\label{thm:d-mod} The system of quantized chromatic ideals $\{\mathcal{I}_\Gamma\}$ is compatible with quantum cluster mutations: if $\Gamma,\Gamma'$ are regular cubic graphs related by a flip at  edge $e_0$ as in Figure~\ref{fig:q-mut}, then we have $\mu_0(\mathcal{I}_{\Gamma',\Gamma}) = \mathcal{I}_{\Gamma,\Gamma'}$.
\end{theorem}

\begin{proof}
Consider the generator $R_{f,\Gamma'}$ of $\mathcal{I}_{\Gamma',\Gamma}$ associated to the left face of the graph $\Gamma'$ in Figure~\ref{fig:q-mut}, as defined in ~\eqref{quantum.face.2h}. We show that it is mapped to the corresponding 
to a generator $R_{f,\Gamma}$ of $\mathcal{I}_{\Gamma,\Gamma'}$ under $\mu_0$. As explained in Remark~\ref{rmk:rebase},
by multiplying $R_{f,\Gamma'}$ by a unit in $\mathcal{T}_{\Gamma'}$ we may assume that the edge $e_0$ at which we mutate is neither $e_1$ nor $e_{n-1}$ in the notations of~\eqref{quantum.face.2h}. Then reading counterclockwise around the left face of the right graph in Figure~\ref{fig:q-mut}, we see that
\begin{align*}
\mu_0(X'_{e_2} +qX'_{e_2}X'_{e_0} + q^2  X'_{e_2}X'_{e_0}X'_{e_1}) &= X_{e_2}(1+qX_{-e_0})^{-1} + qX_{e_2}X_{-e_0}(1+qX_{-e_0})^{-1} \\
&\phantom{=}\; +  q^2X_{e_2}X_{-e_0}(1+qX_{-e_0})^{-1}X_{e_1}(1+qX_{e_0})\\
&= X_{e_2}+ q X_{e_2}X_{e_1},
\end{align*}
where we used that $X_{e_0}X_{e_1} = q^2X_{e_1}X_{e_0}$ by the  relation~\eqref{eq:quantum-torus-relation} applied to the graph on the left of Figure~\ref{fig:q-mut}. From this computation, we see that $\mu_0(R_{f,\Gamma'}) = R_{f,\Gamma}$. The intertwining of the generators of the form~\eqref{quantum.face.0} and~\eqref{quantum.global} follows in exactly the same way.
\end{proof}

\begin{remark}
In this lengthy remark we explain the sense in which Theorem~\ref{thm:d-mod} allows us to build a global quantum Lagrangian from the compatible system of ideals in the different cluster charts.
In this context,  Theorem~\ref{thm:d-mod} implies the Lagrangian  has a well-defined `quantum structure sheaf', which is an object in the category of representations of the quantum cluster variety $\cX_{G,S}$. 
The category of such representations can be defined by means of the gluing procedure explained in~\cite{BBP}. Indeed, let us fix as in Section~\ref{subsec:quant} an initial seed with corresponding basis $\Pi_0=\{e_i\}$ for the lattice $\Lambda$. As in that section, we write $|\Pi|$ for the set of all bases for $\Lambda$ reachable from $\Pi_0$ by some sequence of sign-coherent mutations. For each $\Pi\in|\Pi|$ we take a separate copy $\mathcal{B}_{\Pi}$ of the same abelian category $\mathcal{T}^q_{\Lambda}-\mathrm{mod}$ of left modules over the quantum torus algebra associated to the lattice $\Lambda$, and form the category $\mathcal{B} = \prod_{\Pi\in|\Pi|}\mathcal{B}_{\Pi}$.

As explained in~\ref{subsec:quant}, for each $\Pi\in|\Pi|$ we have an element $\Phi_{\Pi}\in\widehat{\mathcal{T}^q}$ that depends only on the basis $\Pi$ and not the mutation sequence leading from $\Pi_0$ to it. Given a pair $\Pi_1,\Pi_2,$ we set
$$
\Phi_{\Pi_2,\Pi_1} := \Phi_{\Pi_2}\Phi_{\Pi_1}^{-1}.
$$
We use these elements to define $(\mathcal{T}^q,\mathcal{T}^q)$-bimodules $M_{\Pi_2,\Pi_1} := \mathcal{T}^q \cdot \Phi_{\Pi_2,\Pi_1} \cdot \mathcal{T}^q \subset \widehat{\mathcal{T}^q}$. 
When $\Pi_1,\Pi_2$ correspond to cubic graphs $\Gamma,\Gamma'$ differing by a single mutation, this bimodule is the one coming from the ring $\mathcal{T}_{\Gamma,\Gamma'}$  defined earlier by Ore localization.
Tensoring with any bimodule $M_{\Pi,\Pi'}$ defines an endofunctor on $\mathcal{T}^q-\mathrm{mod}$, and together they form the components of an endofunctor $\underline{\Phi}=\left(\Phi_{\Pi',\Pi}\right):\mathcal{B}\rightarrow\mathcal{B}$.
The maps $M_{\Pi_3,\Pi_1}\rightarrow M_{\Pi_3,\Pi_2}\otimes_{\mathcal{T}^q} M_{\Pi_2,\Pi_1}$
coming from the inclusions$$
 \mathcal{T}^q\cdot\Phi_{\Pi_3,\Pi_1}\cdot\mathcal{T}^q= \mathcal{T}^q\cdot\Phi_{\Pi_3,\Pi_2}\Phi_{\Pi_2,\Pi_1}\cdot\mathcal{T}^q\longrightarrow \mathcal{T}^q\cdot\Phi_{\Pi_3,\Pi_2}\cdot\mathcal{T}^q\cdot\Phi_{\Pi_2,\Pi_1}\cdot\mathcal{T}^q 
$$
define a natural transformation $\delta:\underline\Phi\rightarrow \underline\Phi\circ\underline\Phi$ making $\underline\Phi$ into a comonad. The category of representations of the quantum cluster variety can then be defined as the category $\mathcal{B}_{\underline\Phi}$ of comodules for this comonad. Objects of this category are objects $\mathcal{F}=(\mathcal{F}_\Pi)$ of $\mathcal{B}$ together with a morphism $\nabla: \mathcal{F}\rightarrow \underline{\Phi}(\mathcal F)$ satisfying $\underline{\Phi}(\nabla)\circ\nabla = \delta_{\mathcal{F}}\circ\nabla$ and $\epsilon_\mathcal{F} \circ \nabla =\mathbf{1}_{\mathcal{F}}$, where $\epsilon:\Phi\rightarrow \mathrm{Id}_\mathcal{B}$ is the counit transformation projecting to the diagonal components.

The object $(\mathcal{O},\nabla)$ where each $\mathcal{O}_\Pi = \mathcal{T}^q$ and the map $\nabla$ is defined by the inclusions 
\begin{align}
\label{eq:trivial-gluing}
\nabla_{21}:\mathcal{T}^q\rightarrow \mathcal{T}^q\cdot\Phi_{\Pi_2,\Pi_1}\subset \mathcal{T}^q\cdot\Phi_{\Pi_2,\Pi_1}\cdot \mathcal{T}^q
\end{align}
plays the role of the structure sheaf of the quantum cluster variety. Its endomorphism ring $\mathrm{End}_{\mathcal{B}_{\underline\Phi}}(\mathcal{O})$ is naturally identified with the algebra $\mathbb{L}^q_{\cX}$ of universally Laurent elements in $\mathcal{T}^q$. So we have a global sections functor 
$$
\mathcal{B}_{\underline\Phi}\rightarrow \mathbb{L}^q_{\cX}-\mathrm{mod}, \quad \mathcal{F}\rightarrow \Hom_{\mathcal{B}_{\underline\Phi}}(\mathcal{O},\mathcal{F}).
$$
In concrete terms, a global section of $\mathcal{F}=(\bigoplus\mathcal{F}_{\Pi},\nabla)$ is described by a collection of elements $f_\Pi\in\mathcal{F}_\Pi$ satisfying $\nabla_{\Pi,\Pi'}(f_\Pi)  =M_{\Pi,\Pi'}\otimes f_{\Pi'}$.
The global sections functor has a left adjoint which sends a module $V$ over $\mathbb{L}^q_{\cX}$ to its localization, i.e. the object $\mathcal{O}\otimes_{\mathbb{L}^q_{\cX}} V$ of $B_{\underline\Phi}$ with components $\mathcal{T}^q\otimes_{\mathbb{L}^q_{\cX}}V$.

Now suppose we are given a collection of left ideals $\mathcal{I}_\Pi\subset \mathcal{T}^q$ satisfying
$$
M_{\Pi'',\Pi'}\otimes_{\mathcal{T}^q}\mathcal{I}_{\Pi'} = \mathcal{I}_{\Pi''}\otimes_{\mathcal{T}^q} M_{\Pi'',\Pi'}
$$
for all $\Pi',\Pi'' \in|\Pi|$. Setting $\mathcal{F}_\Pi := \mathcal{T}^q/\mathcal{I}_\Pi$, this condition implies that the maps in ~\eqref{eq:trivial-gluing} descend to maps $\mathcal{F}_{\Pi''}\rightarrow  M_{\Pi'',\Pi'}\otimes_{\mathcal{T}^q}\mathcal{F}_{\Pi'}$, and 
we get an object $\mathcal{F} = (\mathcal{T}^q/\mathcal{I}_{\Pi'})$ of $\mathcal{B}_{\underline\Phi}$.
So in this language, Theorem~\ref{thm:d-mod} implies that the system of ideals~$ \mathcal{I}_\Gamma$ defines a representation of the quantum cluster variety $\cX_{G,S}$. In particular, taking global sections of this object defines a global chromatic left ideal $I\subset\mathbb{L}^q_\cX$.

A similar construction can be performed using the admissible framed seed groupoid $\mathbb{G}_{\mathrm{ad}}$, where we replace $\mathcal{T}^q$ by the quantum torus $\mathcal{D}_{2g}$, and the elements $\Phi_{\Pi,\Pi'}\in{\widehat{\mathcal{T}^q}}$ by the ones $\Phi_{\vec a}\in \widehat{\mathcal A}_{2g}$ associated the to arrows in $\mathbb{G}_{\mathrm{ad}}$, as defined in Section~\ref{sec:groupoid-rep}. 
\end{remark}

We now illustrate the constructions of this section in the following simple but fundamental example.
\begin{example}
\label{eg:intertwining}
Consider the framed seed $\underline{\mathbf{i}}_0$ for the $g=1$ necklace graph $\Gamma_0$ shown in Figure~\ref{fig:g1-necklace-alg}. The additive face relation
$$
R = q^{-1} + X_{e_2}
$$ 
corresponding to its left bead is mapped under the corresponding framing isomorphism $\iota_0:\mathcal{T}^q_{\Gamma_0}\rightarrow \mathcal{D}_2$ to the element
$$
\iota_0(R) = q^{-1}(1-V).
$$
Let us now perform a positive mutation at the edge $e_3$ of $\Gamma_0$ to obtain the framed seed for the canoe graph $\Gamma_1$ shown in Figure~\ref{fig:alg-canoe}. Then  we see that 
$$
R = \mu_3\left(R'\right), \qquad R' = q^{-1} + X_{e_2'} + X_{e_2'+e_3'},
$$
where $R'$ is the additive face relation associated to the face of $\Gamma'$ bounded by $(e_1',e_2',e_3')$. Under the new framing isomorphism $\underline{\iota}_1:\mathcal{T}^q_{\Gamma_1}\rightarrow \mathcal{D}_2$, the element $R'$ is mapped to
\begin{equation}
\label{eq:g1-facerel}    
\underline{\iota}_1(R')  = q^{-1} +q^{-1}UV - q^{-1}V.
\end{equation}
The element $\Phi_{\mu_3^+}\in\mathcal{A}_{2}$ is given by
\begin{align*}
\Phi_{\mu_3^+} &= \Phi\left(-q^{-1}U\right)^{-1}\\
&= (U,q^2)_{\infty},
\end{align*}
and hence the operators associated to the face relations $R,R'$ are indeed intertwined under by the action of $\Phi_a$: we have
$$
\Phi_{\mu_3^+}\circ\underline{\iota}_0(R) = \underline{\iota}_1(R') \circ \Phi_{\mu_3^+}. 
$$
\begin{figure}[htbp] 
 \begin{minipage}{0.5\linewidth} 
  \centering 
        \begin{tikzpicture}[scale=.7]


\draw[gray, thick] (3,0) circle (1.5cm);
\draw[gray, thick] (4.5,0) -- (7.5,0);
\draw[gray, thick] (9,0) circle (1.5cm);

\draw[gray, thick] (1.5, 0) .. controls  (0,5) and (12,5) .. (10.5,0);


\filldraw[red] (3,0) circle (2pt) node[anchor=west]{};
\filldraw[red] (9,0) circle (2pt) node[anchor=west]{};
\filldraw[red] (6,2) circle (2pt) node[anchor=west]{};
\filldraw[red] (6,-2) circle (2pt) node[anchor=west]{};


\node at (6,.4) {\small$-q^{-1}U$};
\node[text=blue]  at (6,-.3) {\small$3$};

\node (b) at (6,4) {\small $-qU^{-1}$};
\node[text=blue]  at (6,3.5) {\small$6$};

\node  at (3,1.8) {$-q^{-1}V^{-1}$};
\node[text=blue] at (3,1.2) {\small$1$};
\node[text=blue]  at (3,-1.2) {\small$2$};
\node (e) at (3,-1.9) {$-q^{-1}V$};

\node (d) at (9,1.9) {$-q^{-1}V$};
\node[text=blue] at (9,1.2) {\small$4$};
\node[text=blue]  at (9,-1.2) {\small$5$};
\node (e) at (9,-1.9) {$-q^{-1}V^{-1}$};

\end{tikzpicture} 
  \caption{The standard necklace framed seed $\underline{\mathbf{i}}_0$ for $g=1$} 
  \label{fig:g1-necklace-alg} 
 \end{minipage}%
 \begin{minipage}{0.5\linewidth} 
  \centering 
        \begin{tikzpicture}[scale=.75]


\draw[gray, thick] (6,1.5) -- (6,-1.5);

\draw[gray, thick] (6,1.5) -- (1.5,0);
\draw[gray, thick] (6,-1.5) -- (1.5,0);

\draw[gray, thick] (6,1.5) -- (10.5,0);
\draw[gray, thick] (6,-1.5) -- (10.5,0);


\draw[gray, thick] (1.5, 0) .. controls  (0,5) and (12,5) .. (10.5,0);


\filldraw[red] (4,0) circle (2pt) node[anchor=west]{};
\filldraw[red] (8,0) circle (2pt) node[anchor=west]{};
\filldraw[red] (6,2.5) circle (2pt) node[anchor=west]{};
\filldraw[red] (6,-2.5) circle (2pt) node[anchor=west]{};


\node  at (5.1,0) {\small $-qU^{-1}$};
\node[text=blue]  at (6.3,0) {\small$3$};

\node  at (6,4) {\small $-qU^{-1}$};
\node[text=blue]  at (6,3.5) {\small$6$};

\node  at (3,1) {\small $-q^{-1}V^{-1}$};
\node[text=blue]  at (4.5,1.3) {\small$1$};
\node  at (3,-1) {\small $q^{-1}UV$};
\node[text=blue]  at (4.5,-1.3) {\small$2$};

\node  at (9,1.1) {\small $q^{-1}UV$};

\node[text=blue]  at (7.5,1.3) {\small$4$};
\node  at (9,-1) {\small $-q^{-1}V^{-1}$};

\node[text=blue]  at (7.5,-1.3) {\small$5$};

\end{tikzpicture}
  \caption{The  framed seed $\underline{\mathbf{i}}_1=\mu_3^+(\underline{\mathbf{i}}_0)$ for the canoe graph. } 
  \label{fig:alg-canoe} 
 \end{minipage} 
\end{figure}

\end{example}

%% file: foamsphasescones.tex
\section{Foams, Phases and Framings}
\label{sec:foamsphasesandframings}

We have shown that the moduli space of constructible sheaves with singular support on a Legendrian surface $\Lambda$ is a (quantum) Lagrangian subvariety (ideal) of a symplectic leaf in a (quantized) cluster Poisson variety.  This ideal is defined by a ``wavefunction.''
The purpose of this section is to describe the combinatorics of
non-exact Lagrangian fillings $L\subset \bC^3$ of the Legendrian.  The geometric/combinatoric
set-up will allow us to make conjectures about open Gromov-Witten invariants
of the pair $(\bC^3,L).$

Here are the constructions we describe.  We begin with 
a Legendrian surface $S_\Gamma$ defined by a cubic planar graph $\Gamma\subset S^2,$
as described in previous sections.

\begin{itemize}
\setlength{\itemsep}{2pt}
\item A singular exact Lagrangian filling $L_0$ is constructed from an ideal foam, $\bF.$
\item A deformed foam $\bF'$ determines a non-exact Lagrangian filling, $L$.
\item $L$ is a branched double cover of the three-ball, branched over a tangle, also defined by $\bF'.$
\item A deformation is described by a short \emph{arc} between strands of the tangle at each vertex.
\item The map $\tau: H_1(\Lambda)\twoheadrightarrow H_1(L)$ is determined
combinatorially from $\bF'$ and the arcs.
\item A \emph{splitting} of the map $\tau$ gives a \emph{phase} and \emph{framing}.
\item We further require a maximal cone of $H_1(L).$
\item These constructions allow us to make open Gromov-Witten predictions about $(\bC^3,L).$
\item All these notions can be carried through allowed mutations of the deformed foam $\bF'.$
\end{itemize}

The upshot is that we get open Gromov-Witten predictions from the wavefunction at all
points of the framed seed groupoid accessed by allowed mutations from the necklace foam. 
This is a large class of Lagrangian fillings and framings.

We now proceed as outlined above.

\subsection{Foams}
\label{sec:foams}

A cubic graph $\Gamma$ on the sphere $S$ is dual to a triangulation of $S$.  If $\Gamma$ is three-connected, then by Steinitz's theorem it is the edge graph of a polyhedron.  A \emph{foam} $\bF$ is the dual structure to a tetrahedronization of the polyhedron:  it is a polyhedral decompsition of the three-ball $B$ with $\partial B = S.$ The data of $\bF$ includes the quadruple $(R,F,E,V)$ of regions, faces, edges and vertices.
A face or edge is called \emph{external} if it intersects the boundary,
and \emph{internal} if it does not.
The foam is \emph{ideal} if it is dual to an ideal tetrahedronization of $B$, i.e~one with no internal vertices. Even if $\Gamma$ is not dual to a polyhedron, the notion of
ideal foam makes sense --- see \cite[Definition 3.1]{TZ}.  For example, if there is a bigon between two vertices,
then there is a single edge of the foam whose boundary is those vertices --- see Example \ref{ex:foam-necklace}.

\begin{example}[Foam filling for $\necklace_g$]
\label{ex:foam-necklace}
The necklace graph has a distinguished foam filling, that we in fact believe to be unique.  This foam has no vertices:  $\bF^1$ is already smooth --- in other words there is a unique phase.  See Figure \ref{fig:necklacefoam} below.  In fact, using the local construction at the left of Figure \ref{fig:necklacefoam}, we can construct similar foam fillings of any
iterated sequence of bigon additions (handle attachments for the Legendrian surface),
starting from the genus-zero necklace (theta graph).  We refer to these as necklace-type graphs, and equip them with these canonical foam fillings.  Note that while
generic foams are dual to tetrahedronizations, these foams are
dual to somewhat degenerate tetrahedronizations.  For that reason, we will
mainly focus on foams and not their duals.

\begin{center}

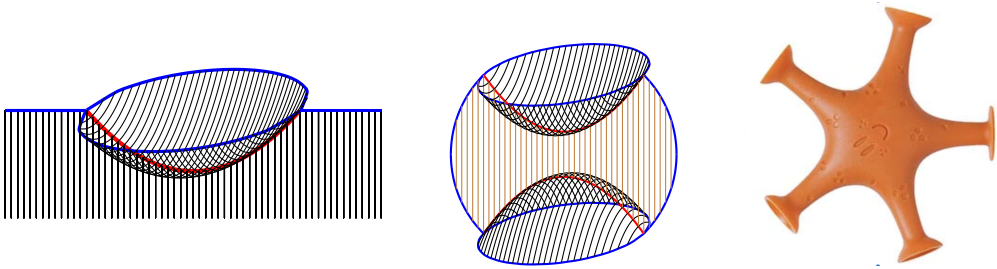
\begin{figure}[H]
\begin{tikzpicture}[scale=2.5]
\draw[white](-.1,0,0)--(.1,0,0);
 {\draw[ultra thick,red] plot[variable=\x,domain=-.5682:.5682,samples=20,smooth] (\x,.5+\x*\x,0);}
  \pgfmathsetmacro{\ymax}{sqrt(1-.3229)} 
  \pgfmathsetmacro{\radius}{\ymax-.5}
 {\draw[very thick,blue] plot[variable=\x,domain=-.5682:.5682,samples=20,smooth] (\x,\ymax,{sqrt(1-\ymax*\ymax-\x*\x)});}
 {\draw[very thick,blue] plot[variable=\x,domain=-.5682:.5682,samples=20,smooth] (\x,\ymax,{-sqrt(1-\ymax*\ymax-\x*\x)});}
 \foreach \i in {0,1,...,40}
 \pgfmathsetmacro{\x}{-.5682 + \i*(2*.5682)/40}
 \pgfmathsetmacro{\maxz}{sqrt(.3229-\x*\x)}
 \pgfmathsetmacro{\minz}{-\maxz}
 {\draw plot[variable=\z,domain=\minz:\maxz,samples=20,smooth] 
 (\x,{1/2+\x*\x+\z*\z},\z);}
 \pgfmathsetmacro{\numsamples}{60}
 \foreach \i in {0,1,...,\numsamples}
 {\pgfmathsetmacro{\x}{-1 + \i*(2*1)/\numsamples}
 \pgfmathsetmacro{\ymax}{sqrt(1-.3229)}
 \pgfmathsetmacro{\maxy}{min((1/2+\x*\x),sqrt(1-.3229))}
 \draw[opacity=.4] (\x,.25,0)--(\x,\maxy,0);}
 \pgfmathsetmacro{\ymax}{sqrt(1-.3229)} 
 \draw[thick,blue] (.5682,\ymax,0)--(1,\ymax,0);
  \pgfmathsetmacro{\ymax}{sqrt(1-.3229)} 
 \draw[very thick,blue] (-.5682,\ymax,0)--(-1,\ymax,0);
\end{tikzpicture}
\qquad
\begin{tikzpicture}[scale=1.5]
\pgfmathsetmacro{\numsamples}{30}
\pgfmathsetmacro{\b}{.2}
\pgfmathsetmacro{\ymax}{(1/2)*(-1+sqrt(1+4*(1+\b)))}
\pgfmathsetmacro{\xmax}{sqrt(1-\ymax*\ymax)}
\pgfmathsetmacro{\fraction}{(1-\xmax)/\xmax)}
\pgfmathsetmacro{\moresamps}{.5*floor(\fraction*\numsamples)}
\pgfmathsetmacro{\radius}{sqrt(1-\ymax*\ymax)}
\draw[ thick, blue] plot[variable=\t,domain=0:360,samples=20,smooth] 
 ({\radius*cos(\t)},\ymax,{\radius*sin(\t)});
\draw[ thick, blue] plot[variable=\t,domain=0:360,samples=20,smooth] 
 ({\radius*cos(\t)},-\ymax,{\radius*sin(\t)});
\draw[ thick, blue] plot[variable=\y,domain=-\ymax:\ymax,samples=20,smooth] 
 ({sqrt(1-\y*\y)},\y,0);
\draw[ thick, blue] plot[variable=\y,domain=-\ymax:\ymax,samples=20,smooth] 
 ({-sqrt(1-\y*\y)},\y,0);
\draw[ thick, red] plot[variable=\x,domain=-\xmax:\xmax,samples=20,smooth]
 (\x,{\b+\x*\x},0);
\draw[ thick, red] plot[variable=\x,domain=-\xmax:\xmax,samples=20,smooth]
 (\x,{-(\b+\x*\x)},0);
 \foreach \i in {0,1,...,\numsamples}
 {\pgfmathsetmacro{\x}{-\xmax + \i*(2*\xmax)/\numsamples}
 \pgfmathsetmacro{\maxz}{sqrt(1-\ymax*\ymax-\x*\x)}
 \pgfmathsetmacro{\minz}{-\maxz}
 \draw plot[variable=\z,domain=\minz:\maxz,samples=\numsamples,smooth] 
 (\x,{\b+\x*\x+\z*\z},\z);
 \draw plot[variable=\z,domain=\minz:\maxz,samples=\numsamples,smooth] 
 (\x,{-(\b+\x*\x+\z*\z)},\z);
 \draw[brown] (\x,{-(\b+\x*\x)},0)--(\x,{(\b+\x*\x)},0); 
 }
 \foreach \i in {1,...,{\moresamps}}
  {
 \pgfmathsetmacro{\stupidtikz}{(1-\xmax)*\i}
 \pgfmathsetmacro{\x}{\xmax+\stupidtikz/\moresamps}
 \draw[brown] (\x,{-sqrt(1-\x*\x)},0)--(\x,{sqrt(1-\x*\x)},0);
 \draw[brown] (-\x,{-sqrt(1-\x*\x)},0)--(-\x,{sqrt(1-\x*\x)},0);}
\end{tikzpicture}
\qquad 
\includegraphics[scale=.2]{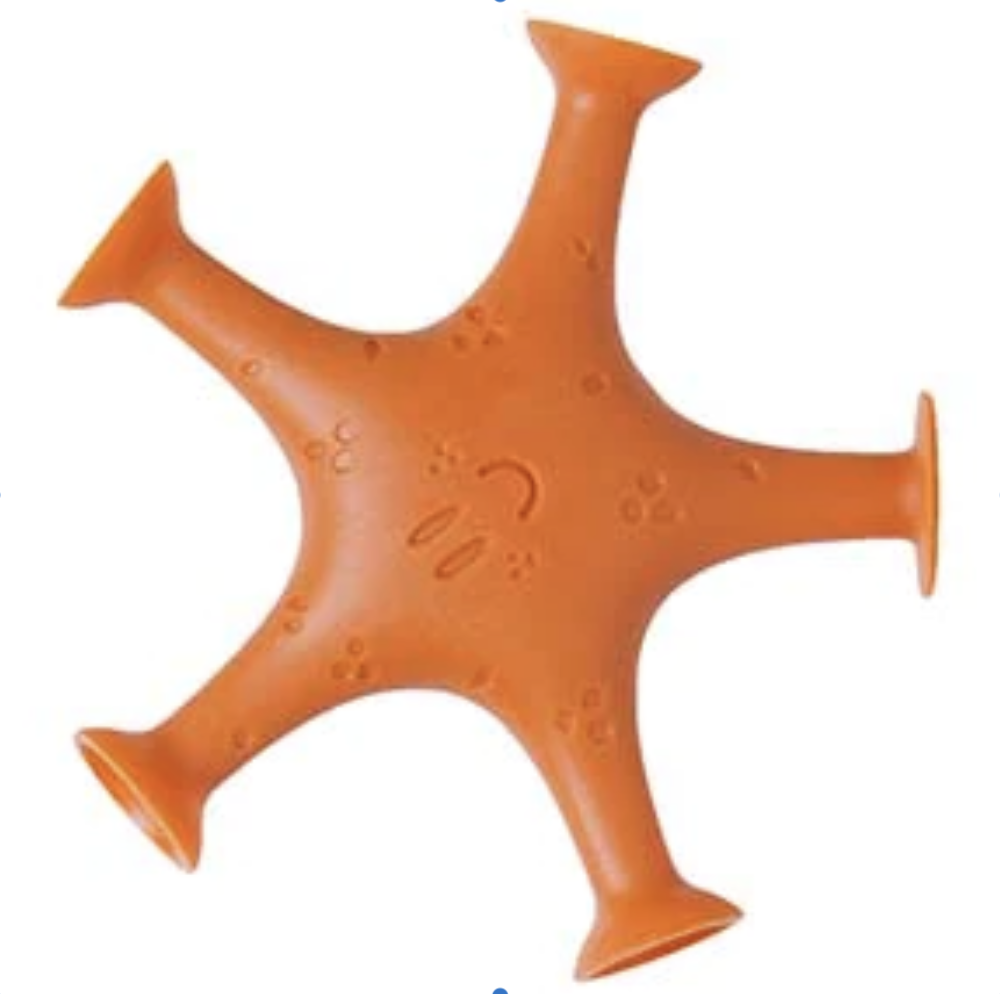}
\caption{The necklace graph ${\necklace_g},$ pictured in blue, and its foam filling.
At left is a local model near a bead, with tangle strand in red. 
In the middle is the foam filling for $g=1$.
The Ubbi toy at right is something close to a foam for ${\necklace_{4}}$.}
\label{fig:necklacefoam}
\end{figure}
\end{center}
\end{example}

\subsubsection{The Harvey-Lawson Foam}
\label{subsec:hlfoam}

The foam $\bF_{\mathrm{HL}}$ of the Harvey-Lawson Lagrangian $L_0\subset \bC^3$
has a single vertex at the origin in $\bR^3,$
four edges $E_i$ equal to the rays $\bR_{\geq 0}\cdot v_i$
where $v_0 = (1,1,1)$ and $v_i = -v_0 + 2e_i$, with $e_i$ the standard basis vectors.
There are $6 = \binom{4}{2}$ faces $F_{ij}$
equal to the cones spanned by unordered pairs of edges,
and $4 = \binom{4}{3}$ regions equal to the cones spanned by triples of edges. 
(It can be succinctly described as the toric fan of ${\bf P}^3.$)

The singular, exact Harvey-Lawson Lagrangian $L_0$ in
$(\bC^3,\omega_{\mathrm{std}}=d\theta_{\mathrm{std}})$
is a branched $2:1$ cover of $\bR^3,$ branched over the over edges.
$L_0$ is a cone over $S^1\times S^1$ with parametrization
$(r,s,t) \mapsto (re^{is},re^{it},re^{-i(s+t)}) \in \bC^3$,
where $r\in \bR_{{\geq 0}}$ and $(s,t)\in S^1\times S^1.$
The covering map is the restriction to $L_0$ of $\bC^3\to \bR^3$
sending a complex triple to its real part:
explicitly $(r,s,t) \mapsto (r\cos(s),r\cos(t),r\cos(s+t)) \in \bC^3.$
The map is $1:1$ over the four rays with $(s,t) = (0,0), (\pi,0), (0,\pi),(\pi,\pi),$
which we think of as a singular tangle.
The six sheets of the foam are defined by $s = 0, s = \pi, t = 0, t = \pi, s+t = 0, s+t = \pi.$

The primitive function $f$ obeying $df = \theta_{\mathrm{std}}\vert_L$ is
$f = \frac{1}{4}r^2\left(\sin(2s)+\sin(2t)-\sin(2s+2t)\right).$  Note that $f$ is
odd under the hyperelliptic-type involution $(s,t)\leftrightarrow (-s,-t)$ and
$f=0$ along the preimages of the sheets of the foam.  Thus $f$ allows us to label the
branches of $L_0$ on the regions $R$.

\begin{center}
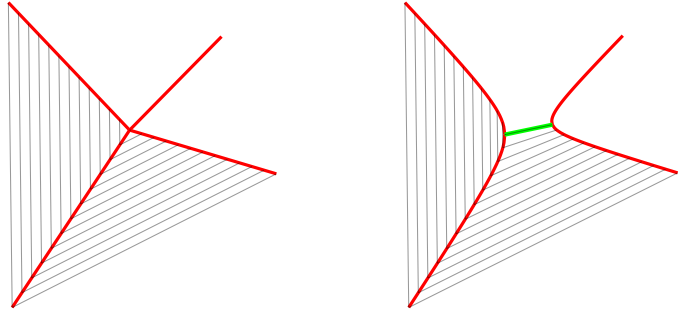
\begin{figure}[H]
\begin{tikzpicture}[rotate around y=15, rotate around z=7]
\pgfmathsetmacro{\numsamples}{12}
\pgfmathsetmacro{\s}{1.5}
\pgfmathsetmacro{\fudge}{1}
\coordinate (O) at (0,0,0);
\coordinate (T) at (\s,\s,\s);
\coordinate (A) at (-\s,-\s,\fudge*\s);
\coordinate (B) at (-\s,\fudge*\s,-\s);
\coordinate (C) at (\fudge*\s,-\s,-\s);
\draw[very thick, red] (A)--(O)--(T);
\draw[very thick, red] (B)--(O)--(C);
\foreach \i in {0,1,...,\numsamples}
{\pgfmathsetmacro{\rval}{\i/\numsamples}
\draw[opacity=.4] (-\rval*\s,-\rval*\s,\fudge*\rval*\s)--(\fudge*\rval*\s,-\rval*\s,-\rval*\s);
\draw[opacity=.4] (-\rval*\s,-\rval*\s,\fudge*\rval*\s)--(-\rval*\s,\fudge*\rval*\s,-\rval*\s);}

\end{tikzpicture}
\qquad\qquad
\begin{tikzpicture}[rotate around y=15, rotate around z=7]
\pgfmathsetmacro{\numsamples}{12} 
\pgfmathsetmacro{\lll}{1/sqrt(1/3)}
\pgfmathsetmacro{\eps}{.3}
\pgfmathsetmacro{\range}{1.5}
\draw[very thick, red] plot[variable=\r,domain=0:\range,samples=20,smooth]
({sqrt(\r*\r+\eps*\eps)},\r,\r);
\draw[very thick, red] plot[variable=\r,domain=0:\range,samples=20,smooth]
({sqrt(\r*\r+\eps*\eps)},-\r,-\r);
\draw[very thick, red] plot[variable=\r,domain=0:\range,samples=20,smooth]
({-sqrt(\r*\r+\eps*\eps)},\r,-\r);
\draw[very thick, red] plot[variable=\r,domain=0:\range,samples=20,smooth]
({-sqrt(\r*\r+\eps*\eps)},-\r,\r);
\draw[ultra thick, green] (-\eps,0,0)--(\eps,0,0);
\foreach \i in {0,1,...,\numsamples}
{\pgfmathsetmacro{\rval}{\i*\range/\numsamples}
\draw[opacity=.4] ({sqrt(\rval*\rval+\eps*\eps)},-\rval,-\rval)--({-sqrt(\rval*\rval+\eps*\eps)},-\rval,\rval);
\draw[opacity=.4] ({-sqrt(\rval*\rval+\eps*\eps)},\rval,-\rval)--({-sqrt(\rval*\rval+\eps*\eps)},-\rval,\rval);}
\end{tikzpicture}
\caption{Left:  the Harvey-Lawson foam, with its four edges but just two of the $\binom{4}{2}=6$ faces drawn.  Right:  the deformed foam, with arc in green, and the two deformed faces drawn.}
\label{fig:hlfoam}
\end{figure}
\end{center}

\begin{center}
\begin{figure}[ht]
\includegraphics[scale=.14]{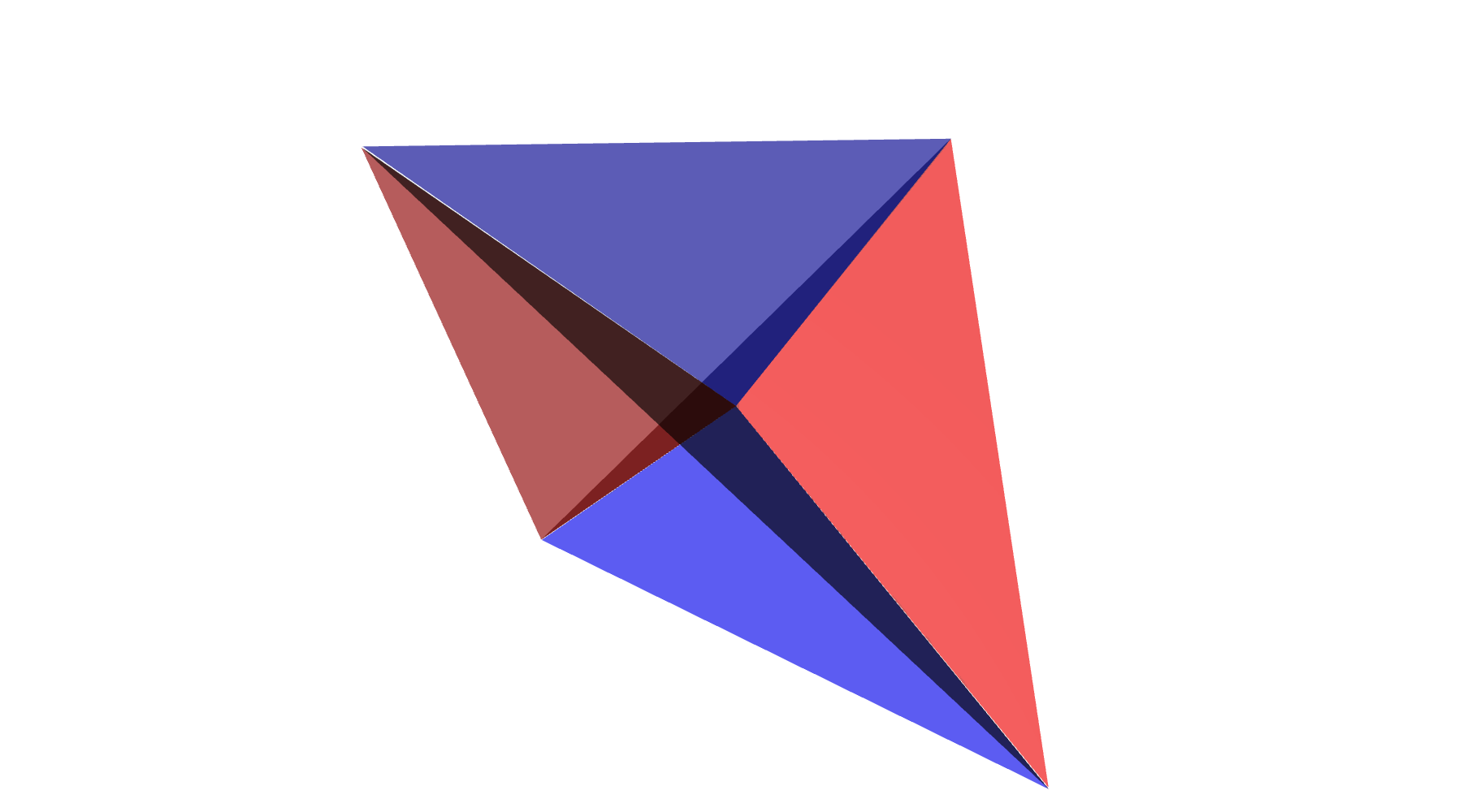}\hskip-.9in
\includegraphics[scale=.15]{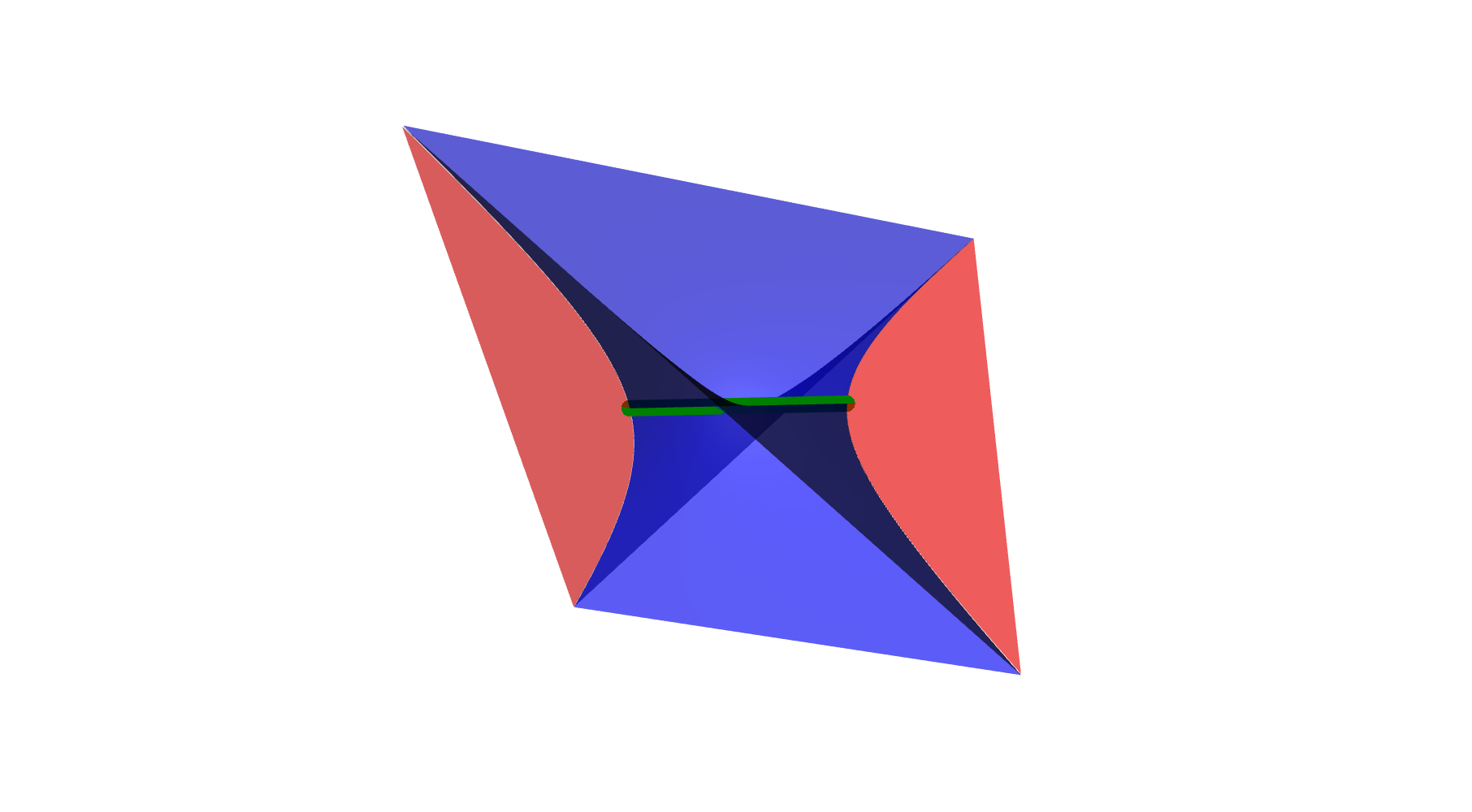}
\caption{At left is the foam of the Harvey-Lawson Lagrangian, with all $6 = \binom{4}{2}$ sheets drawn.  Warning:  diagonally opposite vertices lie in different half-spaces, so despite appearances the two corresponding triangular sheets do \emph{not} intersect outside the origin.  At right is the deformed
foam of its smoothing.  Two sheets (pink) are smoothed to have hyperbolas as boundaries, while the boundaries of the other four consist of two halves of different hyperbolas, as well as the arc (green).}
\end{figure}
\end{center}

\subsubsection{Foams and singular exact Lagrangians}

From a foam $\bF$ we can define a singular exact (not necessarily special) Lagrangian $L_0$ locally modeled on the Harvey-Lawson cone and foam --- see \cite[Section 3.2]{TZ}.  As with the Harvey-Lawson cone
and foam, we can define a multi-valued function $f$ whose sign labels the branches of $L_0$
in the regions of the foam.

\subsubsection{Deformation of the Harvey-Lawson foam}
\label{sec:hlfoamdef}

There are three distinct families of smoothings of $L_0$ corresponding
to the three matchings of the four
edges. 
We will describe the one for the matching $0\leftrightarrow 1, 2\leftrightarrow 3$;
the others are similar and are related by a permutation of
coordinates.
The smoothing $L_\epsilon$
has the topology of $\bR^2 \times S^1$ and has a parametrization in polar
coordinates 
$(r,s,t) \mapsto (\sqrt{r^2+\epsilon^2}e^{is},re^{it},re^{-is-it}) \in \bC^3$, which maps to 
$(\sqrt{r^2+\epsilon^2}\cos(s),r\cos(t),r\cos(s+t))\in \bR^3.$  
These are all diffeomorphic for $\epsilon \neq 0,$ so when we are interested in topological
questions, we can restrict to $L_1$ without loss of generality.
The branched cover is $1:1$ over the points with
$(s,t) = (0,0), (\pi,0), (0,\pi),(\pi,\pi),$ and these parametrize four rays $E_i$ which
trace out \emph{two} hyperbola components ($E_0\cup E_1 = \{x^2-y^2 = 1, y=z, x> 0\}$
and $E_2\cup E_3 = \{x^2-y^2 = 1, y=-z,
x < 0\}$), a smoothing of the singular
tangle of $L_0$. )

There is also the line segment $a \subset \bR^3$
between $(-1,0,0)$ and $(1,0,0)$ which we call an \emph{arc} --- it
is the image of $r=0.$  Note that $L_1\to \bR^3$ is $2:1$ over the arc.
The six sheets $F_{ij}$ now bound either a smooth edge $E_i\cup E_j$ if $(ij) = (01)$ or $(23)$,
or otherwise the union $E_i \cup E_j \cup a$.  This will be our local model of a deformed foam.
More generally,
let $s_i$ be the matching of edges of ${\bf F}_{\mathrm{HL}}$
which pairs $v_0$ and $v_i$.  We write ${\bf F}_{\mathrm{HL},s_i}$ for the local deformed foam of $L_1$
Its arc is the line segment between $-e_i$ and $e_i$.
We write ${\bf F}_{\mathrm{HL},\epsilon,s_i}$ for the deformed foam of $L_\epsilon.$

Away from the origin and the preimage of the arc, the Harvey-Lawson cone and its smoothing are homeomorphic:  $L_0\vert_{r\neq 0}\cong L_1\vert_{r\neq 0}$.  As a result, the same primitive function $f$ can be used to label regions of the foam and of its deformation, at least away from the
arc.  The local geometry of $L$ and the deformed foam near an arc is shown in Figure
\ref{fig:foam-upclose}.

\begin{center}
\begin{figure}[ht]
\begin{tikzpicture}
 \pgfmathsetmacro{\a}{5}
 \pgfmathsetmacro{\b}{-3}
  \pgfmathsetmacro{\c}{.707}
  \pgfmathsetmacro{\d}{.4}
    \pgfmathsetmacro{\s}{.3}
    \pgfmathsetmacro{\lift}{1.2}
    \pgfmathsetmacro{\n}{4}
    \pgfmathsetmacro{\nn}{2*\n}
    \pgfmathsetmacro{\ninv}{\n^(-1)}

\draw[fill,gray,opacity=.5] (-\a+\c*\d,\b-\d)--(\a+\c*\d,\b+\d)--(\a-\c*\d,\b-\d)--(-\a-\c*\d,\b+\d);
\draw[fill,gray,opacity=.5] (-\a+\c*\d,\b-\d)--(\a-\c*\d,\b-\d)--(\a+\c*\d,\b+\d)--(-\a-\c*\d,\b+\d);
\draw[green,very thick] (-\a,\b)--(\a,\b);
\draw[green,very thick] (0,0) ellipse (5 and 1);
\draw[fill,red] (\a,0) circle (.1);
\draw[fill,red] (-\a,0) circle (.1);
\draw[fill,red] (\a,\b) circle (.1);
\draw[fill,red] (-\a,\b) circle (.1);
\foreach \ii in {0,...,1}
{
\pgfmathsetmacro{\jump}{3.6*\ii}
\foreach \i in {0,...,\nn}
{
\pgfmathsetmacro{\k}{(\i - \n)*\ninv}
\pgfmathsetmacro{\kk}{\k*\k}
\pgfmathsetmacro{\testy}{0}
\pgfmathsetmacro{\j}{\a*\k}
\draw[fill,green] (\j,\b+\lift+\jump) circle (.05cm) ;
\draw[thick,black] (\j - \s*1 ,\b+\lift+\jump - \s*\k)--(\j  + \s*1,\b+\lift + \jump + \s*\k);  
\draw[thick,black] (\j - \s*\k ,\b+\lift+\jump - \s*1)--(\j  + \s*\k,\b+\lift + \jump + \s*1);  
\pgfmathsetmacro{\z}{(2*\ii-1)*(.17)*(1-\kk)^(1/2)}
\draw[fill,violet] (\j ,\b + \lift + \jump - \z) circle (.05cm);
}}
\foreach \i in {0,...,\nn}
{\pgfmathsetmacro{\k}{(\i - \n)/\n}
\pgfmathsetmacro{\j}{\a*\k}
\pgfmathsetmacro{\jump}{3.6}
\node at (\j + .35,\b + \lift + \jump + .35) {$+$};
\node at (\j - .35,\b + \lift + \jump - .35) {$+$};
}
\foreach \i in {0,...,\nn}
{
\pgfmathsetmacro{\k}{(\i - \n)/\n}
\pgfmathsetmacro{\j}{\a*\k}
\pgfmathsetmacro{\jump}{0}
\node at (\j + .35,\b + \lift + \jump + .35) {$-$};
\node at (\j - .35,\b + \lift + \jump - .35) {$-$};
}

\end{tikzpicture}
\caption{The neighborhood of an arc (the green line segment) and its lift to the Lagrangian (green oval).  Four sheets, forming two surfaces (gray) meet at the arc.  The cross-sectional planes are shown, along with the sign of the primitive function $f$ on $L$.  The red dots are where the sheets
$s = 0,\pi$ meet the arc, so the sign of $f$ changes as they are crossed.  The purple dot
is the cross section of the oriented loop $\gamma_a$ --- see Definition \ref{def:ordloop}.}
\label{fig:foam-upclose}
\end{figure}
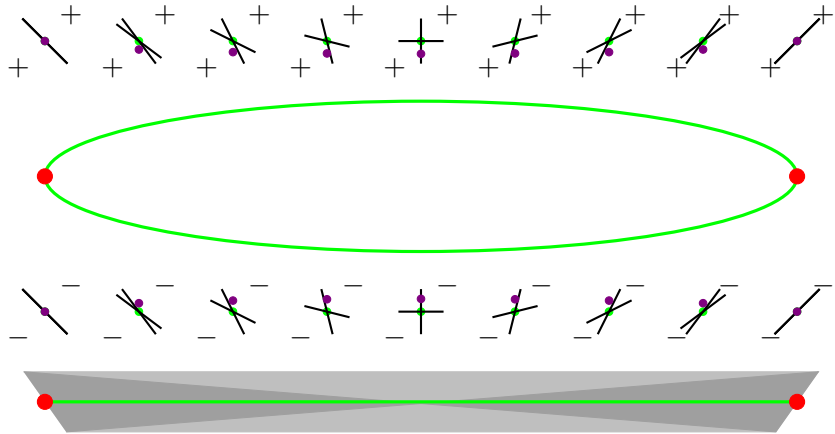
\end{center}

\subsubsection{Deformed Foams}

Given a foam $\bF$, we will define a deformed foam $\bF'$ to be a structure locally modeled near each
vertex on a Harvey-Lawson deformed foam.

\begin{definition}
Let ${\bf F}$ be a foam with vertex set $V$ consisting of $n:= \#V$ vertices.
Write $\cS$ for the set of matchings of half edges at each vertex, so $\#\cS = 3^n.$
Let $s\in \cS$.  {We define a \emph{deformed foam} ${\bf F}_s$ to be any set of vertices,
edges, arcs, faces and regions which agrees with ${\bF}$ outside some $3\epsilon$-neighborhood of $V$,
is homeomorphic to ${\bf F}$ outside of a $2\epsilon$-neighborhood of $V$,
and which is linearly equivalent to the local deformed foam
of ${\bf F}_{\mathrm{HL},\epsilon,s}$ of Section \ref{sec:hlfoamdef} within a $2\epsilon$-ball
of each vertex.}
\end{definition}

The smoothed Lagrangian $L$ is branched over a tangle $T\subset B,$ i.e.~$T$ is a one-manifold. 
The construction of $L$ from a deformed foam identities a particular
set of branch cuts we call \emph{arcs}.

\begin{definition}
Let $s_i$ be the smoothing of the Harvey-Lawson foam
which matches the ray $\bR_{\geq 0}\cdot v_0$
with $\bR_{\geq 0}\cdot v_i$, where $v_0$ and $v_i$ are as in Section \ref{subsec:hlfoam}. 
The \emph{arc} of the deformed Harvey-Lawson foam $\bF_{\mathrm{HL},s_i}$
is the line segment from
$-e_i$ to $e_i$.
An \emph{arc} of a deformed foam $\bF'$ is the locus in $B$ corresponding to the
arc of the Harvey-Lawson foam under the local identification of $\bF'$ with
$\bF_{\mathrm{HL},s_i}$.
We write $A$ for the set of arcs of a deformed foam $\bF'$.
\end{definition}

Given an arc $a$ of a deformed foam $\bF'$ and its associated branched double cover $\pi: L\to B,$
define $\lambda_a = \pi^{-1}(a).$  Since $\pi$ is $2:1$ over the interior of the arc and $1:1$ at its edges, $\lambda_a$ is a circle.  Note that $\lambda_a$ does not (yet) have a distinguished orientation.

\begin{center}
\begin{tikzpicture}
 \pgfmathsetmacro{\a}{2}
 \pgfmathsetmacro{\b}{-1.7}
\draw[green,very thick] (-\a,\b)--(\a,\b);
\draw[green,very thick] (0,0) ellipse (2 and .5);
\draw[fill,red] (\a,0) circle (.1);
\draw[fill,red] (-\a,0) circle (.1);
\draw[fill,red] (\a,\b) circle (.1);
\draw[fill,red] (-\a,\b) circle (.1);
\draw[thick,-stealth] (0,-.8)--(0,-1.4);
\node at (.3,-1.1) {$\pi$};
\node at (.5,-1.9) {$a$};
\node at (1,-.7) {$\lambda_a$};
\end{tikzpicture}
\end{center}

\begin{remark}
\label{rmk:neckfoam}
If $\bF$ is a foam, a deformed foam $\bF'$ is defined by choosing a matching of the four internal edges meeting at each vertex.  In the case of the necklace graph $\Gamma^g_{\mathrm{neck}}$,
the foam $\bF$ has no internal vertices,
and therefore $\bF' = \bF$.
In particular, $\bF$ is already deformed.  In fact,
we will learn that the $g+1$ strands of $\Gamma^g_{\mathrm{neck}}$
can be thought of as the
arcs of $\bF' = \bF$, and the face of $\bF'$ which
they bound gives rise to a single relation among them -- see Definition
\ref{def:facerelation} and Proposition \ref{prop:hladdition}.
Similar
considerations apply whenever $\Gamma$ has a bigon.
\end{remark}

\subsection{Phases and Framings}

\begin{definition}
\label{def:phaseandframing}
Let $H\cong \bZ^{2g}$ be a rank-$2g$ lattice with a non-degenerate, antisymmetric pairing $\omega.$ 
A \emph{phase} is a rank-$g$ isotropic subgroup $K\subset H$.  A \emph{framing} of $K$
is a transverse isotropic subspace.
We call the combination of phase and framing an \emph{isotropic splitting}, or sometimes
just \emph{splitting}.
\end{definition}

We will be studying phases and framings when $H = H_1(\Lambda)$ is the homology
of a genus-$g$ surface, $\Lambda$ and $\omega$ is the intersection pairing.
So let $L$ be an orientable three-manifold with boundary a genus-$g$ surface $\Lambda = \partial L$, with  $H_2(L)=0$. Then it follows from the long exact sequence in homology together with the Poincar\'e-Lefschetz duality isomorphisms $H^1(L)\simeq H_2(L,\Lambda), H_1(L,\Lambda)\simeq H^2(L)$ that $b_1(L)=g,$\footnote{In this section,
homology will be taken with $\bZ$ coefficients unless otherwise stated.}
so that we obtain the short exact sequence
\begin{equation}
\label{eq:ses}
\xymatrix{0 \ar[r]& H^1(L) \ar[r]& H_1(\Lambda) \ar[r]^{\tau} & H_1(L) \ar[r] & 0} 
\end{equation}

The notion of phases and framings will apply to above geometic setting.

\begin{definition}
\label{def:geometricphase}
Suppose $H = H_1(\Lambda)$ is the first homology of a genus-$g$ oriented surface $\Lambda$,
and $L$ is an orientable three-fold with $\partial L = \Lambda$. 
We have the short exact sequence of Equation \ref{eq:ses}.
We say that a phase $K\subset H$ is \emph{geometric} if $K = Ker(\tau)\cong H^1(L).$
An accompanying framing is an $\omega$-isotropic splitting
$\tau: H_1(L)\hookrightarrow H_1(\Lambda)$ of the short exact sequence \eqref{eq:ses}.
\end{definition}

In the context of open Gromov-Witten theory, $\Lambda$ is
Legendrian in a contact manifold and $L$
is Lagrangian in a symplectic filling.
%
%
%
%

\begin{remark} In \cite{TZ}, the above geometric phases
were called ``OGW framings" to
connote open Gromov-Witten theory.  The definition was generalized from \cite{AKV},
where mirror symmetry was used to make
conjectures in open Gromov-Witten theory.\footnote{In \cite{I} a similar definition of framing is made, but without the isotropic condition.}  The terminology
stems from the connection to Chern-Simons theory through large-N duality,
where Lagrangians are knot conormals and framing relates to the framing of knots.
We describe the connection to open Gromov-Witten theory later in this section.  
\end{remark}

\begin{remark}
We need phases and framings to define a framed seed as in Definition \ref{def:framedseed}, from which we will construct wavefunctions
and conjectural enumerative inormation --- see Conjecture \ref{conj:ogw}.  The geometry behind this is as follows. 
Let $\cX = \cX_{PGL_2,S^2}$ be the cluster variety of framed local systems on a sphere.  Let $\cP$ be the symplectic leaf of unipotents, and let $\cM$ be the Lagrangian subvariety defined by trivial monodromy.
Let $\Gamma$ be a cubic graph on the sphere $S = S^2$ and $S_\Gamma \subset J^1(S)$ the associated a Legendrian surface up to isotopy. 
We write $\cP_\Gamma = H^1(S_\Gamma;\bC^*)$ for the corresponding cluster chart.\footnote{As explained in Section \ref{sec:cpgafm}, the cluster charts $\cP_\Gamma$ of $\cP$ are spaces of rank-one local systems with fixed
monodromy $-1$ {around the critical points of
the branched double cover} $S_\Gamma \to S^2.$ 
Since $\cP_\Gamma$ is a torsor over $H^1(S_\Gamma;\mathbb{C}^*)$, its tangent space  at any point is canonically $H^1(S_\Gamma;\bC)$ and its Poisson structure is determined by the intersection
form on $S_\Gamma,$ independent of choice of base point. Hereafter, we often
omit the distinction and refer to cluster charts as the tori $H^1(S_\Gamma;\bC^*).$}
A splitting allows us to write $H^1(S_\Gamma;\bC^*)$ as
$T^*(H^1(L;\bC^*))/H_1(L)$.
When we lift $\cM$ to $T^*H^1(L;\bC^*)$ we can write it locally as the
graph of the differential of a function $W_\Gamma$ on $H^1(L;\bC^*)$, from which we will
extract enumerative information --- see Section \ref{sec:classical-limit}.
\end{remark}

Recall from \cite{TZ} the combinatorial model of the first homology of a Legendrian
$\Lambda:=S_\Gamma$
defined from a cubic planar graph $\Gamma$ on a sphere, $S$.
The faces of $\Gamma$ define a relation $\sim_\Gamma$ on the edge
lattice $\bZ^{E_\Gamma}$, namely
$\sum_{e\in \partial f}e \sim_\Gamma 0.$
We then have
$H_1(S_\Gamma) \cong \bZ^{E_\Gamma}/\!\sim_\Gamma$.
We have an antisymmetric pairing $\overline{\omega}$ on $\bZ^{E_\Gamma}$,
depending only on the orientation of $S$, defined by $\overline{\omega}(e,e') = \pm 1$ if $e$ and $e'$ are adjacent to a vertex $v$ with $e$ preceding/following $e'$ in the cyclic ordering at the vertex
\begin{tikzpicture}[scale=.2]
\draw (-1.4,1)--(0,0)--(1.4,1);
\draw (0,0)--(0,-1.5);
\node at (-2.2,1) {$e'$};
\node at (.9,-.5) {$v$};
\draw[fill] (0,0) circle (.3cm) ;
\node at (2.2,1) {$e$};
\end{tikzpicture}, and zero otherwise.  
Since $\sum_{e\in \partial f}e$ generates the kernel of this pairing,
$\overline{\omega}$ descends to a nondegenerate, antisymmetric intersection pairing $\omega$
on 
\begin{equation}
\label{eq:h1}
H_1(S_\Gamma) \cong \bZ^{E_\Gamma}/\!\sim_\Gamma.
\end{equation}

We now have a combinatorial model of $H_1(S_\Gamma).$  We next
build combinatorial models of $H_1(L)$ for $L$ arising from a deformed
foam, and of the map $H_1(S_\Gamma)\to H_1(L).$

\subsection{Combinatorics of Tangles from Deformed Foams}

We continue our study of smooth Lagrangians arising from deformed foams.
Cutting to the chase, the loops defined by the edge set $E$ and arc set $A$
will generate $H_1(S_\Gamma)$ and $H_1(L)$,
with relations determined by faces.  In total, we find 
$$\xymatrix{\bZ^{E_\Gamma} \ar[r]^{\sim_\Gamma \quad}\ar[d]^{\iota}& H_1(S_\Gamma)\ar[d]^{\tau}\\
\bZ^{E_{\Gamma}\cup A}\ar[r]^{\sim_{\bf F'}}&H_1(L)}.$$
Here $\iota$ is induced by the inclusion $E_\Gamma\to E_\Gamma\cup A,$ and
the top line was defined in the previous section.  The bottom line will
be defined in this section.

Let $\Gamma \subset S$ be a cubic graph on the sphere, let $\bF$ be an ideal foam on the three-ball $B$, whose regions, faces and edges respectively bound the faces, edges and vertices of $\Gamma$.
Let $\cS$ be the discrete set of smoothings of $\bF$, i.e.~the set of matchings of edges incident to each vertex of $\bF$ --- so $\#\cS = 3^{\#V}.$  Let $s$ be a smoothing whose resulting tangle $T$
has no circle components.  Let $L$ be a smooth Lagrangian corresponding to the
deformed foam $\bF_s$.

Recall that for an arc $a$ we werite $\lambda_a := \pi^{-1}(a).$  We now define
an orientation on $\lambda_a$, thus definining an element $\gamma_a\in H_1(L).$

Since the construction of the smoothing is local, we need only look at the Harvey-Lawson
smoothing $L_1$ and its unique arc $a$, which we can lift to the \emph{parametrized}
curve $(e^{it},0,0)$ and take the induced orientation.  This is the
orientation induced from the unique holomorphic disk in $\bC^3$ bounding $L$,
i.e.~$|z_1|\leq 1.$  We can also
give a more combinatorial construction that does not require an explicit local
model, as follows.

\begin{definition}
\label{def:ordloop}
We choose a canonical orientation for $\lambda_a$ by
orienting the arc arbitrarily and taking
a push-off $\tilde{\lambda}_a$ of the path along the
arc that has some combinatorial properties, using
the primitive function, $f$.  We require that near the start of the push-off,
in the chosen
orientation, that $f$ has a negative value and lies in one of
the two fat regions (see Figure \ref{fig:foam-upclose}) ---
in particular, outside of two sheets which meet at the arc's origin ---
then crosses once at the midpoint of the arc in a counterclockwise direction
(in the induced orientation of the transverse plane).  The
remainder of $\tilde{\lambda}_a$ traverses the arc backwards after
crossing the origin of the transverse plane at the arc's endpoint, and has
the same combinatorial recipe as the first half of $\tilde{\lambda}_a.$  This completes the
description of the push-off, $\tilde{\lambda}_a.$
There are actually two such push-offs, but the
resulting paths are homotopic. 
Likewise, the opposite orientation of the arc leads to a homotopic path (just shifted). 
For an arc $a$, write $\gamma_a$ for the resulting element of $H_1(L).$
\end{definition}

\begin{remark}
The prescription in Definition \ref{def:ordloop} is similar to how
unoriented edges of a cubic planar graph
lead to an oriented loop of the associated branched cover --- see \cite[Section 4.6]{TZ} ---
albeit somewhat more intricate.
\end{remark}
%
%
%

For each face of the deformed foam ${\bf F'},$ we will
define a relation among the edges $e$
and arc loops $\gamma_a$ along its boundary.
Together these relations will characterize $H_1(L)$ as $\bZ^{E_\Gamma\cup A}/\sim_{\bf F'}.$
%
%
%
To define the relation, we need a careful discussion of the sign of an arc relative to a face. 

\begin{definition}
\label{def:arc-sign}
Let $F$ be a face of a deformed foam $\bF'$ bounding an arc $a\in A$.
Then we have a homeomorphism of a
neighborhood of $a$ with a neighborhood of the lone arc of the Harvey-Lawson deformed foam
$\bF_{\mathrm{HL},s_i}$
defined by some smoothing $s_i$ which pairs the edges containing vectors $v_0$ and $v_i$ ---
see Section \ref{sec:hlfoamdef}.  Let $F_{ij}'$ be the face of $\bF_{\mathrm{HL},s_i}$
corresponding to $F$,
which deforms the face $F_{ij}$ of $\bF_{\mathrm{HL}}$ containing $v_i$ and $v_j$ (note $i,j \in \{0,1,2,3\}$).
Let us orient the arc from the end bounding the strand of the tangle deforming
the edge of $\bF_{\mathrm{HL}}$ with $v_i$
to the end bounding the tangle strand deforming the edge with $v_j$.
Call a vector along the arc in this orientation $v$.
We define the sign of the arc relative to the face by
$$\sigma(F,a) := \mathrm{sgn}\det(v,v_i,v_j) = \mathrm{sgn}\det(-v,v_j,v_i)$$

\begin{center}
\begin{tikzpicture}[scale=1.3]
\draw[fill,gray,opacity=.5] (0,0)--(1.5,0)--(1.5,1)--(-.7,.7)--(0,0);
\draw[red, very thick] (0,0)--(-.7,.7);
\draw[green, very thick] (0,0)--(1.5,0);
\draw[red, very thick] (1.5,0)--(1.5,1);
\draw[black,thick, -stealth] (0,0)--(.3,0);
\draw[black,thick,-stealth] (0,0)--(-.2,.2);
\draw[black,thick,-stealth] (1.5,0)--(1.5,.3);
\draw[black,thick,-stealth] (1.5,0)--(1.2,0);
\draw[black] (-.07,.07)--(0.03,.07)--(.1,0);
\draw[black] (1.4,0)--(1.4,.1)--(1.5,.1);

\node at (.2,-.2) {$v$};
\node at (1.3,-.2) {$-v$};
\node at (.75,.5) {$F$};
\node at (.75,-.4) {$a$};
\node at (-.5,.2) {$v_i$};
\node at (1.8,.2) {$v_j$};
\end{tikzpicture}
\end{center}

Note that the opposite orientation on the arc leads
to $\mathrm{sgn} \det (-v,v_j,v_i),$ which is the same.
The definition therefore only depends on the orientation of $B$.
\end{definition}

\begin{definition}
\label{def:facerelation}
Let $\bF$ be foam with boundary $\Gamma,$ and let ${\bF'}$ be
a deformed foam with arc set $A$.
We define a relation $\sim_{\bF'}$ on $\bZ^{E_\Gamma \cup A}$ by setting
\begin{equation}
\label{eq:facerelation}
\sum_{e\in \partial F} e + \sum_{a\in \partial F} \sigma(F,a) \cdot a
\,\sim_{\bF'}\, 0,
\end{equation}
for each face $F$ of $\bF'$.
\end{definition}

\subsection{Face relations for foams}

On general grounds, $L$ with $\partial L = S_\Gamma$ and
$b_1(L)=g = \frac{1}{2}b_1(S_\Gamma)$ defines a phase as the kernel of
the surjection $\tau:  H_1(S_\Gamma)\twoheadrightarrow H_1(L).$  Here we want to
understand this combinatorially when $L$ arises from a deformed foam $\bF'$,
in terms of its arcs and the edges of $\Gamma.$

\begin{proposition}
\label{prop:hladdition}
Let $\bF$ be an ideal foam filling a cubic graph $\Gamma$ with edge set $E_\Gamma$, and let $L$ be a smoothing
associated to a deformed foam $\bF'$ with arc set $A$, such that the corresponding tangle has no circle components.  Let
$\sim_{\bF'}$ be as in Definition \ref{def:facerelation}.
Then $H_2(L)=0$, and we have an isomorphism 
$$
H_1(L)\cong \bZ^{E_\Gamma\cup A}/\sim_{\bf F'}
$$
such that the homology pushforward $H_1(S_\Gamma)\rightarrow H_1(L)$ is identified with the map induced by the inclusion $\bZ^{E_\Gamma}\hookrightarrow \bZ^{E_\Gamma\cup A}$. 
\end{proposition}

Before the proof, a remark.
\begin{remark}
\label{rmk:nobigons}
If $\Gamma$ has no bigons, then
each edge $e$ is equivalent to a sum of arcs under Equation \ref{eq:facerelation}
by the external face of $\bF'$ containing $e$ in its
boundary.  
Then after taking the partial quotient of $\bZ^{E_\Gamma \cup A}\to \bZ^A$ by
the external faces of $\bF'$,
we may think of $H_1(L)$ as $\bZ^A/\sim_{\bF'}.$
\end{remark}

\begin{proof}
We prove the proposition by induction on the number of internal vertices of a foam $\bF.$

The base cases (no internal vertices) then consist of any of the canonical foam filling of necklace-type graphs of any genus, as in
Example \ref{ex:foam-necklace}.
Each such graph is itself obtained from the genus-$0$ necklace (theta graph) by bigon addition, or Legendrian
one-handle attachment of the corresponding
Legendrian surface --- see \cite[Theorem 4.10(1)]{CZ}) --- so we treat the base cases themselves by induction on the genus.
The genus-$0$ foam consists of the three filled semicircles in the unit ball at azimuthal angles $0,2\pi/3,4\pi/3.$
The edge lattice modulo face relations is zero, as is $H_1(L)$ for the filling, and the proposition is true.
Now we induct on the genus of the base case by adding bigons.
Each bigon addition adds three edges and one face to the boundary,
as seen here,
\begin{center}
\begin{tikzpicture}
\draw[very thick,blue] (-2,0)--(-1,0);
\node at (-3/2,.2) {$e$};
\node at (1/3,.2) {$e_-$};
\node at (5/3,.2) {$e_+$};
\node at (-1/2,0) {$\rightarrow$};
\draw[very thick,blue] (0,0)--(2/3,0);
\draw[very thick,blue] (1,0) circle (1/3);
\draw[very thick,blue] (4/3,0)--(2,0);
\end{tikzpicture}
\end{center}
thereby increasing 
$H_1$ of the Legendrian by $2$ and the genus by $1$.  Two faces are added to the foam,
which end in the two edges of the bigon.  The bigon edges sum to zero in
homology of the Legendrian, by the relation from the bigon face.
The foam face relations then show that
these edges are trivial in $H_1$ of the filling, thus in the kernel
of the homology map corresponding to inclusion of the boundary --- see Figure \ref{fig:necklacefoam}. 
The difference $e_+-e_-$ is in no boundary and therefore is an
additional nontrivial class
in $H_1$ of the Lagrangian filling the new Legendrian.  This establishes the base
case of no internal vertices, for every genus.

We now induct on the number of internal vertices by attaching a Harvey-Lawson foam. 
We can attach at a single vertex, along an edge, or a face. 
%
%
To verify the inductive step in the first
case, let $\bF'$ be a smoothed ideal foam, whose boundary is a cubic graph $\Gamma_{\bF}$ of genus $g$.

Now suppose $\Delta$ is a single tetrahedron together with a smoothed Harvey-Lawson foam $\bF_\Delta'$ in it. Let us choose a vertex $v$ of $\Gamma_{\bF}$ along with a vertex $w$ of the cubic graph $\Gamma_\Delta$ on the boundary of the tetrahedron. Let us write $e_1,e_2,e_3$ for the three edges of $\Gamma_{\bF}$ incident to the vertex $v$ listed in cyclic order determined by the orientation, and similarly write  $\epsilon_1,\epsilon_2,\epsilon_3$ for the edges of $\Gamma_\Delta$ incident to $w$, but listed in opposite cyclic order. Each of these edges determines an external face of the corresponding foam, which we denote by $f_{e_i}$ or $f_{\epsilon_i}$. We glue a neighborhood of the vertex $w$ to ${\bF}'$ by identifying the tetrahedron (dual) face corresponding to $w$ with the boundary (dual) face of ${\bF}$ corresponding to $v$ as indicated in Figure~\ref{fig:single-tetra-glue}, so that each edge $e_i$ is glued to the corresponding $\epsilon_i$ to form a new edge $\overline{e}_i$.  As a result of this gluing, we obtain a new ideal foam $\widetilde{\bF}'$. 
\begin{center}
\begin{figure}[ht]

\begin{tikzpicture}

        \draw[fill=brown]  (-3,4.5) -- (-1, 1.5) -- (-5, 1.5) ;
        \draw[fill=brown]  (1,2)--(5,2)--(3,0.5);

		\node [] (0) at (-3, 4.5) {};
		\node [] (1) at (-1, 1.5) {};
		\node [] (2) at (-3, 4.5) {};
		\node [] (3) at (-5, 1.5) {};
		\node [label={\small $e_1$}] (4) at (-5, 3.5) {};
		\node [circle,fill=blue,inner sep=0pt,minimum size=3pt,label={$v$}] (5) at (-3, 2.5) {};
		\node [] (6) at (-3, 0.5) {};
		\node [label={\small $e_2$}] (40) at (-3, -0.15) {};
		\node [label={\small $e_3$}] (7) at (-1, 3.5) {};

		\node [] (8) at (3, 0.5) {};
		\node [] (9) at (1, 2) {};
		\node [] (10) at (5, 2) {};
		\node [] (11) at (3, 4.5) {};
		\node [circle,fill=blue,inner sep=0pt,minimum size=3pt] (12) at (2.25, 2.5) {}; 
		\node [circle,fill=blue,inner sep=0pt,minimum size=3pt] (13) at (3.75, 2.5) {}; 
		\node [circle,fill=blue,inner sep=0pt,minimum size=3pt] (14) at (2.75, 1.45) {}; 
		\node [label={$w$}] (41) at (2.55, 0.85) {}; 
		\node [circle,fill=blue,inner sep=0pt,minimum size=3pt] (15) at (3.25, 3.05) {}; 
		\node [] (16) at (2, 3.25) {}; 
		\node [] (17) at (4, 3.25) {}; 
		\node [] (18) at (2, 1.25) {};
		\node [] (19) at (4, 1.25) {}; 

		\node [] (20) at (9.475, 4.5) {};
		\node [] (21) at (10.725, 1) {};
		\node [] (22) at (9.475, 4.5) {};
		\node [] (23) at (6.725, 1.75) {};
		\node [] (24) at (8.975, 2.5) {};
		\node [] (25) at (7.975, 2.75) {};
		\node [] (26) at (9.975, 2.75) {};
		\node [] (27) at (8.825, 1.25) {};
		\node [] (28) at (6.975, 3.25) {};
		\node [label=$\overline{e}_1$] (29) at (6.975, 3.25) {};
		\node [label=$\overline{e}_3$] (30) at (10.975, 3.25) {};
		\node [] (31) at (8.425, 0.5) {};
		\node [label=$\overline{e}_2$] (43) at (8.125, 0.25) {};
		\node [circle,fill=blue,inner sep=0pt,minimum size=3pt] (32) at (8.475, 2.75) {};
		\node [circle,fill=blue,inner sep=0pt,minimum size=3pt] (33) at (8.975, 2) {};
		\node [circle,fill=blue,inner sep=0pt,minimum size=3pt] (34) at (9.475, 2.75) {};

		\draw (0.center) to (1.center);
		\draw (2.center) to (3.center);
		\draw (3.center) to (1.center);
		\draw[blue] (4.center) to (5.center);
		\draw[blue] (5.center) to (7.center);
		\draw[blue] (5.center) to (6.center);
		\draw (9.center) to (8.center);
		\draw (8.center) to (10.center);
		\draw (11.center) to (9.center);
		\draw (11.center) to (10.center);
		\draw (9.center) to (10.center);
		\draw (11.center) to (8.center);
		\draw[blue] (12.center) to (16.center);
		\draw[blue,dashed] (16.center) to (15.center);
		\draw[blue,dashed] (15.center) to (17.center);
		\draw[blue,dashed] (15.center) to (14.center);
		\draw[blue] (12.center) to (13.center);
		\draw[blue] (17.center) to (13.center);
		\draw[blue,dashed] (14.center) to (18.center);
		\draw[blue,dashed] (14.center) to (19.center);
		\draw[blue] (12.center) to (18.center);
		\draw[blue] (13.center) to (19.center);
		
		\draw (20.center) to (21.center);
		\draw (22.center) to (23.center);
		\draw (23.center) to (21.center);
		\draw (22.center) to (24.center);
		\draw (24.center) to (23.center);
		\draw (24.center) to (21.center);
		\draw[blue] (29.center) to (25.center);
		\draw[blue] (30.center) to (26.center);
		\draw[blue] (27.center) to (31.center);
		\draw[blue] (25.center) to (32.center);
		\draw[blue] [in=135, out=45, looseness=1.25] (32.center) to (34.center);
		\draw[blue]	 [in=30, out=-90] (34.center) to (33.center);
		\draw[blue] [in=165, out=-90] (32.center) to (33.center);
		\draw[blue] (34.center) to (26.center);
		\draw[blue] (33.center) to (27.center);
\end{tikzpicture}

\caption{Attaching a Harvey-Lawson foam to ${\bF}$ in a neighborhood $\widehat{f}$ of a vertex, shaded in brown. (Dually, $\widehat{f}$ is a face of the dual triangulation to
$\Gamma_{\bF}$.)}
\label{fig:single-tetra-glue}
\end{figure}
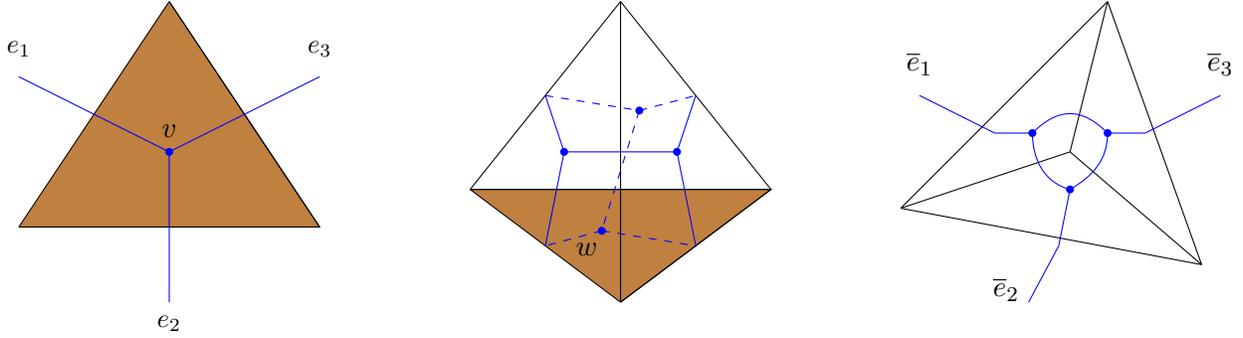
\end{center}

The set of faces of the new deformed foam $\widetilde{\bF}'$
may be described as follows. The internal faces of $\widetilde{\bF}'$ are the same as those of $\bF'$. The set of external faces of $\widetilde{\bF}'$ consists of all those external faces of $\bF'$ and $\bF_\Delta'$ that correspond to edges of $\Gamma_{\bF}$ or $\Gamma_\Delta$ not incident to $v,w$, along with three faces $f_1,f_2,f_3$ obtained by gluing each $f_{e_i}$ to the corresponding $f_{\epsilon_i}$.

We now turn our attention to the effect of this gluing at the level of the double covers of the ball. Let us write $\pi:L\rightarrow B$,
$\pi_\epsilon:\mathrm{HL}_{\epsilon}\rightarrow B$ for the branched double covers corresponding to ${\bF}'$ and $\bF_\Delta'$ respectively, so that we have $\widetilde{L} = L \cup_{\pi^{-1}(\hat{f})} \mathrm{HL}_{\epsilon}$, where $\hat{f}$ is the neighborhood of the vertex $v$ where we attached $\bF_{\Delta}'$, i.e.~the face of
along which the dual tetrahedron $\Delta$ was glued. This gluing is illustrated in Figure~\ref{fig:single-face-glue-cover}. 
\begin{center}
\begin{figure}[ht]

\begin{tikzpicture}[scale=.75]
	
		\node [] (0) at (-11.25, 1) {};
		\node [] (1) at (-3.25, 2) {};
		\node [] (2) at (-11.25, 9) {};
		\node [] (3) at (-3.25, 10) {};
		\node [] (4) at (-7.25, 1.5) {};
		\node [] (5) at (-7.25, 9.5) {};
		\node [] (6) at (-11.25, 5) {};
		\node [] (7) at (-3.25, 6) {};
		\node [] (8) at (-8, 4.5) {};
		\node [] (9) at (-6.75, 4){};
		\node [] (10) at (-6, 5) {};
		\node [] (11) at (-6.5, 6.5) {};
		\node [] (12) at (-7.75, 6.75) {};
		\node [] (13) at (-8.5, 5.75) {};
		\node [] (14) at (0, 10) {};
		\node [] (15) at (8, 9) {};
		\node [] (16) at (0, 2) {};
		\node [] (17) at (8, 1) {};
		\node [] (18) at (4, 9.5) {};
		\node [] (19) at (4, 1.5) {};
		\node [] (20) at (0, 6) {};
		\node [] (21) at (8, 5) {};
		\node [] (22) at (3.25, 6.5) {};
		\node [] (23) at (4.5, 7) {};
		\node [] (24) at (5.25, 6) {};
		\node [] (25) at (4.75, 4.5) {};
		\node [] (26) at (3.5, 4.25) {};
		\node [] (27) at (2.75, 5.25) {};
		\node [] (28) at (-7.5, 7.75) {\tiny $e_3$};
		\node [] (29) at (-9.6, 7.25) {\tiny $e_1$};
		\node [] (32) at (8.25, 2.5) {};
		\node [] (34) at (4.25, 8) {\tiny $\epsilon_3'$};
		\node [] (35) at (6.75, 7.5) {\tiny $\epsilon_1'$};
		\node [] (36) at (2, 1.5) {};
		\node [] (37) at (-5, 6.10) {\tiny $e_2$};
		\node [] (40) at (1.5, 6.10) {\tiny $\epsilon_2'$};
		\node [] (41) at (-8.45, 7.75) {$-$};
		\node [] (42) at (-9.75, 6.25) {$+$};
		\node [] (43) at (-8.5, 3.75) {$-$};
		\node [] (44) at (-6, 3) {$+$};
		\node [] (45) at (-4.5, 4.5) {$-$};
		\node [] (46) at (-5.75, 7.25) {$+$};
		\node [] (47) at (2.5, 7.25) {$+$};
		\node [] (48) at (1.5, 4.5) {$-$};
		\node [] (49) at (2.75, 3.25) {$+$};
		\node [] (50) at (5.5, 3.75) {$-$};
		\node [] (51) at (6.75, 6.25) {$+$};
		\node [] (52) at (5.25, 8) {$-$};

		\draw[fill=brown]  (-8, 4.5) -- (-6.75, 4) -- (-6, 5)  -- (-6.5, 6.5)--(-7.75, 6.75) -- (-8.5, 5.75);
				\draw[fill=brown]  (3.25,6.5) -- (4.5,7) -- (5.25,6)  -- (4.75,4.5)--(3.5,4.25) -- (2.75,5.25);

		\draw[teal] (0.center) to (1.center);
		\draw[teal] (2.center) to (3.center);
		\draw[blue,<-] (4.center) to (5.center);
		\draw[magenta] (0.center) to (2.center);
		\draw[magenta] (1.center) to (3.center);
		\draw[red,->] (6.center) to (7.center);
		\draw[pink] (4.center) to (6.center);
		\draw[green,->] (1.center) to (2.center);
		\draw[pink] (7.center) to (5.center);
		\draw (8.center) to (9.center);
		\draw (9.center) to (10.center);
		\draw (10.center) to (11.center);
		\draw (11.center) to (12.center);
		\draw (13.center) to (12.center);
		\draw (8.center) to (13.center);
		\draw[teal] (14.center) to (15.center);
		\draw[teal] (16.center) to (17.center);
		\draw[blue,<-] (18.center) to (19.center);
		\draw[magenta] (14.center) to (16.center);
		\draw[magenta] (15.center) to (17.center);
		\draw[red,->] (20.center) to (21.center);
		\draw[pink] (18.center) to (20.center);
		\draw[green,->] (15.center) to (16.center);
		\draw[pink] (21.center) to (19.center);
		\draw (22.center) to (23.center);
		\draw (23.center) to (24.center);
		\draw (24.center) to (25.center);
		\draw (25.center) to (26.center);
		\draw (27.center) to (26.center);
		\draw (22.center) to (27.center);
\end{tikzpicture}

\caption{Gluing two copies of $\mathrm{HL}_{\epsilon}$ along the disk $\pi^{-1}(\hat{f})\subset \partial \mathrm{HL}_{\epsilon}\simeq T^2$. The set $\pi^{-1}(\hat{f})$ is shaded brown, and opposite pairs of boundary edges of the square are identified in the figure. The arrows on edges indicate the canonical lifts of the edges of the cubic graphs on $S^2$ to  cycles in $H_1(\Lambda)$. External faces of the foam are labelled by the corresponding sign of the primitive function.} 
\label{fig:single-face-glue-cover}
\end{figure}
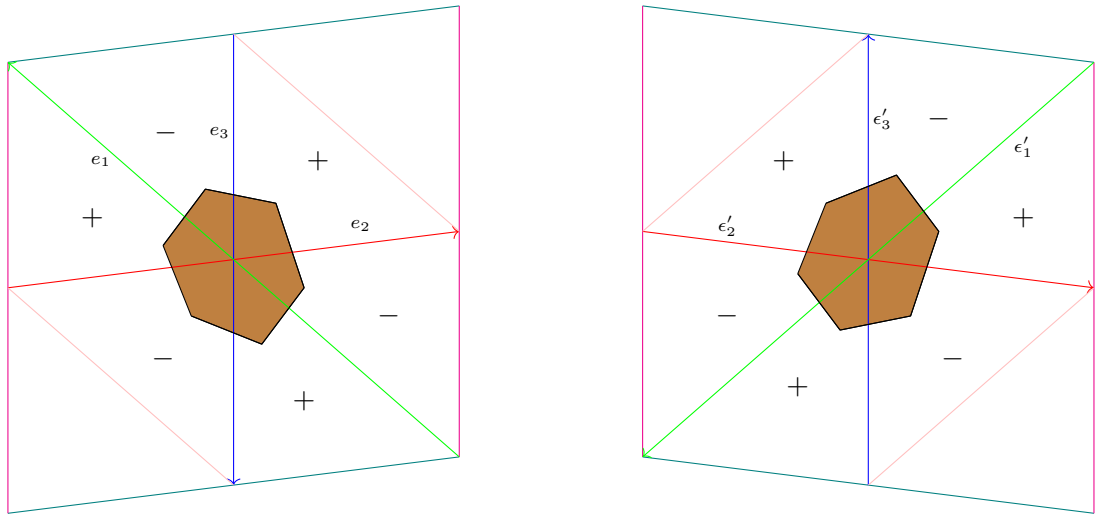
\end{center}
Now since the space $\pi^{-1}(\hat{f})$ is homeomorphic to a disk, the Mayer-Vietoris long exact sequence shows that $H_2(\widetilde{L})\simeq H_2(L)\oplus H_2({\mathrm{HL}_{\epsilon}})=0$. Similarly, it delivers an isomorphism
\begin{equation}
\label{mv-iso-1-face}
i_*+\iota_*: H_1(L)\oplus H_1({\mathrm{HL}_{\epsilon}}) \rightarrow H_1(L'),
\end{equation}
such that $[\overline{e}_i] = i_*([e_i])+\iota_*([\epsilon_i])$.
Hence all that remains is to verify the face relations for the faces
$f_1,f_2,f_3$ of $\widetilde{\bF}'$ obtained by gluing faces of $\bF'$ to those of $\bF_\Delta'$. But recall from \ref{def:arc-sign} that the definition of the sign of an arc $a$ relative to a face $f$ is entirely local,  depending only on the tangent vectors $v_i,v_j$ to the two edges of the deformed foam that meet $a$ and bound $f$. So if $a_0$ is the unique arc at the vertex of ${\bF}'$ that is connected to $v$ by an edge,
and $a_1$ the corresponding arc in $\bF_\Delta'$ (connected to $w$), the sign of $a_0$ with respect to face $f_{e_i}$ in $\bF$ is identical to its sign with respect to face $f_i$ of the glued foam $\bF'$. Similarly, the sign of $a_1$ with respect to $f_{\epsilon_i}$ coincides with its sign with respect to $f_i$. The face relation for $f_i$ is therefore obtained as the sum of those for $f_{\epsilon_i}$ and $f_{e_i}$ under the isomorphism~\eqref{mv-iso-1-face}. 

We next consider the case of gluing in a
Harvey-Lawson cone along an edge. 
We have a foam $\bF$ with boundary $\Gamma_{\bF}$
and a Harvey-Lawson foam $\bF_\Delta$ with boundary
a tetrahedron graph $\Gamma_{\Delta}.$
Suppose that we fix an edge $e_0$ of $\Gamma_{\bF}$ connecting two vertices $v_1,v_2$, and correspondingly fix an edge $\epsilon_0$ of $\Gamma_{\Delta}$ connecting vertices $w_1,w_2$ of $\Delta$. Let us denote the edges of $\Gamma_{\bF}$ incident to $v_1$ by $e_0,e_1,e_2$, cyclically ordered in accordance with the orientation of $\partial B$, and  similarly write $e_0,e_3,e_4$ for the edges incident to $v_2$. We denote by $\epsilon_0,\epsilon_1,\epsilon_2$ the edges incident to $w_1$ but ordered with respect to the opposite of the orientation on $\Delta$, and similarly write $\epsilon_0,\epsilon_3,\epsilon_4$. We write $\epsilon_5$ for the remaining edge of $\Gamma_\Delta$ which is incident to neither $w_1$ nor $w_2$. We now glue the
foam $\bF_\Delta$ to $\bF$ by identifying the
edges so that each edge $e_i$ is glued to the corresponding $\epsilon_i$.  We denote by $\widetilde{\bF}$ the ideal foam produced as a result of this gluing.
\begin{center}
\begin{figure}[ht]
\begin{tikzpicture}
		\node [] (0) at (-6, 0) {};{rl};
		\node [] (1) at (0, 0) {};{rr};
		\node [] (2) at (-3.25, 2.75) {};{rt};
		\node [] (3) at (-2.5, -2.25) {};{rb};
		\node [circle,fill=blue,inner sep=0pt,minimum size=3pt] (4) at (-3.5, 1) {};
		\node [] (5) at (-4.75, 1.25) {};
		\node [] (6) at (-1.25, 1.25) {};
		\node [circle,fill=blue,inner sep=0pt,minimum size=3pt] (7) at (-4, 0.75) {};
		\node [] (20) at (-3.85, 0.55) { $v_1$};;
		\node [] (8) at (-0.75, -0.75) {};
		\node [] (10) at (-5.25, -0.5) {};
		\node [circle,fill=blue,inner sep=0pt,minimum size=3pt] (11) at (-3, -0.25) {};
		\node [] (12) at (-6.5, 2.5) {};{xt};
		\node [] (13) at (-6.5, -1.75) {};{xt};
		\node [] (14) at (0.5, 2.75) {};{xt};
		\node [] (15) at (0.75, -1.75) {};{xt};
		\node [circle,fill=blue,inner sep=0pt,minimum size=3pt] (16) at (-2, 0.75) {};
		\node [] (21) at (-2.1, 0.55) {$v_2$};
		\node [] (22) at (-3.6,1.2) {$w_1$};
		\node [] (23) at (-3,-.5) {$w_2$};

		\draw[thick,red] (2.center) to (0.center);
		\draw[thick,red] (0.center) to (3.center);
		\draw[thick,red] (3.center) to (1.center);
		\draw[thick,red] (2.center) to (1.center);
		\draw[thick,red] (0.center) to (1.center);
		\draw[thick,red,dashed] (2.center) to (3.center);
		\draw[thick,blue] (12.center) to (5.center);
		\draw[thick,blue] (13.center) to (10.center);
		\draw[thick,blue] (14.center) to (6.center);
		\draw[thick,blue] (8.center) to (15.center);
		\draw[thick,blue] (5.center) to (4.center);
		\draw[thick,blue] (4.center) to (11.center);
		\draw[thick,blue] (10.center) to (11.center);
		\draw[thick,blue] (11.center) to (8.center);
		\draw[thick,blue] (4.center) to (6.center);
		\draw[thick,blue,dashed] (5.center) to (7.center);
		\draw[thick,blue,dashed] (10.center) to (7.center);
		\draw[thick,blue,dashed] (7.center) to (16.center);
		\draw[thick,blue,dashed] (6.center) to (16.center);
		\draw[thick,blue,dashed] (16.center) to (8.center);
\end{tikzpicture}
\caption{The cubic graph (shown solid in blue) produced by gluing a
Harvey-Lawson foam an edge is related to the
orignal (blue, dotted) by a diagonal exchange.
The Harvey-Lawson dual tetrahedron is shown in red.}
\label{fig:double-tetra-flip}
\end{figure}
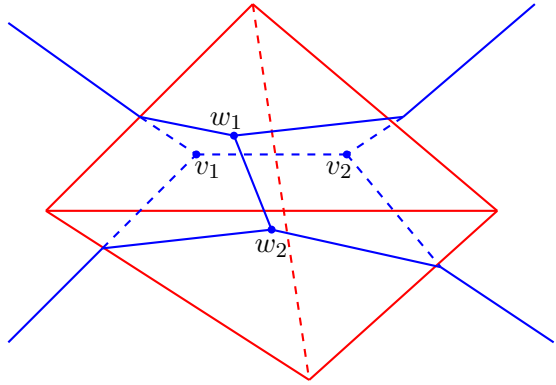
\end{center}

\begin{center}
\begin{figure}[ht]
\begin{tikzpicture}

		\node [] (0) at (0.5, 1) {};
		\node [] (1) at (8.5, 2) {};
		\node [] (2) at (0.5, 9) {};
		\node [] (3) at (8.5, 10) {};
		\node [] (4) at (4.5, 1.5) {};
		\node [] (5) at (4.5, 9.5) {};
		\node [] (6) at (0.5, 5) {};
		\node [] (7) at (8.5, 6) {};
		\node [] (8) at (3.75, 4.5) {};
		\node [] (9) at (5, 4) {};
		\node [] (10) at (5.75, 5) {};
		\node [] (11) at (5.25, 6.5) {};{h1l};
		\node [] (12) at (4, 6.75) {};{h0l};
		\node [] (13) at (3.25, 5.75) {};
		\node [] (23) at (0.5, 6.5) {};
		\node [] (24) at (1.75, 2.5) {};
		\node [] (25) at (6.75, 1.75) {};
		\node [] (26) at (2.5, 1.25) {};
		\node [] (29) at (7, 8.5) {};
		\node [] (30) at (7, 9.75) {};
		\node [] (31) at (8.5, 8.75) {};
		\node [] (32) at (2.5, 9.25) {};
		\node [] (33) at (0.5, 2.5) {};
		\node [] (34) at (8.5, 4.5) {};

		\draw[fill=brown]  (3.75, 4.5)-- (5, 4) --(5.75, 5)--(5.25, 6.5)--(4, 6.75)--(3.25, 5.75);
		\draw[fill=brown]  (4, 6.75)-- (2.5, 9.25) --(7, 9.75)--(7, 8.5)--(5.25, 6.5);
		\draw[fill=brown]  (3.75, 4.5)-- (1.75, 2.5) --(2.5, 1.25)--(6.75, 1.75)--(5, 4);

		\draw[blue,->>] (0.center) to (1.center);
		\draw[blue, ->>] (2.center) to (3.center);
		\draw[blue] (4.center) to (5.center);
		\draw[blue,->] (0.center) to (2.center);
		\draw[blue,->] (1.center) to (3.center);
		\draw[blue] (6.center) to (7.center);
		\draw[blue] (4.center) to (6.center);
		\draw[blue] (1.center) to (2.center);
		\draw[blue] (7.center) to (5.center);
		\draw[] (8.center) to (9.center);
		\draw[] (9.center) to (10.center);
		\draw (10.center) to (11.center);
		\draw (11.center) to (12.center);
		\draw (13.center) to (12.center);
		\draw (8.center) to (13.center);
		\draw (12.center) to (32.center);
		\draw (11.center) to (29.center);
		\draw (30.center) to (29.center);
		\draw (29.center) to (31.center);
		\draw (13.center) to (23.center);
		\draw (8.center) to (24.center);
		\draw (24.center) to (33.center);
		\draw (24.center) to (26.center);
		\draw (9.center) to (25.center);
		\draw (10.center) to (34.center);

\end{tikzpicture}
\caption{Shaded in brown is the subset $\pi^{-1}(\hat{\sigma})$ of $\partial \mathrm{HL}_{\epsilon} = T^2$. Edges of the cubic graph are shown in blue, and those of the dual triangulation in black. Since the opposite pairs of blue boundary edges are identified, the space $\pi^{-1}(\hat{\sigma})$ is homeomorphic to a cylinder.}
\label{fig:cylinder-glue}
\end{figure}
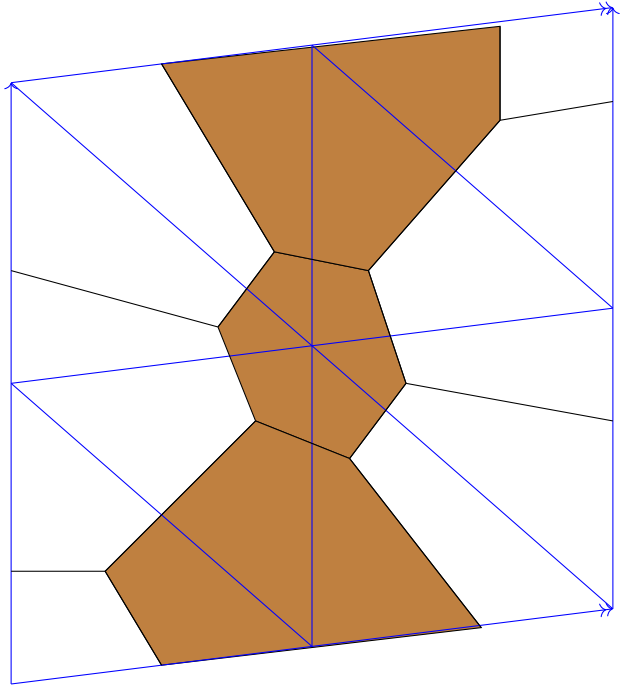
\end{center}


Note that the cubic graphs $\Gamma_{\widetilde{\bF}}$ and $\Gamma_{\bF}$ have the same genus: indeed, the two are related by a single diagonal exchange/edge mutation, as illustrated in Figure~\ref{fig:double-tetra-flip}. 
(We will return to this point in Proposition \ref{prop:mutationtransfer}.)
The set of external faces of $\widetilde{\bF}$ is thus in natural bijection with that of $\bF$: the latter contains the external faces $\{f_i, i=1,\ldots 4\}$ obtained by gluing each face $f_{e_i}$ to the corresponding $f_{\epsilon_i}$, along with the external face with
boundary $f_{\epsilon_5}$.  On the other hand, we now have a new \emph{internal} face $f_0$ created by gluing $f_{e_0}$ to $f_{\epsilon_0}$. By assumption, the smoothing of the foam in $\Delta$ is chosen such that the tangle in the glued
deformed foam $\widetilde{\bF}'$ has no circle components; this is equivalent to requiring that at least one of the faces $f_{e_0},f_{\epsilon_0}$ contains an arc as part of its boundary. 

We now consider the gluing of double covers $\pi:L\rightarrow B$
and $\pi_\epsilon:\mathrm{HL}_{\epsilon}\rightarrow B$. We have $\widetilde{L} = L \cup_{\pi^{-1}(\hat{\sigma})} \mathrm{HL}_{\epsilon}$, where $\hat{\sigma}$ is 
a neighborhood in $B$ of the edge $\epsilon_0$, or dually 
the quadrilateral along which the dual tetrahedron $\Delta$ is glued to $B$. As shown in Figure~\ref{fig:cylinder-glue} the space $\pi^{-1}(\hat{\sigma})$ is homeomorphic to a cylinder $C$. We fix the isomorphism $H_1(C)\simeq \mathbb{Z}\{\gamma\}$, where we take the generator $\gamma$ to be the oriented loop on $\partial \mathrm{HL}_\epsilon$ given by canonical lift of the edge $\epsilon_0$ of $\Gamma_{\Delta}$.  The relevant part of Mayer-Vietoris sequence then reads
\begin{equation}
\label{mv-iso-2-face}
0\rightarrow H_2(\widetilde{L})\rightarrow \mathbb{Z}\{\gamma\}\rightarrow H_1(L)\oplus H_1({\mathrm{HL}_{\epsilon}}) \rightarrow H_1(\widetilde{L})\rightarrow 0.
\end{equation}
By our assumption that at least one of the faces $f_{e_0},f_{\epsilon_0}$ contains an arc as part of its boundary, we see that the map 
$$
i_*\oplus (-\iota_*) : \mathbb{Z}\{\gamma\}\rightarrow H_1(L)\oplus H_1({\mathrm{HL}_{\epsilon}})
$$ is injective. Hence $H_2(\widetilde{L})=0$, and 
$$
H_1(\widetilde{L}) \simeq \frac{H_1(L)\oplus H_1({\mathrm{HL}_{\epsilon}})}{\mathbb{Z}\{\gamma\}}.
$$
The face relations for all external faces of $\widetilde{\bF}'$ now follow from this description of $H_1(\widetilde{L})$ exactly as in the case of the single-triangle gluing. Finally, since 
$$
\left(i_*\oplus (-\iota)_*\right)(\gamma)  = ([e_0],[\epsilon_0]),  
$$
we see from the isomorphism~\eqref{mv-iso-2-face} that the relation in $H_1(\widetilde{L})$ corresponding to the new internal face $f_0$ is also obtained as the sum of the relation corresponding to $f_{e_0}$ in $H_1(L)$ with that corresponding to $f_{\epsilon_0}$ in $H_1({\mathrm{HL}_{\epsilon}})$. 

It remains to consider attaching along a triangular face.
The proof is very similar to the above case, so we only comment briefly. 
In this case, the attachment is along a punctured torus, so $H_1(\pi^{-1}(\hat{\sigma}))$ has rank two. 
It still injects into $H_1(L) \oplus H_1(HL_\epsilon)$, and otherwise the exact sequence looks the same.
Therefore $H_1(\widetilde{L})$ has rank one less than $H_1(L)$.  The rest of the proof is as above.

This completes the proof of the Proposition.

\end{proof}


\subsection{Example -- triangular prism}
\label{sec:prism}

Let $\Gamma$ be the edge graph of a triangular prism and let $\bF'$
be the deformed foam pictured here:

\begin{center}
\begin{tikzpicture}[scale=1.3]
\pgfmathsetmacro{\Ax}{0}
 \pgfmathsetmacro{\Ay}{3}
 \pgfmathsetmacro{\Bx}{2.3}
 \pgfmathsetmacro{\By}{3.3}
 \pgfmathsetmacro{\Cx}{1}
 \pgfmathsetmacro{\Cy}{4}
 \pgfmathsetmacro{\Dx}{0}
 \pgfmathsetmacro{\Dy}{0}
 \pgfmathsetmacro{\Ex}{2.3}
 \pgfmathsetmacro{\Ey}{0.3} 
\pgfmathsetmacro{\Fx}{1}
 \pgfmathsetmacro{\Fy}{1}
 \coordinate (A) at (\Ax,\Ay);
 \coordinate (B) at (\Bx,\By);
 \coordinate (C) at (\Cx,\Cy);
 \coordinate (D) at (\Dx,\Dy);
 \coordinate (E) at (\Ex,\Ey);
  \coordinate (F) at (\Fx,\Fy);
 \coordinate (Amid) at (.7,1.5);
  \coordinate (Amidtext) at (.7-.3,1.5);
 \coordinate (Amid2) at (.5,1.7);
 \coordinate (Bmid) at (1.5,1.6);
 \coordinate (Bmid2) at (1.7,1.7);
  \coordinate (Bmidtext) at (1.6,1.5);
 \coordinate (Cmid) at (1.32,2.3);
 \coordinate (Cmid2) at (1.2,2.1);
  \coordinate (Cmidtext) at (1.15,2.25);

 \coordinate (ACmid) at (\Ax/2+\Cx/2-.4,\Ay/2+\Cy/2);
 \coordinate (BCmid) at (\Bx/2+\Cx/2+.3,\By/2+\Cy/2+.1);
 \coordinate (ABmid) at (\Ax/2+\Bx/2+.4,\Ay/2+\By/2-.3);
 \coordinate (DFmid) at (\Dx/2+\Fx/2+.2,\Dy/2+\Fy/2-.1);
 \coordinate (EFmid) at (\Ex/2+\Fx/2,\Ey/2+\Fy/2+.2);
 \coordinate (DEmid) at (\Dx/2+\Ex/2+.2,\Dy/2+\Ey/2-.25);
 \coordinate (ADmid) at (\Ax/2+\Dx/2-.3,\Ay/2+\Dy/2);
 \coordinate (CFmid) at (\Fx/2+\Cx/2-.3,\Fy/2+\Cy/2+.1);
 \coordinate (BEmid) at (\Bx/2+\Ex/2+.3,\By/2+\Ey/2);

\draw[fill,black,opacity=.4] (Amid2)--(Amid)--(Bmid)--(Bmid2)--(Cmid2)--(Amid2);
\draw[fill,red,opacity=.4] (A)--(Amid2)--(Cmid2)--(Cmid)--(C)--(A);

\draw[red, thick] (D)--(Amid)--(Bmid)--(B);
\draw[red, thick] (A)--(Amid2)--(Cmid2)--(Bmid2)--(E);
\draw[red, thick] (C)--(Cmid);
\draw[red, thick, dashed] (Cmid)--(F);
\draw[green, ultra thick] (Amid)--(Amid2);
\draw[green, ultra thick] (Cmid)--(Cmid2);
\draw[green, ultra thick] (Bmid)--(Bmid2);

\draw[blue,thick] (A)--(B)--(C)--(A)--(D)--(E)--(B); 
\draw[blue,thick,dashed] (D)--(F)--(E);
\draw[blue,thick,dashed] (F)--(C);

\node at (Amidtext) {$b$};
\node at (Bmidtext) {$c$};
\node at (Cmidtext) {$a$};

\node at (ACmid) {$T_1$};
\node at (BCmid) {$T_2$};
\node at (ABmid) {$T_3$};
\node at (DFmid) {$B_1$};
\node at (EFmid) {$B_2$};
\node at (DEmid) {$B_3$};
\node at (ADmid) {$L_1$};
\node at (CFmid) {$\; L_2$};
\node at (BEmid) {$L_3$};

\end{tikzpicture}
\end{center}

Write $G$ for the gray face and $P$ for the pink face.  Then $\sigma(P,a) = 1,$
and from Equation \ref{def:facerelation},
$P$ gives the relation $$\tau(T_1)+\gamma_a = 0.$$  
In total, the external face relations give
$$\begin{array}{rclccrclccrcl}
\tau(T_1)+\gamma_a &=& 0 && \tau(L_1)+\gamma_b &=& 0&& \tau(B_1)-\gamma_b - \gamma_a &=& 0\\
\tau(T_2)-\gamma_a-\gamma_c &=& 0 && \tau(L_2)&=&0&& \tau(B_2)+ \gamma_a &=& 0\\
\tau(T_3)-\gamma_b &=& 0&&\tau(L_3)+\gamma_c &=& 0&& \tau(B_3)- \gamma_c &=& 0
\end{array}$$
We also have the internal (gray) face relation, and since $\sigma(G,b)=\sigma(G,c)=1,$ we see
$$\gamma_b+\gamma_c = 0.$$
The relations are consistent with the face relations from $\Gamma$.  For example, the sum
$T_1 + T_2 + T_3 \sim_{\Gamma} 0$, and this implies $\tau(T_1+T_2+T_3)=0,$ or
$\gamma_b + \gamma_c = 0,$ and this is true by the internal gray face relation of $\bF'$.
The other face relations are consistent, as well. 

So $\tau$ indeed descends from a map from $\bZ^{E_\Gamma\cup A}$ to
one from $H_1(S_\Gamma)\cong \bZ^{E_\Gamma}/\sim,$ giving a map to $H_1(L)=\bZ^A/\sim$.
$H_1(S_\Gamma)$ is rank-$4$ and we can take Darboux generators $T_1,T_2;B_2,B_1$
(careful about the cyclic order on the back side of the prism:  $\overline{\omega}(B_2,B_1)=1$).
$H_1(L)$ is rank-$2$ and we can take generators $\gamma_a,\gamma_b$.
With these generators,
$$\tau(T_1) = -\gamma_a, \quad
 \quad \tau(T_2) = \gamma_a - \gamma_b,\quad
 \tau(B_2) = -\gamma_a,\quad
 \tau(B_1) = \gamma_a+\gamma_b
.$$
We can see that the kernel of $\tau$ is generated by
$\mu_1 := -(T_1+T_2+B_1+B_2)$ and
$\mu_2 := T_1 - B_2,$
and is indeed isotropic.
A framing $H_1(L)\hookrightarrow H_1(S_\Gamma)$
must send $\gamma_a$ to $-T_1 + \alpha \mu_1 + \beta \mu_2$
and $\gamma_b$ to $B_1 + B_2 + \gamma \mu_1 + \delta \mu_2$
The image is isotropic if $\beta =\gamma,$ so the different framings
for this phase are parametrized
by symmetric $2\times 2$ integer matrices $\binom{\tiny{\alpha\;\beta}}{\tiny{\beta\;\delta}}.$

\subsection{Associated Cones, Geometric Cones}

To formulate open Gromov-Witten conjectures, we want to
express a wavefunction in a power series about a limit point
of the moduli space.  A phase and framing define
an algebraic torus, but pinning down a limit point for the expansion requires
the notion of an associated cone, which we define below
after setting notation.

We write $\Gamma$ for the underlying cubic graph, $S_\Gamma$ for the associated
Legendrian surface, $\bF'$ for the deformed foam, and
$L$ for the corresponding Lagrangian.  

%


\begin{definition}
Given a splitting, i.e.~a phase
$K = \mathrm{Ker}(\,\tau: H_1(S_\Gamma\twoheadrightarrow H_1(L)$)
and a framing $F \subset H_1(S_\Gamma)$
$\omega$-isotropic and transverse to $K$ (so $\tau: F
\stackrel{\textstyle\sim}{\smash{\longrightarrow}\rule{0pt}{0.4ex}}
 H_1(L)$), an \emph{associated cone}
(or just \emph{cone}) is an open integral convex cone
$C_F \subset F\cong H_1(L)$ containing no lines. 
\end{definition}

With Remark \ref{rmk:nobigons} in mind,
if $\Gamma$ is simple (in particular has no bigons)
and we are given a splitting, then
we can specify an associated cone
by choosing a spanning set of arcs.

\begin{definition}
Define a \emph{geometric cone} of a deformed foam $\bF'$ with
$\partial\bF' = \Gamma$
to be the $\bZ_{\geq 0}$ span of a spanning set of
arcs and edges in $\bZ^{E_\Gamma \cup A}/\sim_{\bF'}.$
When $\Gamma$ is simple, without loss of generality we
take a spanning set of arcs.
\end{definition}

\begin{example}
\label{ex:neckseed}
Let $\bF' = \bF$ be the foam for the necklace graph
$\Gamma_{\mathrm{neck}}^g$ (see Remark \ref{rmk:neckfoam}),
and $L$ the corresponding Lagrangian.
Label the beads $0$ through $g$ in clockwise order
from some chosen starting point,
and let $b_i,$ $i= 0,...,g$,
be edges along the outer edges of the corresponding bead.
Label the strands of the necklace
$a_i$, $i = 0,...,g,$ so that
the $i$th strand succeeds the $i$th bead in clockwise order
and $\omega(a_i,b_i) = 1.$ 
The $b_i$ span the kernel of $\tau:
H_1(S_{\Gamma_{\mathrm{neck}}^g})\to H_1(L)$, so define a phase.
(The strands $a_i$ function as arcs, albeit there are no vertices,
as they connect tangle components.)  The map $a_i \mapsto a_i$ defines
a splitting $H_1(L)\to H_1(S_{\Gamma_{\mathrm{neck}}^g}).$ 
We note the following
relations in $H_1(S_{\Gamma_{\mathrm{neck}}^g})$:  $\sum a_i = 0, \sum b_i = 0.$
So any $g$-element subset of $\{a_i\}$ determines
a geometric cone, and by symmetry we may as well take this to be
$a_1,...,a_g.$  
The necklace therefore has a unique (up to symmetry)
phase, framing and geometric cone.

\end{example}

\begin{example}
A Harvey-Lawson smoothing has a single arc
and therefore a unique geometric cone.
The blue edges are equivalent under the face relations
and span the kernel (phase) of $\pi: H_1(S_{\Gamma_\Delta})\to H_1(L).$
A splitting is defined by mapping the green arc
to a transverse element of $H_1(\Lambda),$
and the unique
associated cone is the $\bZ_{\geq 0}$ span of this vector.
\begin{center}
\begin{tikzpicture}
\pgfmathsetmacro{\a}{.8};
\pgfmathsetmacro{\b}{.3};
\pgfmathsetmacro{\c}{.2};
\pgfmathsetmacro{\d}{.6};

\coordinate (ar) at (.2,0);
\coordinate (al) at (.1,.1);

\coordinate (b) at (\a,-\b);
\coordinate (c) at (-\a,-\c);
\coordinate (d) at (-\c,-\d);
\coordinate (e) at (\c,1);
\draw[thick] (b)--(c)--(d)--(e)--(b)--(d);
\draw[thick] (c)--(e);
\draw[thick, red] (d)--(ar)--(b);
\draw[thick, red] (e)--(al)--(c);
\draw[very thick, green] (al)--(ar);
\draw[very thick,blue] (b)--(d);
\draw[very thick,blue] (c)--(e);

\end{tikzpicture}
\end{center}

\end{example}

\begin{example}
Let $\bF'$ be as in the Example of Section \ref{sec:prism}.
Given that $\Gamma$ is simple, per Remark \ref{rmk:nobigons}
and the fact that the internal face relation is
$a+b\sim_{\bf'} 0,$
we conclude that
there are two choices of geometric cones: $\{a,b\}$
and $\{a,c\}.$
\end{example}

\subsection{Mutations of Foams and Cones}

We now show that for a large class of mutations of
the boundary graph, the foam filling can be mutated,
along with a phase, framing and cone.

\begin{definition}
\label{def:allowablemutation}
A mutation of a deformed foam
at an edge $e\in \Gamma$ is \emph{allowable} if
$e$ is not the boundary of a single tangle strand.
\end{definition}
The reason for this definition is to exclude the case
where the class $[e]\in H_1(S_\Gamma)$ is in the kernel of $\tau:
H_1(S_\Gamma)\to H_1(L)$,
rendering the action on the wavefunction zero.
At the level of tangles, the condition ensures that the new tangle has no circle component.

\begin{proposition}
\label{prop:foammutation}
Let $\Gamma$ be a cubic planar graph bounding
a deformed foam $\bF'$.
Let $\Gamma_e$ be the graph defined by performing
an allowable mutation
at the edge $e\in E(\Gamma)$.  
Then there is either one or two
canonically defined deformed
foams $\bF'_{e,+}$ and $\bF'_{e,-}$
with boundary $\Gamma_e,$
corresponding to positive and negative mutations, respectively. 
\end{proposition}

\begin{proof}
The proposition follows immediately from the proof of Proposition \ref{prop:hladdition}:
the allowed mutations correspond to attaching a Harvey-Lawson foam along an edge, with
the allowable condition corresponding to the hypothesis that the tangle
of the deformed foam have no circle components.
Nevertheless,
for the convenience of the reader, we provide a separate description in the language
of triangulations --- though they are not as general as foams (the foams of necklace-type graphs
are degenerate tetrahedronizations),
they are often easier to visualize.

We first mutate the ideal foam $\bF$,
then worry about its deformation $\bF'$.
On the surface $\partial B = S^2 \supset \Gamma,$
the geometry near the dual edge $e^\vee$ of $e$ is a quadrilateral as
pictured here:
\begin{center}
\begin{tikzpicture}
\draw[very thick] (0,0)--(3,0)--(1.5,1.5)--(0,0);
\draw[very thick] (0,0)--(1.5,-1.5)--(3,0);
\node[above] at (1.5,0) {$e^\vee$};
\node[right] at (3,0) {$B$};
\node[above] at (1.5,1.5) {$C$};
\node[left] at (0,0) {$A$};
\node[below] at (1.5,-1.5) {$D$};
\draw[red,fill] (1.5,.75) circle (.1);
\draw[blue,fill] (1.5,-.75) circle (.1);

\end{tikzpicture}
\end{center}
Define $\widetilde{\triangle}$ as follows:  if the two faces
above are part of a tetrahedron $T$, then $\widetilde{\triangle} = 
\triangle \setminus T.$  Otherwise, let $T$ be the
tetrahedron with two faces as pictured above and the other
two $ACD$ and $BCD$, and set $\widetilde{\triangle} = \triangle
\cup T,$  In both cases, the geometry of the quadrilateral
$ABCD$ at the boundary of $\widetilde{\triangle}$ is
\begin{center}
\begin{tikzpicture}
\draw[very thick] (1.5,-1.5)--(1.5,1.5);
\draw[very thick] (0,0)--(1.5,-1.5)--(3,0)--(1.5,1.5)--(0,0);
\node[above right] at (1.5,0) {$\widetilde{e}{}^\vee$};
\node[right] at (3,0) {$B$};
\node[above] at (1.5,1.5) {$C$};
\node[left] at (0,0) {$A$};
\node[below] at (1.5,-1.5) {$D$};
\draw[blue,fill] (.85,0) circle (.1);
\draw[red,fill] (2.25,0) circle (.1);
\end{tikzpicture}
\end{center}
Call the foam so constructed $\bF_e$.
It remains to describe how to deform $\bF_e$ to $\bF'_{e,\pm}.$

Suppose $\widetilde{\triangle}$ is formed from $\triangle$
by adding a tetrahedron as in the proof of Proposition \ref{prop:hladdition}, and so $\bF_e$ is formed
from $\bF$ by attaching a Harvey-Lawson foam $\bF_{\text{HL}}$.  We
define $\bF'_e$ by extending $\bF'$ together with
a choice of one of the three possible smoothings of
$\bF_{\text{HL}}$ One of these three pairs the two
tangles with endpoints
at the centers of triangles $ABC$ and $ABD$
(pictured in red and blue)
with one another, creating a new short tangle component.
This is the disallowed smoothing.
The other two rotate pair these with the centers of the
two new triangles $ACD$ and $BCD$.  The matching
corresponding to the deformed foam of the positive mutation
$\bF'_{e,+}$ is shown
above.  $\bF'_{e,-}$ is defined similarly.

Now suppose otherwise that $\widetilde{\triangle}$ is formed
from $\triangle$ by deleting a tetrahedron $T$.
Then $\triangle$ was the result of a mutation $\mu_{\widetilde{e},\pm}$ of $\widetilde{\triangle}$
and only the inverse mutation $\mu_{e,\mp}$ is possible.
Since the case where $e$ bounds a single tangle component is not allowed,
the tangle components after deleting the tetrahedron are clear:  they
are truncations of the original tangle strands.

\end{proof}

Proposition \ref{prop:foammutation}
will allow us to transport foams
across mutations, along with phases, framings and cones.
This will allow us to connect open Gromov-Witten conjectures
for Lagrangian fillings related by allowed mutations
which have corresponding cones, phases and framings.

\begin{proposition}
\label{prop:mutationtransfer}
Let $\bF'$ and $\widetilde{\bF}'$ be deformed foams corresponding to
an allowed mutation $\Gamma\to \widetilde{\Gamma}$ of their boundaries. 
Then there is a canonical isomorphism
$$\varphi:\bZ^{E_\Gamma\cup A_{\bF'}}/\sim_{\bF'} \;\cong\; \bZ^{E_{\widetilde{\Gamma}}\cup A_{\widetilde{\bF}'}}/\sim_{\widetilde{{\bF}}'}$$
\end{proposition}
\begin{proof}
We can assume that $\widetilde{\Gamma}$ is obtained by
attaching a tetrahedron, as removal will give 
rise to the inverse isomorphism.  A local study near the
attachment will suffice to establish $\varphi.$  We label the relevant edges
and vertices as in the figure below, with $\Gamma$ indicated
by dashed lines.
\begin{center}
\begin{figure}[ht]
\begin{tikzpicture}[rotate around x=-90, rotate around x=5,scale=2]
\pgfmathsetmacro{\eps}{.1}
\coordinate (a) at (1,0,-1);
\coordinate (b) at (0,1,-1);
\coordinate (c) at (-1,0,-1);
\coordinate (d) at (0,-1,-1);
\coordinate (e) at (1/3,1/3,-1);
\coordinate (ebelow) at (1/3,1/3,-2);
\coordinate (emid) at (0,0,-1);
\coordinate (f) at (-1/3,-1/3,-1);
\coordinate (fbelow) at (-1/3,-1/3,-2);
\coordinate (A) at (1,0,1);
\coordinate (B) at (0,1,1);
\coordinate (C) at (-1,0,1);
\coordinate (D) at (0,-1,1);
\coordinate (E) at (1/3,-1/3,1);
\coordinate (Emid) at (0,0,1);
\coordinate (F) at (-1/3,1/3,1);
\coordinate (g) at (0,0,0);
\coordinate (g1) at (-\eps,0,0);
\coordinate (g2) at (\eps,0,0);
\node[below] at (emid) {$e$};
\node[above right] at (Emid) {$\widetilde{e}$};
\node[right] at (g2) {$\alpha$};
\node[right] at (a) {$a$};
\node[above] at (b) {$b$};
\node[right] at (A) {$\widetilde{a}$};
\node[right] at (B) {$\widetilde{b}$};
\node[left] at (c) {$c$};
\node[below] at (d) {$d$};
\node[left] at (C) {$\widetilde{c}$};
\node[left] at (D) {$\widetilde{d}$};
\node[above right] at (e) {$v$};
\node[below right] at (E) {$\widetilde{v}$};
\node[below left] at (f) {$w$};
\node[above left] at (F) {$\widetilde{w}$};
\draw[thick,green] (g1)--(g2);
\draw[thick,blue] (E)--(F);
\draw[thick,blue,dashed] (e)--(f);
\draw[thick,red] (e)--(ebelow);
\draw[thick,red] (f)--(fbelow);
\foreach \i in {B,C}
{\draw[thick,blue] (F)--(\i);}
\foreach \i in {A,D}
{\draw[thick,blue] (E)--(\i);}
\foreach \i in {a,b}
{\draw[thick,blue,dashed] (e)--(\i);}
\foreach \i in {c,d}
{\draw[thick,blue,dashed] (f)--(\i);}
\draw[thick,red] plot [smooth,tension=.5] coordinates {(e) (g2) (E)};
\draw[thick,red] plot [smooth,tension=.5] coordinates {(f) (g1) (F)};
\end{tikzpicture}
\end{figure}
\end{center}
Consider the $\binom{4}{2} = 6$ sheets of the deformed Harvey-Lawson foam,
after gluing to $\bF'$ and deforming.  They correspond to unordered pairs from among
the vertices $\{v,w,\widetilde{v},\widetilde{w}\}$.  
Write $f_{v,w}$ for the face determined by $v$ and $w$, and likewise for the
others.  
Let $\gamma$ be the arc of the Harvey-Lawson deformed foam, and write $\alpha = 
\sigma(f_{v,w},\gamma) \cdot \gamma = \pm \gamma$
for the signed contribution to the relation from $f_{v,w}$,
as defined in Definition \ref{def:facerelation}.
Now suppose the face relations
on $\bF'$ relate give $e + s_e \sim 0$, $a + s_a \sim 0$, and so on.
Let us list the unordered pairs along with the relations from the corresponding
glued face.
$$
\begin{array}{llll}
f_{v,w}:& \alpha + s_e&\Rightarrow&\alpha = e
\\ f_{\widetilde{v},\widetilde{w}}:& \widetilde{e} + \alpha  &\Rightarrow&\widetilde{e}=-\alpha = -e
\\ f_{v,\widetilde{v}}:& \widetilde{a} + s_a&\Rightarrow&\widetilde{a} = a
\\ f_{w,\widetilde{w}}:& \widetilde{c} + s_c&\Rightarrow&\widetilde{c} = c
\\ f_{v,\widetilde{w}}:& \widetilde{b} - \alpha + s_b&\Rightarrow&\widetilde{b}=b+e
\\ f_{\widetilde{v},w}:& \widetilde{d} - \alpha + s_d&\Rightarrow&\widetilde{d}=d+e
\end{array}
$$
This gives the positive mutation. 
The other allowed matching $v\leftrightarrow \widetilde{w}$ gives the negative mutation,
as follows from the interchange $\widetilde{v}\leftrightarrow\widetilde{w}.$  

\end{proof}

\begin{corollary}
\label{cor:mutationtransfer}
Suppose $\widetilde{\Gamma}$ is obtained by an allowed mutation of $\Gamma$.
Let $\nu:\bZ^{E_\Gamma}\to \bZ^{E_{\widetilde{\Gamma}}}$ be the corresponding isomorphism of edge lattices, respecting the
antisymmetric pairing.  Let $\varphi:\bZ^{E_\Gamma\cup A_{\bF'}}/\sim_{\bF'} \;\cong\; \bZ^{E_{\widetilde{\Gamma}}\cup A_{\widetilde{\bF}'}}/\sim_{\widetilde{{\bF}}'}$
be the isomorphism provided by Proposition \ref{prop:mutationtransfer} above.
Then under the isomorphisms of Equation \ref{eq:h1} and Proposition \ref{prop:hladdition}, the maps
$\nu$ and $\varphi$
intertwine $\tau:H_1(S_\Gamma)\to H_1(L)$
with $\widetilde{\tau}:H_1(S_{\widetilde{\Gamma}})\to H_1(\widetilde{L}).$
\end{corollary}
\begin{proof}
It only remains to note that $\nu$ respects the antisymmetric pairing of edges.
\end{proof}

We immediately obtain the following.
\begin{corollary}
\label{cor:seedtransfer}
The maps $\nu$ and $\varphi$ map phases, framings and cones to phases, framings and cones.
\end{corollary}

%% file: wavefunction.tex
\section{The wavefunction}
\label{sec:wavefunction}
\subsection{Construction of the wavefunction}
Suppose that $\mathbf{F}'$ is a deformed ideal foam obtained from the standard necklace foam by a sequence of admissible mutations, and $\mathbf{f}$ is a framing for $\mathbf{F}'$. As explained in Section~\ref{sec:foamsphasesandframings}, the pair $(\mathbf{F}',\mathbf{f})$ gives rise to a framed seed $\underline{\mathbf i}=\underline{\mathbf i}(\mathbf{F}',\mathbf{f})$. It is convenient to visualize the framed seed as a labelling of the edges of the cubic graph $\Gamma$ by monomials in the standard quantum torus generated over $\mathbb{Z}[q^{\pm}]$ by $\{U^{\pm}_i,V^{\pm}_i\}_{i=1,\ldots g}$.

In this section, we will show that there is a canonical wavefunction $\Psi_{\underline{\mathbf{i}}}\in\mathcal{K}$ associated to such a  framed seed, thereby providing a prediction for the generating function of all-genus open Gromov-Witten invariants of the corresponding Lagrangian $L_{\mathbf{F}'}\subset \mathbb{C}^3$.

We begin with the definition of $\Psi$ in the case of the standard necklace framed seed $\underline{\mathbf{i}}_{\mathrm{neck}}$. The corresponding foam gives rise to an exact Lagrangian filling of the Chekanov surface, so that by Stokes' theorem all its open Gromov-Witten invariants will be zero. We therefore take the wavefunction for the standard necklace to be $\Psi_{\mathrm{neck}}=1$.  Let us note that the necklace wavefunction depends only on the underlying deformed foam, and is completely independent of the choice of framing $\mathbf{f}$.


Now recall the sub-category $\mathbb{G}_{\mathrm{ad}}$ of the framed seed groupoid $\mathbb{G}$ whose morphisms are given by the admissible ones, and let $\mathbb{G}_{\mathrm{ad}}(\underline{\mathbf{i}}_{\mathrm{neck}})$ be the connected component of $\mathbb{G}_{\mathrm{ad}}$ containing the framed seed $\underline{\mathbf{i}}_{\mathrm{neck}}$. 

Part of the data of a representation of a quantum cluster variety (as defined in~\cite{FG2}) consists of a functor from the cluster modular groupoid to the category whose objects are Hilbert spaces, and whose morphisms are unitary equivalences. 

 By analogy, let us define an \emph{algebraic representation} of the admissible groupoid $\mathbb{G}_{\mathrm{ad}}(\underline{\mathbf{i}}_{\mathrm{neck}})$ to be a functor from $\mathbb{G}_{\mathrm{ad}}(\underline{\mathbf{i}}_{\mathrm{neck}})$ to the category $\mathrm{Vect}_{\mathbb{Q}(q)}$ of $\mathbb{Q}(q)$-vector spaces with morphisms given by $\mathbb{Q}(q)$-linear maps. 

The results established so far allow us to construct an algebraic representation of $\mathbb{G}_{\mathrm{ad}}(\underline{\mathbf{i}}_{\mathrm{neck}})$ as follows. Consider the functor which 
 assigns to each object $\underline{\mathbf{i}}$ of $\mathbb{G}_{\mathrm{ad}}$ the same vector space $\mathcal{K}$, and assigns to each arrow $a:\underline{\mathbf{i}}_{\mathrm{neck}}\rightarrow \underline{\mathbf{i}}$ in $\mathbb{G}_{\mathrm{ad}}$
the automorphism $\Phi_{\vec a}$ of $\mathcal{K}$ defined in Section~\ref{sec:groupoid-rep}.
 In this language, Lemma~\ref{lem:repwelldef} implies
   \begin{lemma}
   \label{lem:grpd-rep}
   The assignment 
   $$
  \underline{\mathbf{i}}\mapsto \mathcal{K}, \quad \vec{a}\mapsto\Phi_{\vec a} 
  $$ 
  defines an algebraic representation of the admissible groupoid $\mathbb{G}_{\mathrm{ad}}(\underline{\mathbf{i}}_{\mathrm{neck}})$.
 \end{lemma}
   
For each  object $\underline{\mathbf{i}}$ of $\mathbb{G}_{\mathrm{ad}}(\underline{\mathbf{i}}_{\mathrm{neck}})$, we now explain how to construct a canonical vector $\Psi_{\underline{\mathbf{i}}}\in\mathcal{O}\subset\mathcal{K}$ which we call the \emph{wavefunction} of the framed seed $\underline{\mathbf{i}}$. This vector is constructed as follows: choose an arbitrary path $\vec a:\underline{\mathbf{i}}_{\mathrm{neck}}\rightarrow \underline{\mathbf{i}}$ in $\mathbb{G}_{\mathrm{ad}}$. By Lemma~\ref{lem:grpd-rep}, the morphism $\vec a$ gives rise to an automorphism $\Phi_{\vec a}$ of $\mathcal{K}$, which we apply to $\Psi_{\mathrm{neck}}$ to produce a candidate for $\Psi_{\underline{\mathbf{i}}}$:
\begin{align}
\label{eq:wavefunction-def}
\Psi_{\underline{\mathbf{i}}}:= \Phi_{\vec a}\cdot \Psi_{\mathrm{neck}}.
\end{align}
What must be checked in order for this definition to make sense is that the wavefunction $\Psi_{\underline{\mathbf{i}}}$ depends only on the endpoint of the path $\vec a$ in the framed seeds groupoid. This path-independence is the content of the following Theorem.

\begin{theorem} 
\label{thm:wavefunction}
The map 
\begin{align*}
\Psi ~:~ \mathrm{Ob}\left(\mathbb{G}_{\mathrm{ad}}(\underline{\mathbf{i}}_{\mathrm{neck}})\right)& \longrightarrow \mathcal{O},\qquad
\underline{\mathbf{i}} \longmapsto  \Psi_{\underline{\mathbf{i}}}
\end{align*}
is well-defined, i.e is independent of the choice of path $\vec a:\underline{\mathbf{i}}_{\mathrm{neck}}\rightarrow \underline{\mathbf{i}}$ in ~\eqref{eq:wavefunction-def}. Moreover, if there exists such a path consisting entirely of primitive mutations, the wavefunction $\Psi_{\underline{\mathbf{i}}}$ satisfies the Ooguri-Vafa integrality constraint~\eqref{eq:OV-intro}.
\end{theorem}

\begin{proof}
The key observation is the following immediate consequence of Theorem~\ref{thm:d-mod}: if $\Psi_{\mathbf{i}}$ satisfies the face relations in framed seed $\underline{\mathbf{i}}$ and $\underline{\mathbf{i}}' = a(\underline{\mathbf{i}})$ where $a$ is an admissible mutation or framing shift, then $\Phi_a\cdot\Psi_{\mathbf{i}}$ satisfies the face relations for $\underline{\mathbf{i}}'$. Now suppose we have two sequences of admissible mutations and framing shifts $\vec{a}_1$ and $\vec{a}_2$ as in the statement of the theorem. Then it suffices to show that
\begin{equation}
\Phi_{\vec{a}_1}^{-1}\Phi_{\vec{a}_2}\cdot\Psi_{\mathrm{neck}}=\Psi_{\mathrm{neck}}.
\end{equation}
To this end, consider the framed seed $\underline{\mathbf{i}}' =  \vec{a}_1^{-1} \vec{a}_2(\underline{\mathbf{i}}_{\mathrm{neck}})$. Its underlying cubic graph $\Gamma'$ is the image of the original necklace graph $\Gamma_{\mathrm{neck}}$ under an element of the mapping class group of the $(g+3)$-times punctured sphere, and moreover the labelling of the edges of $\Gamma'$ by monomials in the $U_i, V_i$ induced by its phase and framing are identical to that in the standard necklace framed seed. In particular, the face relations for $\underline{\mathbf{i}}'$ and $\underline{\mathbf{i}}_{\mathrm{neck}}$ are identical, and from the binomial face relations corresponding to the beads we deduce that
$$
\left(1-V_i\right)\cdot \left(\Phi_{\vec{a}_1}^{-1}\Phi_{\vec{a}_2}\cdot\Psi_{\mathrm{neck}}\right)=0, \quad i=1,\ldots, g.
$$
It follows that $\Phi_{\vec{a}_1}^{-1}\Phi_{\vec{a}_2}\cdot\Psi_{\mathrm{neck}}=\Psi_{\mathrm{neck}}=1$, which completes the proof that the map $\Psi$ is well-defined. The Ooguri-Vafa integrality follows from Proposition~\ref{prop:integrality}, which is established in Section~\ref{sec:integrality}.
\end{proof}

\subsection{Examples of wavefunctions}
We now proceed to compute the wavefunction defined in the previous section in some fundamental examples.

\begin{example} 
\label{eg:canoe}
The calculation in Example~\ref{eg:intertwining} shows that the wavefunction $\Psi_{\underline{\mathbf{i}}_1}$ associated to the framed seed $\underline{\mathbf{i}}_1$ for the canoe graph shown in Figure~\ref{fig:alg-canoe} is given by 
$$
\Psi_{\underline{\mathbf{i}}_1} = (X;q^2)_{\infty}.
$$
It satisfies the $q$-difference equation
$$
(1 +UV - V)\Psi_{\underline{\mathbf{i}}_1}=0,
$$
which is a scalar multiple of the face relation $R'$ in~\eqref{eq:g1-facerel}.
As an exercise, let us compute the effect on the wavefunction of applying the framing shift operator $T_{-1}$, which we recall acts on $\mathcal{A}_{2g}$ by $U\mapsto q^{-1}UV^{-1}$. The resulting framed seed is illustrated in Figure~\ref{fig:alg-reframed-canoe}.
\begin{figure}[htbp] 
 \begin{minipage}{0.5\linewidth} 
  \centering 
        \begin{tikzpicture}[scale=.75]


\draw[gray, thick] (6,1.5) -- (6,-1.5);

\draw[gray, thick] (6,1.5) -- (1.5,0);
\draw[gray, thick] (6,-1.5) -- (1.5,0);

\draw[gray, thick] (6,1.5) -- (10.5,0);
\draw[gray, thick] (6,-1.5) -- (10.5,0);


\draw[gray, thick] (1.5, 0) .. controls  (0,5) and (12,5) .. (10.5,0);


\filldraw[red] (4,0) circle (2pt) node[anchor=west]{};
\filldraw[red] (8,0) circle (2pt) node[anchor=west]{};
\filldraw[red] (6,2.5) circle (2pt) node[anchor=west]{};
\filldraw[red] (6,-2.5) circle (2pt) node[anchor=west]{};


\node  at (5.1,0) {\small $qVU^{-1}$};
\node[text=blue]  at (6.3,0) {\small$3$};

\node  at (6,4) {\small $qVU^{-1}$};
\node[text=blue]  at (6,3.5) {\small$6$};

\node  at (3,1) {\small $-q^{-1}V^{-1}$};
\node[text=blue]  at (4.5,1.3) {\small$1$};
\node  at (3,-1) {\small $-q^{-1}U$};
\node[text=blue]  at (4.5,-1.3) {\small$2$};

\node  at (9,1.1) {\small $-q^{-1}U$};

\node[text=blue]  at (7.5,1.3) {\small$4$};
\node  at (9,-1) {\small $-q^{-1}V^{-1}$};

\node[text=blue]  at (7.5,-1.3) {\small$5$};

\end{tikzpicture}
  \caption{The  framed seed $\underline{\mathbf{i}}_2=(\sigma_{-1}\circ T_{-1})(\underline{\mathbf{i}}_1)$ for the canoe graph. } 
  \label{fig:alg-reframed-canoe} 
 \end{minipage}%
 \begin{minipage}{0.5\linewidth} 
  \centering 
\begin{tikzpicture}[scale=.72]
\draw[gray, thick] (3,0) circle (1.5cm);
\draw[gray, thick] (4.5,0) -- (7.5,0);
\draw[gray, thick] (9,0) circle (1.5cm);

\draw[gray, thick] (1.5, 0) .. controls  (0,5) and (12,5) .. (10.5,0);


\filldraw[red] (3,0) circle (2pt) node[anchor=west]{};
\filldraw[red] (9,0) circle (2pt) node[anchor=west]{};
\filldraw[red] (6,2) circle (2pt) node[anchor=west]{};
\filldraw[red] (6,-2) circle (2pt) node[anchor=west]{};


\node at (6,.3) {\small$-q^{-1}U$};
\node[text=blue]  at (6,-.3) {\small$2$};

\node (b) at (6,4) {\small $-qU^{-1}$};
\node[text=blue]  at (6,3.5) {\small$4$};

\node  at (3,1.8) {\small$-q^{-1}V^{-1}$};
\node[text=blue] at (3,1.2) {\small$1$};
\node[text=blue]  at (3,-1.2) {\small$6$};
\node (e) at (3,-1.8) {\small$-q^{-1}V$};

\node (d) at (9,1.8) {\small$-q^{-1}V$};
\node[text=blue] at (9,1.2) {\small$3$};
\node[text=blue]  at (9,-1.2) {\small$5$};
\node (e) at (9,-1.8) {\small$-q^{-1}V^{-1}$};
\end{tikzpicture}
  \caption{The  framed seed $\mu_4^+(\underline{\mathbf{i}}_2)\simeq \underline{\mathbf{i}}_0$. } 
  \label{fig:alg-finalneck} 
 \end{minipage} 
\end{figure}

\begin{lemma}
\label{lem:reframed-canoe}
We have 
\begin{align}
\label{eq:g1-canoe-wav}
    \Psi_{\underline{\mathbf{i}}_2}&=(\sigma_{-1} \circ T_{-1})\cdot (X;q^2)_{\infty}\\
    &=\nonumber (X;q^2)_{\infty}^{-1}.
\end{align}
\end{lemma}
\begin{proof}
Since
$$
 (\sigma_{-1} \circ T_{-1})(1+UV-V) = (1-U-V) (\sigma_{-1} \circ T_{-1}),
$$
the Lemma follows by observing that both sides satisfy the $q$-difference equation
$$
(1-U - V)\Psi = 0,
$$
which is easily seen to have a unique formal power series solution of the form $\Psi\in 1+\mathfrak{m}\in\mathcal{O}$.
\end{proof}
Now observe that applying to the framed seed $\underline{\mathbf{i}}_2$ the positive mutation at edge 4 returns us to the framed seed shown in Figure~\ref{fig:alg-finalneck}, which coincides with the standard necklace $\underline{\mathbf{i}}_0$ up to a permutation of the numbering of its edges. Hence we have a loop in the framed seed groupoid
\begin{equation}
\label{eq:loop-g1}
\begin{tikzcd}[row sep=huge]
 &
\underline{\mathbf{i}}_0 \arrow[dl,swap,"\mu_3^+"]  &
\\
\underline{\mathbf{i}}_1 \arrow[rr,"\sigma_{-1} \circ T_{-1}"]& & \underline{\mathbf{i}}_2 \arrow[ul,swap,"\mu_4^+"]
\end{tikzcd}
\end{equation}
and we indeed see that 
\begin{align*}
\Psi_{\mu_4^+(\underline{\mathbf{i}}_2)} &= \Phi(-q^{-1}U)\cdot (X;q^2)^{-1}_{\infty}\\
&= (X;q^2)_{\infty}\cdot (X;q^2)^{-1}_{\infty}\\
&=1,
\end{align*}
in accordance with Theorem~\ref{thm:wavefunction}.

More generally, we can consider the framed seed $\underline{\mathbf{i}}_{\mathrm{canoe}}$ obtained from the standard genus $g$ necklace framed seed $\underline{\mathbf{i}}_{\mathrm{neck}}$ by performing positive mutations at all $g$ beads labelled $-q^{-1}U_j, ~j=1,\ldots g$ under the framing isomorphism. The corresponding wavefunction is then
\begin{equation}
\label{eq:canoewavefunction}
\Psi_{\mathbf{i}_{\mathrm{canoe}}} = \prod_{i=1}^g(X_i; q^2)_\infty,
\end{equation}
which is annihilated by the left ideal in $\mathcal{D}_{2g}$ generated by
\begin{equation}
\label{eq:genus-g-necklace}
{R}_i = 1 + U_iV_i -  V_i,  \hskip 15mm i =1, \ldots, g.
\end{equation}
Let us write $\Psi_{\mathbf{i}_{\mathrm{canoe}}^{(1)}}$ for the wavefunction obtained by applying the operator $\sigma_{(-1,\ldots,-1)}\circ T_{-I_g}$ to $\Psi_{\mathbf{i}_{\mathrm{canoe}}}$, where $I_g$ is the $g\times g$ identity matrix. Then we again have
\begin{equation}
\label{super.pants}
\Psi_{\mathbf{i}_{\mathrm{canoe}}^{(1)}} = \prod_{i=1}^g(X_i; q^2)^{-1}_\infty
\end{equation}

\begin{lemma} 
\label{lem:canoepowerseries}
The explicit power series of the wavefunction \eqref{super.pants} is 
\[
\Psi_{\mathbf{i}_{\mathrm{canoe}}^{(1)}} = \sum _{{\bf v}\in \mathbb{Z}_{\geq 0}^g} \frac{1}{(q^2)_{\bf v}} X^{\bf v},  \hskip 1cm \mbox{where } (q^2)_{\bf v}=\prod_{i=1}^g \prod_{k=1}^{v_i}(1-q^{2k}). 
\]
\end{lemma}
\begin{proof} Set
\[
\Psi := \sum _{{\bf v}\in \mathbb{Z}_{\geq 0}^g} \frac{C_{\bf v}(q)}{(q^2)_{\bf v}} X^{\bf v},  \hskip 1cm \mbox{where } C_{0}(q)=1
\]
Let $e_i  \in \mathbb{Z}^g$ be the $i^{th}$ unit vector. At the level of the coefficients of $X^{\bf v}$, the equation $(1-U_i-V_i)\Psi=0$ is equivalent to the recurrence
\[
\frac{C_{\bf v}(q)}{(q^2)_{\bf v}} -   \frac{C_{{\bf v}-e_i}(q)}{(q^2)_{{\bf v}-e_i}}  - \frac{C_{\bf v}(q)}{(q^2)_{\bf v}} q^{2v_i} =0
\]
Note that $(q^2)_{\bf v}= (q^2)_{{\bf v}-e_i} (1-q^{2v_i})$. Therefore we have
\[
C_{\bf v}(q) = C_{{\bf v}-e_i}(q)=\ldots = C_0=1.
\]
\end{proof}
\end{example}
More generally, given a $g\times g$ integer symmetric matrix $A$,
let us consider the framed seed $\mathbf{i}_{\mathrm{canoe}}^{(A)}$ obtained by applying $\sigma_{(-1,\ldots,-1)}\circ T_{-A}$ to $\mathbf{i}_{\mathrm{canoe}}$. Then by~\eqref{eq:frame-change-action} we have
\begin{align}
    \label{eq:canoe-general-framing}
    \Psi_{\mathbf{i}_{\mathrm{canoe}}^{(A)}} = \sum _{{\bf v}\in \mathbb{Z}_{\geq 0}^g} \frac{1}{(q^2)_{\bf v}} q^{{\bf v}^t{\bf v}-{\bf v}^tA{\bf v}}X^{\bf v},  \hskip 1cm \mbox{where } (q^2)_{\bf v}=\prod_{i=1}^g \prod_{k=1}^{v_i}(1-q^{2k})
\end{align}

\begin{example}[Non-existence of algebraic wavefunctions]
\label{eg:non-algebraic}
Algebraic wavefunctions may not exist for framed seeds which cannot be obtained from the standard necklace by a sequence of admissible mutations. A simple counterexample is given by the framed seed obtained from the standard $g=1$ necklace by performing positive mutations at both of its strands to produce another necklace graph. The arguments of the corresponding quantum dilogarithms are $-q^{-1}U$ and $-qU^{-1}$, only the first of which corresponds to an admissible mutation. The face relation associated to the bead of the resulting framed seed imposes the difference equation $(1+VU)\Psi=0$, and it is easy to see this admits no nonzero solutions in the ring $\mathbb{Q}((q))((X))$ (or for that matter in the opposite completion $\mathbb{Q}((q))((X^{-1}))$.
\end{example}

We conclude this section with an example from \cite{AENV} concerning the unknot conormal after the
conifold transition.  

\begin{example}[Partition function for unknot conormal \cite{AENV}] 
\label{ex:aenv}
Set
\[
\Psi(X)= (X; q^2)_{\infty}^{-1}(QX; q^2)_\infty,
\]
where the closed-string parameter $Q$ is a formal variable commuting with all the other variables.
Then $\Psi(X)$ is annihilated by 
\[
\cL=(1-U)-(1-QU)V= 1-U-V+ QUV
\]
\begin{lemma} We have
\[
\Psi(X):=\sum_{k\geq 0} \frac{(Q; q^2)_k}{(q^2;q^2)_k} X^k.
\]
\end{lemma}
\begin{proof} Set
\begin{align}
\label{eq:wav-conifold}
\Psi(X):=\sum_{k\geq 0} C_k(Q, q^2)X^k, \hskip 2cm  \mbox{where } C_0=1
\end{align}
By computing the coefficients of $\cL \cdot \Psi$, we get
\[
C_k - C_{k-1} - q^{2k} C_k + Qq^{2(k-1)} C_{k-1} = (1-q^{2k}) C_k - (1- Qq^{2(k-1)}) C_{k-1}=0
\]
Therefore
\[
C_k =\frac{\prod_{i=1}^{k} (1-Qq^{2(i-1)})}{\prod_{i=1}^k(1-q^{2i})}=  \frac{(Q; q^2)_k}{(q^2;q^2)_k}.
\]
\end{proof}
Note that under the specialization $Q=1$ of the closed string parameter we recover the necklace wavefunction $\Psi_{\mathrm{neck}}=1$, while the specialization $Q=0$ delivers the wavefunction $\Psi_{\underline{\mathbf{i}}_2}$ in~\eqref{eq:g1-canoe-wav} associated to the framed seed~\eqref{fig:alg-reframed-canoe} for the canoe graph. Hence the closed-string parameter $Q$ describes an interpolation between these two framed seeds. 
\end{example}

\subsection{Open Gromov-Witten Conjectures}

We can now propose an interpretation of the wavefunction of a geometric seed: 
it is the generating function of open Gromov-Witten invariants of the Lagrangian
filling defined by the deformed foam.  To be more precise, we recall the geometric
framework.

Let $\underline{\bf i} \in {\mathbb G}_{\mathrm{ad}}$, meaning there is a
path $\vec{a}: \underline{\bf i}_{\mathrm neck} \to \underline{\bf i}$
in the admissible framed seed groupoid.
By Theorem \ref{thm:wavefunction}, there is a well-defined
wavefunction $\Psi_{\underline{\bf i}} = \Phi_{\vec{a}}\cdot
\Psi_{{\underline{\bf i}_{\mathrm{neck}}}}
= \Phi_{\vec{a}}\cdot 1.$

The framing of $\underline{\bf i}$ has geometric content.
Recall from Example \ref{ex:neckseed} that $\underline{\bf i}_{\mathrm{neck}}$ is canonical.
By Proposition \ref{prop:mutationtransfer} and especially
Corollary \ref{cor:seedtransfer}, we learn $\underline{\bf i} = \vec{a}\cdot \underline{\bf i}_{\mathrm{neck}}$ is a geometric seed, i.e.~has a geometric phase, as well as a framing and cone.  That is, there is
a corresponding cubic
graph $\Gamma$, deformed foam filling $\bF'$, and Lagrangian $L$, along with phase
$K = {\mathrm Ker}(\tau: H_1(S_\Gamma)\to H_1(L))$ and transverse isotropic
framing $F\subset H_1(S_\Gamma),$ as well as a cone $C \subset F.$ 
We choose a basis $e_i,$ $i = 1,...,g$ for $C$. 
The sequence $0\to K\to H_1(S_\Gamma)\to \pi(F) \to 0$ and basis $e_i$
then defines a framing
for ${\bf i}$ in the sense of Section \ref{sec:polarizationsandframings}.

The geometric seed identifies the quantum torus $\cT^q_{\bf i}$ with
the quantization of the symplectic lattice $H^1(S_\Gamma)$ endowed with
its intersection form.
In particular, a monomial $X^d = \prod_{i=1}^g X_i^{d_i}$
in the ring of power series $\bC[[\{X_i\}]]$
has exponent $d$ lying in $\bZ_{\geq 0}^g \cong C \subset H_1(L).$
Each such $d$ determines an open Gromov-Witten problem of counting holomorphic
maps from Riemann surfaces with one boundary component 
mapping to the pair $(\bC^3,L)$, such that the image of the boundary lies
in homology class $d$.  Such open Gromov-Witten problems depend on additional
data known as a framing. While there is not yet a rigorous definition
of these open Gromov-Witten invariants, it is anticipated that it will involve
framings as constructed here,
generalizing the well-studied cases of Aganagic-Vafa branes
\cite{AKV,KL,L,FL}, fixed points of anti-symplectic
involutions and rational cohomology spheres \cite{ST}.\footnote{We thank Jake Solomon and Sara Tukachinsky
for relaying their expectations for the more general class of Lagrangians considered in this paper.  See also \cite{I} for a more general
definition of framings.}

We then conjecture that
the wavefunction $\Psi_{\underline{\bf i}} \in \bC[q,q^{-1}][[\{X_i\}]]$
is the all-genus generating function of open Gromov-Witten invariants 
and obeys Ooguri-Vafa integrality, which expresses
the invariants in terms of the quantum dilogarithm $\Phi(z) = \prod_{n\geq 0}(1+q^{2n+1}z)^{-1}$.

\begin{conjecture}
\label{conj:ogw}
Let $\underline{\bf i}\in {\mathbb G}_{\mathrm{ad}}$ be a framed seed
with wavefunction 
$\Psi_{\underline{\bf i}}.$  Write $A$ for the framing and $L$ for the Lagrangian
of the deformed foam.  Then
$$\Psi_{\underline{\bf i}} = \prod_{d\in \bZ_{\geq 0}\setminus \{0\}}\prod_{s\in\mathbb{Z}}
\Phi(X^d (-q)^s)^{n^{(A)}_{d,s}},$$
with $n^{(A)}_{d,s}\in \bZ$ the Ooguri-Vafa invariants.
\end{conjecture}

\begin{remark}
The Ooguri-Vafa invariants are related to open Gromov-Witten invariants as follows.
Write $q^2 = e^\lambda$ and expand $\Psi_{\underline{\bf i}}$
as a power series in $\lambda$ (and the $X_i$).
Then the coefficient of $X^d \lambda^h$ is the genus-$h$ open Gromov-Witten invariant
of $L$ in framing $A$, in class $d\in H_1(L) \cong H_2(\bC^3,L)$.
See \cite[Sections 2 and 4]{Za} for further discussion of these variables.
\end{remark}


\begin{remark}
Recall $\Phi(z)\sim e^{{\rm Li}_2(z)/\lambda}$ as $\lambda \to 0$.
Conjecture \ref{conj:ogw} therefore reduces to the conjecture of \cite{TZ}
for disk invariants, described in the Introduction in Section \ref{sec:intro-ogw}.
More specifically, writing $\Psi_{\underline{\bf i}} \sim e^{W_{\underline{\bf i}}/\lambda}$,
in the semiclassical limit, meaning $W$ is a local
potential for the Lagrangian subspace $\cM_\Gamma \subset \cP_\Gamma.$
\end{remark}

\begin{remark}
In the next section, we provide evidence for the conjecture
by arguing that the wavefunctions $\Psi_{\underline{\bf i}}$ obey integrality.
The Harvey-Lawson brane in $\bC^3$ with its various framings,
as studied in \cite[Section 6.1]{AKV}, gives further evidence.
This example enjoys a $U(1)$ symmetry, permitting localization techniques for
open Gromov-Witten calculations \cite{KL}, while the Lagrangians for
for cubic graphs $\Gamma$ generally do not.  Further tests of the
conjecture must therefore await rigorous defitions of open Gromov-Witten invariants
and the development of new techniques.
\end{remark}

\subsection{Integrality of the wavefunction}
\label{sec:integrality}
In this section we will complete the proof of Theorem~\ref{thm:wavefunction} by showing that the wavefunctions $\Psi_{\underline{\bf i}}$ satisfy Ooguri-Vafa integrality.

The wavefunction $\Psi$ constructed in the previous section is an element of the commutative local ring $\mathcal{O}_{\mathbb{Q}}:=\mathbb{Q}\otimes_\mathbb{Z}\mathcal{O}$ of formal power series in $X_1,\ldots, X_g$ with coefficients in the field $\mathbb{Q}((q))$. Let $\mathfrak{m}$ be the unique maximal ideal in the ring $\mathcal{O}_{\mathbb{Q}}$.
By considering the quotients $\mathcal{O}_{\mathbb{Q}}/\mathfrak{m}^k$, it is easy to show that every $F\in 1+\mathfrak{m}$ admits a unique factorization
\begin{equation}
\label{infty.product}
F= \prod_{{\bf v}\in \mathbb{Z}_{\geq 0}^g- \{0\}}\,  \prod_{k\in \mathbb{Z}} \left((-q)^{k}X^{\bf v}; q^2\right)_{\infty}^{c_{\bf v},k}, \hskip 15mm c_{{\bf v}, k}\in \mathbb{Q}.
\end{equation}
The coefficients $c_{{\bf v},k}$ for each ${\bf v}\in  \mathbb{Z}_{\geq 0}^g\setminus \{0\}$ can be packaged in a Laurent series 
\[
P_{F, {\bf v}}(t):= \sum_{k\in \mathbb{Z}} c_{{\bf v}, k} t^k \in \mathbb{Q}((t)).
\]
Following \cite[\S6.1]{KS},  a series $F\in 1+\mathfrak{m}$ is called {\it admissible} if  the $P_{F, {\bf v}}(t)$ are Laurent polynomials with integral coefficients for all ${\bf v}$.
\vskip 2mm

Recall the logarithm
\[
\log:  ~~1+\mathfrak{m} \lra \mathfrak{m}, \hskip 15mm \log (1+f) := \sum_{k\geq 1}\frac{(-1)^k f}{k}. 
\]
\begin{lemma}
For each admissible series $F$, we have
\[
\lim_{q\rightarrow -1} \, (q^2-1)\, \log \, F  = \sum_{{\bf v}\in \mathbb{Z}_{\geq 0}^n-\{0\}} P_{F, {\bf v}}(1) {\rm Li}_2(X^{\bf v}).
\]
\end{lemma}
\begin{proof} By the third formula of \eqref{Pochhammer.infty}, we have 
\begin{equation}
\label{log}
\lim_{q\rightarrow -1} \, (q^2-1)\, \log \,(x; q^2)_\infty= \lim_{q\rightarrow -1} \sum_{k=1}^\infty \frac{(q^2-1)x^k}{k(q^{2k}-1)}=\sum_{k=1}^\infty \frac{x^k}{k^2}= {\rm Li}_2(x).  
\end{equation}
The rest is clear.
\end{proof}
The property of admissibility is clearly preserved under the action of the coordinate rescaling automorphisms $\sigma_{\delta}$ introduced in~\eqref{eq:sigma-transf}. 
In~\cite{KS}, Kontsevich and Soibelman proved that admissibility is also preserved under another, much less trivial family of automorphisms: the changes of framing.
\begin{theorem}[{\cite[Th.6.1]{KS}}] 
\label{admissible.series}
A power series $F\in \mathcal{O}$ is admissible if and only if $T_\Omega\cdot F$ is admissible for all integral symmetric matrices $\Omega$.
\end{theorem}

The integrality of the wavefunction defined in Theorem~\ref{thm:wavefunction} now follows from the following Proposition.
\begin{proposition}
\label{prop:integrality}
Suppose that the mutation $a$ is both admissible and primitive, and that $F\in\mathcal{O}$ is an admissible formal power series. Then the power series $\Phi_a\cdot F$ is also admissible. 
\end{proposition}
\begin{proof}
Since the exponent vector $\mathbf{m}$ in ~\eqref{eq:exp-vector} is primitive we can choose a basis for $\mathbb{Z}^g$ containing $\mathbf{m}$ as one of its elements. Hence ({\it cf.} Remark~\ref{rmk:dual-basis-transformation}) we may reduce to proving the Lemma in the case that $
\iota_\mathbf{f}(e_k) = (-q)^r\prod_{j=1}^g U_1V_j^{n_j}.
$
In this case, let $\Omega$ be any symmetric matrix whose first column is $-\mathbf{n} =(-n_1,\ldots, -n_g) $. By Theorem~\ref{admissible.series}, it suffices to show $T_\Omega\cdot \Phi_a\cdot F$ is admissible, so we compute
\begin{align*}
   T_\Omega\cdot \Phi_a\cdot F&=  
   T_\Omega\cdot\Phi\left(\iota_\mathbf{f}(e_k)\right)^{\epsilon_k}\cdot  F \\
   &=\Phi\left((-q)^r U_1\right)^{\epsilon_k}\cdot T_\Omega\cdot F 
\end{align*}
But $T_\Omega\cdot F$ is admissible by Theorem~\ref{admissible.series}, and so is $\Phi\left((-q)^r U_1\right)^{\epsilon_k}$ by the last formula in~\eqref{Pochhammer.infty}. Since the product of two admissible series is clearly admissible, this  implies that $\Phi_a\cdot F$ is also admissible, thereby proving the Lemma.
\end{proof}

%% file: analytic.tex
\section{Towards an analytic wavefunction}
\label{sec:analyticaspects}

In this section we will discuss the problem of promoting the algebraic construction of the wavefunction from Section~\ref{sec:wavefunction} to an analytic one. Doing this in general would necessitate extending the theory of representations of quantum cluster varieties beyond the ``principal series'', a task we do not take on in the present work. Nonetheless, we will present several examples which we believe provide nontrivial evidence for the existence of a well-defined analytic wavefunction associated to a smoothed ideal foam.   

Let us first recall some of the elements of the theory of unitary representations of quantum cluster varieties as developed in~\cite{FG2}. A representation of a quantum cluster variety is, by definition, a functor $\mathcal{G}_{\mathscr{X}}\rightarrow\mathrm{Hilb}$ from the cluster modular groupoid to the category of Hilbert spaces with morphisms given by unitary isomorphisms. The representations constructed by Fock and Goncharov depend on a quantization parameter $\hbar\in\mathbb{R}$. To each object $\mathbf{i}$ of $\mathcal{G}_{\mathscr{X}}$ is associated a pair of quantum tori $\mathcal{T}^q$ and $\mathcal{T}^{q^\vee}$, generated respectively by $\{e^{2\pi\hbar x_k}\}$ and $\{e^{2\pi\hbar^{-1} x_k}\}$, where $\{x_k\}$ are the logarithmic cluster variables associated to the seed $\mathbf{i}$. For each $\mathbf{i}$, the generators of these quantum tori act by unbounded, self-adjoint operators in the Hilbert space $\mathcal{H}_{\mathbf{i}}$. The latter space comes equipped with the additional data of a dense subspace $\mathcal{S}_{\mathbf{i}}$, the Fock-Goncharov Schwartz space, defined to be the maximal joint domain of the algebras $\mathbb{L}^q_{\mathscr{X}},\mathbb{L}^{q^\vee}_{\mathscr{X}}$. The unitary isomorphism $\mathcal{K}_{\mathbf{i}\rightarrow\mathbf{i}'}$ corresponding to an arrow $a:\mathbf{i}\rightarrow\mathbf{i}'$ in the cluster modular groupoid preserves the corresponding Schwartz spaces, where it intertwines the action of $\mathcal{G}_{\mathscr{X}}$ on $\mathbb{L}^q_{\mathscr{X}}$ by cluster transformations. 

When the skew form on $\Lambda$ has a nontrivial kernel, the Fock-Goncharov unitary representations of the quantum cluster variety are labelled by \emph{central characters} $\lambda\in\mathrm{Hom} (\mathrm{ker}(\epsilon),\mathbb{R})$, and thus can be thought of as a kind of principal series. The reality condition is required to ensure that all elements of the underlying Heisenberg algebra act by self-adjoint operators. This self-adjointness is crucial for the entire construction: indeed, it guarantees that for each logarithmic cluster variable $x_k$ its noncompact quantum dilogarithm $\varphi(x_k)$ defines a unitary automorphism of $\mathcal{H}_{\mathbf{i}}$, which forms the key ingredient in defining the intertwiner $\mathcal{K}_{\mathbf{i}\rightarrow\mathbf{i}'}$.

In the context of moduli spaces of framed local systems on surfaces with punctures, recall that the central characters parametrize the eigenvalues of the local system's monodromy around the punctures. As we have seen in Section~\ref{sec:quantization}, however, the quantization~\eqref{quantum.face.0} of the defining constraints for $\mathcal{P}_g$, which impose that the monodromy around each puncture be unipotent, forces a sum of logarithmic cluster variables to act by a pure imaginary scalar, a constraint which cannot be satisfied if each such variable acts by a self-adjoint operator.  

Thus we cannot appeal to the standard theory of principal series in order quantize the chromatic Lagrangian -- a new kind of representation of the quantum cluster variety is required.  Although we do not currently know how such representations should be defined, let us sketch out some features we would desire of them in order to define an analytic wavefunction. 

To a framed seed $\underline{\mathbf{i}}(\mathbf{F}')$ with underlying deformed foam $\mathbf{F}'$, we would like to associate a space of meromorphic functions $\mathcal{V}_{\underline{\mathbf{i}}}$ in $g$ variables $z_1,\ldots, z_g$, defined by appropriate conditions on their asymptotic behavior along with the possible locations of their poles. The framed seed determines a natural action of $\mathbb{L}^q_{\mathscr{X}}$ by $q$-difference operators on the space of all meromorphic functions on $\mathbb{C}^g$, and the subspace $\mathcal{V}_{\underline{\mathbf{i}}}$ should be preserved under this action. 

To each admissible mutation or framing shift $a:\underline{\mathbf{i}}_1\rightarrow \underline{\mathbf{i}}_2$, there should correspond an isomorphism between the spaces $\mathcal{V}_{\underline{\mathbf{i}}_1},\mathcal{V}_{\underline{\mathbf{i}}_2}$. These isomorphisms should again intertwine the action of $\mathbb{L}^q_{\mathscr{X}}$, and their composites corresponding to trivial cluster transformations should act by scalar multiples of the identity. 

Given a representation of the quantum cluster variety $\mathcal{P}_{g}$ in this sense, one could then attempt to define the wavefunction associated to a framed seed obtained from the standard necklace by a sequence of admissible mutations as in~\eqref{eq:wavefunction-def}. To verify that this prescription is indeed well-defined would amount to showing that the action of the mapping class group of the punctured sphere (which is generated by the half Dehn twists around pairs of punctures) fixes the necklace wavefunction $\Psi_{\mathrm{neck}}=1$. 

Let us note  that there is another regime for $\hbar$ which is nicely compatible with the analytic properties of the noncompact quantum dilogarithm -- namely, when $|\hbar|=1$. In this regime, Faddeev \cite{F} has constructed discrete series-type representations of the modular double of $U_{q}(\mathfrak{sl}_2)$ whose central characters also correspond to a sum of logarithmic cluster variables acting by a pure imaginary scalar.
Thus the regime $|\hbar|=1$ may in fact be the most suitable one in which to try to carry out the construction of such representations of quantum cluster varieties associated to punctured surfaces.

In the following subsections, we present some explicit calculations in $g=1,2$ which indicate how one might try to define the action of admissible mutations and framing shifts in the non-unitary case, and provide examples of candidate analytic wavefunctions.

\subsection{Analytic wavefunctions for \texorpdfstring{$g=1$}{g=1}}
\label{subsec:g1-analytic}
We begin with the standard necklace framed seed
$\underline{\mathbf{i}}_0$ for $g=1$ as shown in Figure~\ref{fig:std-g1-necklace}. In the analytic setting, it is natural to regard a framed seed as associating to each edge a Heisenberg algebra element, which can then be  exponentiated to yield elements of either of the two modular dual quantum tori. 

\begin{figure}[htbp] 
 \begin{minipage}{0.5\linewidth} 
  \centering 
        \begin{tikzpicture}[scale=.7]


\draw[gray, thick] (3,0) circle (1.5cm);
\draw[gray, thick] (4.5,0) -- (7.5,0);
\draw[gray, thick] (9,0) circle (1.5cm);

\draw[gray, thick] (1.5, 0) .. controls  (0,5) and (12,5) .. (10.5,0);


\filldraw[red] (3,0) circle (2pt) node[anchor=west]{};
\filldraw[red] (9,0) circle (2pt) node[anchor=west]{};
\filldraw[red] (6,2) circle (2pt) node[anchor=west]{};
\filldraw[red] (6,-2) circle (2pt) node[anchor=west]{};


\node at (6,.3) {\small$u-c_\hbar$};
\node[text=blue]  at (6,-.3) {\small$3$};

\node (b) at (6,4) {\small $-u+c_\hbar$};
\node[text=blue]  at (6,3.5) {\small$6$};

\node  at (3,1.8) {$-v-c_\hbar$};
\node[text=blue] at (3,1.2) {\small$1$};
\node[text=blue]  at (3,-1.2) {\small$2$};
\node (e) at (3,-1.8) {$v-c_\hbar$};

\node (d) at (9,1.8) {$v-c_\hbar$};
\node[text=blue] at (9,1.2) {\small$4$};
\node[text=blue]  at (9,-1.2) {\small$5$};
\node (e) at (9,-1.8) {$-v-c_\hbar$};

\end{tikzpicture} 
  \caption{The standard necklace framed seed $\underline{\mathbf{i}}_0$ for $g=1$} 
  \label{fig:std-g1-necklace} 
 \end{minipage}%
 \begin{minipage}{0.5\linewidth} 
  \centering 
        \begin{tikzpicture}[scale=.75]


\draw[gray, thick] (6,1.5) -- (6,-1.5);

\draw[gray, thick] (6,1.5) -- (1.5,0);
\draw[gray, thick] (6,-1.5) -- (1.5,0);

\draw[gray, thick] (6,1.5) -- (10.5,0);
\draw[gray, thick] (6,-1.5) -- (10.5,0);


\draw[gray, thick] (1.5, 0) .. controls  (0,5) and (12,5) .. (10.5,0);


\filldraw[red] (4,0) circle (2pt) node[anchor=west]{};
\filldraw[red] (8,0) circle (2pt) node[anchor=west]{};
\filldraw[red] (6,2.5) circle (2pt) node[anchor=west]{};
\filldraw[red] (6,-2.5) circle (2pt) node[anchor=west]{};


\node  at (5.1,0) {\small $-u+c_\hbar$};
\node[text=blue]  at (6.3,0) {\small$3$};

\node  at (6,4) {\small $-u+c_\hbar$};
\node[text=blue]  at (6,3.5) {\small$6$};

\node  at (3,1) {\small $-v-c_\hbar$};
\node[text=blue]  at (4.5,1.3) {\small$1$};
\node  at (3,-1) {\small $u+v-2c_\hbar$};
\node[text=blue]  at (4.5,-1.3) {\small$2$};

\node  at (9,1.1) {\small $u+v-2c_\hbar$};

\node[text=blue]  at (7.5,1.3) {\small$4$};
\node  at (9,-1) {\small $-v-c_\hbar$};

\node[text=blue]  at (7.5,-1.3) {\small$5$};

\end{tikzpicture}
  \caption{The  framed seed $\underline{\mathbf{i}}_1=\mu_3^+(\underline{\mathbf{i}}_0)$ for the canoe graph. } 
  \label{fig:std-canoe} 
 \end{minipage} 
\end{figure}
Hence in Figure~\ref{fig:std-g1-necklace}, we have decorated the edges of the cubic graph by the Heisenberg algebra elements corresponding to the logarithmic cluster variables, and we have set (see Appendix~\ref{sec:qdl-appendix}) 
$$
c_\hbar := \frac{i}{2}\left(\hbar+\hbar^{-1}\right) \in i\mathbb{R}.
$$
The pair of modular dual quantum torus elements corresponding to edge $e_3$, for example, are given by
$$
X_3 \mapsto e^{2\pi \hbar(u-c_\hbar)} = -q^{-1}e^{2\pi\hbar u}, \qquad X_3^\vee \mapsto e^{2\pi \hbar^{-1}(u-c_\hbar)} = -(q^\vee)^{-1}e^{2\pi\hbar^{-1} u}.
$$
Now consider the following loop in the framed seed groupoid. First, observe that the positive mutation $\mu_3^+$ at edge 3 yields the canoe framed seed $\underline{\mathbf{i}}_1$ shown in Figure~\ref{fig:std-canoe}. Performing the change of framing $\sigma$ conjugating all Heisenberg algebra elements by $e^{\pi i(v-c_\hbar)^2 }$, thereby effecting the shift $\sigma: u\mapsto u-v+c_\hbar$, we arrive at the framed seed $\underline{\mathbf{i}}_2=\sigma(\underline{\mathbf{i}}_1)$ shown in Figure~\ref{fig:reframed-canoe}.

\begin{figure}[htbp] 
 \begin{minipage}{0.5\linewidth} 
  \centering 
\begin{tikzpicture}[scale=.72]

\draw[gray, thick] (6,1.5) -- (6,-1.5);

\draw[gray, thick] (6,1.5) -- (1.5,0);
\draw[gray, thick] (6,-1.5) -- (1.5,0);

\draw[gray, thick] (6,1.5) -- (10.5,0);
\draw[gray, thick] (6,-1.5) -- (10.5,0);


\draw[gray, thick] (1.5, 0) .. controls  (0,5) and (12,5) .. (10.5,0);


\filldraw[red] (4,0) circle (2pt) node[anchor=west]{};
\filldraw[red] (8,0) circle (2pt) node[anchor=west]{};
\filldraw[red] (6,2.5) circle (2pt) node[anchor=west]{};
\filldraw[red] (6,-2.5) circle (2pt) node[anchor=west]{};


\node  at (5.3,0) {\small $v-u$};
\node[text=blue]  at (6.3,0) {\small$3$};

\node  at (6,4) {\small $v-u$};
\node[text=blue]  at (6,3.5) {\small$6$};

\node  at (3,1) {\small $-v-c_\hbar$};
\node[text=blue]  at (4.5,1.3) {\small$1$};
\node  at (3,-1) {\small $u-c_\hbar$};
\node[text=blue]  at (4.5,-1.3) {\small$2$};

\node  at (9,1) {\small $u-c_\hbar$};

\node[text=blue]  at (7.5,1.3) {\small$4$};
\node  at (9,-1) {\small $-v-c_\hbar$};

\node[text=blue]  at (7.5,-1.3) {\small$5$};

\end{tikzpicture}
  \caption{The  framed seed $\underline{\mathbf{i}}_2=\sigma(\underline{\mathbf{i}}_1)$ for the canoe graph. } 
  \label{fig:reframed-canoe} 
 \end{minipage}%
 \begin{minipage}{0.5\linewidth} 
  \centering 
\begin{tikzpicture}[scale=.72]
\draw[gray, thick] (3,0) circle (1.5cm);
\draw[gray, thick] (4.5,0) -- (7.5,0);
\draw[gray, thick] (9,0) circle (1.5cm);

\draw[gray, thick] (1.5, 0) .. controls  (0,5) and (12,5) .. (10.5,0);


\filldraw[red] (3,0) circle (2pt) node[anchor=west]{};
\filldraw[red] (9,0) circle (2pt) node[anchor=west]{};
\filldraw[red] (6,2) circle (2pt) node[anchor=west]{};
\filldraw[red] (6,-2) circle (2pt) node[anchor=west]{};


\node at (6,.3) {\small$u-c_\hbar$};
\node[text=blue]  at (6,-.3) {\small$2$};

\node (b) at (6,4) {\small $-u+c_\hbar$};
\node[text=blue]  at (6,3.5) {\small$4$};

\node  at (3,1.8) {$-v-c_\hbar$};
\node[text=blue] at (3,1.2) {\small$1$};
\node[text=blue]  at (3,-1.2) {\small$6$};
\node (e) at (3,-1.8) {$v-c_\hbar$};

\node (d) at (9,1.8) {$v-c_\hbar$};
\node[text=blue] at (9,1.2) {\small$3$};
\node[text=blue]  at (9,-1.2) {\small$5$};
\node (e) at (9,-1.8) {$-v-c_\hbar$};
\end{tikzpicture}

  \caption{The  framed seed $\mu_4^+(\underline{\mathbf{i}}_2)\simeq \underline{\mathbf{i}}_0$.} 
  \label{fig:final-g1-necklace} 
 \end{minipage} 
\end{figure}

Finally, performing a positive mutation $\mu^+_4$ at edge 4 in $\underline{\mathbf{i}}_2$ results in the framed seed shown in Figure~\ref{fig:final-g1-necklace}, which represents the same framed seed as the initial one $\underline{\mathbf{i}}_0$.  We therefore have a loop in the framed seeds groupoid
\begin{equation}
\label{eq:loop-g1-second}
\begin{tikzcd}[row sep=huge]
 &
\underline{\mathbf{i}}_0 \arrow[dl,swap,"\mu_3^+"]  &
\\
\underline{\mathbf{i}}_1 \arrow[rr,"\sigma"]& & \underline{\mathbf{i}}_2 \arrow[ul,swap,"\mu_4^+"]
\end{tikzcd}
\end{equation}

We once again take the wavefunction for the standard necklace framed seed $\underline{\mathbf{i}}_0$ to be $\psi_{\underline{\mathbf{i}}_0}=1$, but now regarded as an entire function on $\mathbb{C}$ rather than as a formal power series. We now explain how the mutations and framing shifts in~\eqref{eq:loop-g1-second} give rise to operators on spaces of meromorphic functions with appropriate analytic properties, and verify that the composite of these operators indeed preserves $\psi_{\underline{\mathbf{i}}_0}$ up to a phase.

By analogy with the Fock-Goncharov construction in the unitary case, we take the positive mutation $\mu_3^+$ at edge 3 carrying Heisenberg element by $u-c_\hbar$ to correspond to the operator of multiplication by the meromorphic function $\varphi(z-c_\hbar)$, which has simple poles at  $\{in\hbar+im\hbar^{-1}\}_{n,m\in\mathbb{Z}_{\geq1}}$. We thus obtain
\begin{align*}
\psi_{\underline{\mathbf{i}}_1} &= \varphi(z-c_\hbar)\cdot \psi_{\underline{\mathbf{i}}_0}\\
&=\varphi(z-c_\hbar)
\end{align*}
which now satisfies the dual pair of face relations
\begin{align*}
e^{-2\pi\hbar z}\psi_{\underline{\mathbf{i}}_1}(z)+(1-e^{-2\pi\hbar z})\psi_{\underline{\mathbf{i}}_1}(z+i\hbar) &= 0\\
e^{-2\pi\hbar^{-1} z}\psi_{\underline{\mathbf{i}}_1}(z)+(1-e^{-2\pi\hbar^{-1} z})\psi_{\underline{\mathbf{i}}_1}(z+i\hbar^{-1}) &=0.
\end{align*}
Let us regard the function $\psi_{\underline{\mathbf{i}}_1}$ as an element of the space $\mathcal{V}_{\underline{\mathbf{i}}_1}$ consisting of functions $f(z)$ analytic outside of the cone $\{in\hbar+im\hbar^{-1}\}_{n,m\in\mathbb{R}_{\geq1}}$, and having prescribed asymptotic behavior
\begin{align}
\label{eq:asymptotic-conditions}
f(z)\big |_{z\rightarrow\infty}\sim\begin{cases}
A_- \quad & |\arg(z)|>\frac{\pi}{2}+\arg(\hbar)\\
A_+e^{\pi i (z-c_\hbar)^2} \quad & |\arg(z)|<\frac{\pi}{2}-\arg(\hbar)
\end{cases}
\end{align}
for some $A_\pm\in\mathbb{C}$. 


Let us now consider the effect of performing the change of framing $\sigma$.  As in the case of mutation, we again define its action on our wavefunction by analytic continuation of the integral transform representing the action of $e^{\pi iv^2 }$ in the unitary case. Indeed, consider the integral
\begin{equation}
\label{eq:framing-kernel}
f\longmapsto \int e^{-\pi i (z-t-c_\hbar)^2}f(t+2c_\hbar)dt,
\end{equation}
where the contour of integration stays within the domain of analyticity of $f$ and escapes to infinity in the sectors $|\arg(t)|>\frac{\pi}{2}+\arg(\hbar)$ and $|\arg(t)|<\frac{\pi}{2}-\arg(\hbar)$. It follows from the asymptotics~\eqref{eq:asymptotic-conditions} that the integral converges absolutely for $|\arg(z)-\frac{\pi}{2}|<\pi-\arg(\hbar)$, so that $\hat{f}(z)$ defines an analytic function on the complement of the cone $\{-in\hbar-im\hbar^{-1}\}_{n,m\in\mathbb{R}_{\geq0}}$.

Applying~\eqref{eq:framing-kernel} to $\psi_{\mathrm{canoe}}$ and using the inversion and Fourier transformation properties~\eqref{inv} and ~\eqref{eq:qdl-fourier} of the noncompact quantum dilogarithm, we obtain
\begin{align*}
 \psi_{\underline{\mathbf{i}}_2}&=\int e^{-\pi i (z-t-c_\hbar)^2}\psi_{\underline{\mathbf{i}}_1}(t+2c_\hbar)dt\\
 &= e^{-\pi i z^2 +2\pi i c_\hbar z }\int e^{2\pi i zt}e^{-\pi i (t+c_\hbar)^2}\varphi(t+c_\hbar)dt\\
 &=\zeta_{inv}e^{-\pi i z^2 +2\pi i c_\hbar z }\int \frac{e^{2\pi i zt}}{\varphi(-t-c_\hbar)}dt\\
 &=\zeta_{inv}\zeta e^{-\pi i z^2 +2\pi i c_\hbar z }\varphi(-z+c_\hbar)\\
  &=e^{\pi i c_\hbar^2}\zeta_{inv}^2\zeta \cdot \varphi(z-c_\hbar)^{-1}.
\end{align*}
Finally, the positive mutation $\mu_4^+$ at edge 4 of $\underline{\mathbf{i}}_2$ which carries Heisenberg element $u-c_\hbar$ acts by the operator of multiplication by the meromorphic function $\varphi(z-c_\hbar)$, and hence under our proposed action for framing shifts and admissible mutations the loop~\eqref{eq:loop-g1-second} does indeed act trivially on the analytic wavefunction  $\psi_{\underline{\mathbf{i}}_0}$, up to a constant phase. 

Let us conclude our discussion of the analytic picture for $g=1$ case with an example of a genuinely non-algebraic wavefunction. Recall the framed seed from Example~\ref{eg:non-algebraic} obtained from that in Figure~\ref{fig:std-g1-necklace} by performing positive mutations at the edges labelled $3$ and $6$, for which we showed no algebraic wavefunction exists. On the other hand, following the prescription above, we obtain the corresponding analytic wavefunction associated to the this framed seed:
\begin{align*}
\Psi &= \varphi(z-c_\hbar)\varphi(-z+c_\hbar)\\
&=\zeta_{inv}e^{\pi i (z-c_\hbar)^2},
\end{align*}
where we again used the inversion formula~\eqref{inv} for the noncompact quantum dilogarithm.

\subsection{Analytic wavefunctions for \texorpdfstring{$g=2$}{g=2}}
In the genus 2 case, the combinatorics of ideal foams and framed seeds becomes richer. To illustrate this, we will describe a loop in the framed seeds groupoid that reflects a 3-2 Pachner move for deformed foams, and verify that this loop acts trivially on our proposed analytic wavefunction.

Again we begin with the standard necklace framed seed $\underline{\mathbf{i}}_{\mathrm{neck}}$ for $g=2$, for which $\psi_{\mathrm{neck}}=1$. Performing positive mutations at the (commuting) edges labelled $u_1-c_\hbar, ~u_2-c_\hbar$, we obtain the framed seed for the canoe graph shown $\underline{\mathbf{i}}_2$ in Figure~\ref{fig:2flips}, whose underlying deformed foam consists of two tetrahedra.  \begin{figure}[htbp] 
 \begin{minipage}{1\linewidth} 
  \centering 
\begin{tikzpicture}[scale=2.9]
\pgfmathsetmacro{\Ax}{0}
 \pgfmathsetmacro{\Ay}{3}
 \pgfmathsetmacro{\Bx}{2.3}
 \pgfmathsetmacro{\By}{3.3}
 \pgfmathsetmacro{\Cx}{1}
 \pgfmathsetmacro{\Cy}{4}
 \pgfmathsetmacro{\Dx}{0}
 \pgfmathsetmacro{\Dy}{0}
 \pgfmathsetmacro{\Ex}{2.3}
 \pgfmathsetmacro{\Ey}{0.3} 
\pgfmathsetmacro{\Fx}{1}
 \pgfmathsetmacro{\Fy}{1}
 \coordinate (A) at (\Ax,\Ay);
 \coordinate (B) at (\Bx,\By);
 \coordinate (C) at (\Cx,\Cy);
 \coordinate (D) at (\Dx,\Dy);
 \coordinate (E) at (\Ex,\Ey);
  \coordinate (F) at (\Fx,\Fy);
 \coordinate (Amid) at (.7,1.5);
  \coordinate (Amidtext) at (.7-.3,1.5);
 \coordinate (Amid2) at (.5,1.7);
 \coordinate (Bmid) at (1.5,1.6);
 \coordinate (Bmid2) at (1.7,1.7);
  \coordinate (Bmidtext) at (1.6,1.5);
 \coordinate (Cmid) at (1.32,2.3);
 \coordinate (Cmid2) at (1.2,2.1);
  \coordinate (Cmidtext) at (1.15,2.25);


 \coordinate (ACmid) at (\Ax/2+\Cx/2-.5,\Ay/2+\Cy/2);
 \coordinate (BCmid) at (\Bx/2+\Cx/2+.3,\By/2+\Cy/2+.1);
 \coordinate (ABmid) at (\Ax/2+\Bx/2+.4,\Ay/2+\By/2-.15);
 \coordinate (DFmid) at (\Dx/2+\Fx/2+.3,\Dy/2+\Fy/2-.1);
 \coordinate (EFmid) at (\Ex/2+\Fx/2,\Ey/2+\Fy/2+.25);
 \coordinate (DEmid) at (\Dx/2+\Ex/2+.2,\Dy/2+\Ey/2-.25);
 \coordinate (ADmid) at (\Ax/2+\Dx/2-.6,\Ay/2+\Dy/2);
 \coordinate (CFmid) at (\Fx/2+\Cx/2+.45,\Fy/2+\Cy/2-.4);
 \coordinate (BEmid) at (\Bx/2+\Ex/2+.5,\By/2+\Ey/2);
 





\draw[blue,very thick] (A)--(B)--(C)--(A)--(D)--(E)--(B); 
\draw[blue,very thick,dashed] (D)--(F)--(E);
\draw[blue,very thick,dashed] (F)--(C);


\node at (ACmid) {\small $u_1+v_1-2c_\hbar$}; 
\node at (BCmid) {\small $-v_1-c_\hbar$};
\node at (ABmid) {\small$-u_1+c_\hbar$}; 
\node at (DFmid) {\small $-v_2-c_\hbar$};
\node at (EFmid) {\small$u_2+v_2-2c_\hbar$};
\node at (DEmid) {\small$-u_2+c_\hbar$};
\node at (ADmid) {\small$u_2+v_2-v_1-2c_\hbar$}; 
\node at (CFmid) {\small$-u_1-u_2+3c_\hbar$}; 
\node at (BEmid) {\small$u_1+v_1-v_2-2c_\hbar$}; 

\end{tikzpicture}
  \caption{The framed seed $\underline{\mathbf{i}}_2$} 
  \label{fig:2flips} 
 \end{minipage}%
\end{figure}
The corresponding wavefunction is 
$$
\psi_{\underline{\mathbf{i}}_2} = \varphi(z_1-c_\hbar)\varphi(z_2-c_\hbar).
$$
On the other hand, consider the framed seed $\underline{\mathbf{i}}_3$ obtained from $\underline{\mathbf{i}}_{\mathrm{neck}}$ by instead performing \emph{negative} mutations at the edges labelled $u_1-c_\hbar, u_2-c_\hbar$, followed by a \emph{positive} mutation at the edge labelled $-u_1-u_2+3c_\hbar$. This framed seed is illustrated in Figure~\ref{fig:3flips}. 
\begin{figure}[htbp] 
 \begin{minipage}{1\linewidth} 
  \centering 
\begin{tikzpicture}[scale=2.9]
\pgfmathsetmacro{\Ax}{0}
 \pgfmathsetmacro{\Ay}{3}
 \pgfmathsetmacro{\Bx}{2.3}
 \pgfmathsetmacro{\By}{3.3}
 \pgfmathsetmacro{\Cx}{1}
 \pgfmathsetmacro{\Cy}{4}
 \pgfmathsetmacro{\Dx}{0}
 \pgfmathsetmacro{\Dy}{0}
 \pgfmathsetmacro{\Ex}{2.3}
 \pgfmathsetmacro{\Ey}{0.3} 
\pgfmathsetmacro{\Fx}{1}
 \pgfmathsetmacro{\Fy}{1}
 \coordinate (A) at (\Ax,\Ay);
 \coordinate (B) at (\Bx,\By);
 \coordinate (C) at (\Cx,\Cy);
 \coordinate (D) at (\Dx,\Dy);
 \coordinate (E) at (\Ex,\Ey);
  \coordinate (F) at (\Fx,\Fy);
 \coordinate (Amid) at (.7,1.5);
  \coordinate (Amidtext) at (.7-.3,1.5);
 \coordinate (Amid2) at (.5,1.7);
 \coordinate (Bmid) at (1.5,1.6);
 \coordinate (Bmid2) at (1.7,1.7);
  \coordinate (Bmidtext) at (1.6,1.5);
 \coordinate (Cmid) at (1.32,2.3);
 \coordinate (Cmid2) at (1.2,2.1);
  \coordinate (Cmidtext) at (1.15,2.25);


 \coordinate (ACmid) at (\Ax/2+\Cx/2-.5,\Ay/2+\Cy/2);
 \coordinate (BCmid) at (\Bx/2+\Cx/2+.5,\By/2+\Cy/2+.1);
 \coordinate (ABmid) at (\Ax/2+\Bx/2+.4,\Ay/2+\By/2-.15);
 \coordinate (DFmid) at (\Dx/2+\Fx/2+.3,\Dy/2+\Fy/2-.1);
 \coordinate (EFmid) at (\Ex/2+\Fx/2,\Ey/2+\Fy/2+.32);
 \coordinate (DEmid) at (\Dx/2+\Ex/2+.2,\Dy/2+\Ey/2-.25);
 \coordinate (ADmid) at (\Ax/2+\Dx/2-.4,\Ay/2+\Dy/2);
 \coordinate (CFmid) at (\Fx/2+\Cx/2-.3,\Fy/2+\Cy/2-.4);
 \coordinate (BEmid) at (\Bx/2+\Ex/2+.35,\By/2+\Ey/2);
 





\draw[blue,very thick] (A)--(B)--(C)--(A)--(D)--(E)--(B); 
\draw[blue,very thick,dashed] (D)--(F)--(E);
\draw[blue,very thick,dashed] (F)--(C);


\node at (ACmid) {\small $v_1-c_\hbar$};  %
\node at (BCmid) {\small $u_1+v_2-v_1-2c_\hbar$};%
\node at (ABmid) {\small$-u_1-v_2+c_\hbar$}; %
\node at (DFmid) {\small $-u_2-v_1+c_\hbar$};%
\node at (EFmid) {\small$u_2+v_1-v_2-2c_\hbar$};%
\node at (DEmid) {\small$v_2-c_\hbar$}; %
\node at (ADmid) {\small$u_1+u_2-3c_\hbar$}; 
\node at (CFmid) {\small$-u_1+c_\hbar$}; 
\node at (BEmid) {\small$-u_2+c_\hbar$}; 

\end{tikzpicture}
  \caption{The framed seed $\underline{\mathbf{i}}_3$} 
  \label{fig:3flips} 
 \end{minipage}%
\end{figure}
The corresponding deformed foam now consists of three deformed Harvey-Lawson tetrahedra, and the wavefunction is
$$
\psi_{\underline{\mathbf{i}}_3} = \frac{\varphi(-z_1-z_2+3c_\hbar)}{\varphi(-z_1+c_\hbar)\varphi(-z_2+c_\hbar)}.
$$
Introducing the following composite of framing shift and coordinate rescaling operators
\begin{align*}
\sigma ~:~ u_1 &\longmapsto u_1 + v_2- v_1+3c_\hbar\\
  u_2 &\longmapsto u_2 - v_2+v_1+3c_\hbar
\end{align*}
and the change of coordinates
$$
\tau ~:~ u_j\mapsto -u_j, \quad v_j\mapsto -v_{j},
$$
we observe that the framed seed $(\tau\circ\sigma)\cdot\underline{\mathbf{i}}_2$ coincides with 
$\underline{\mathbf{i}}_3$ up to a re-labelling of edges of the cubic graph. We will now confirm that the corresponding wavefunctions are indeed projectively equal. 

The action of the operator $\sigma = e^{-6\pi ic_\hbar(v_1+v_2)}e^{\pi i(v_2-v_1)^2}$ on $\psi_{\underline{\mathbf{i}}_2}$ can be understood with the help of the following Lemma:
\begin{lemma}
\label{lem:23}
We have
$$
e^{-6\pi ic_\hbar(v_1+v_2)}e^{\pi i(v_2-v_1)^2}\cdot\varphi(z_1-c_\hbar)\varphi(z_2-c_\hbar)\equiv \frac{\varphi(z_1+z_2+3c_\hbar)}{\varphi(z_2+c_\hbar)\varphi(z_1+c_\hbar)},
$$
where the symbol $\equiv$ denotes projective equality modulo phase constants.
\end{lemma}

\begin{proof}
Commuting the operator $e^{-6\pi ic_\hbar(v_1+v_2)}e^{\pi i(v_2-v_1)^2}$ past $\varphi(z_1-c_\hbar)\varphi(z_1-c_\hbar)$, we see that
\begin{align*}
\sigma\cdot\varphi(z_1-c_\hbar)\varphi(z_2-c_\hbar) &= \varphi(u_1 +v_2-v_1+2c_\hbar)\varphi(u_2-v_2+v_1+2c_\hbar)\cdot \sigma\cdot1\\
&=\varphi(u_1 +v_2-v_1+2c_\hbar)\varphi(u_2-v_2+v_1+2c_\hbar)\cdot1.
\end{align*}
The action of the latter operators are once again understood by means of the Fourier self-duality ~\eqref{eq:qdl-fourier}, so that we have, e.g. 
\begin{align*}
\varphi(u_1 +v_2-v_1+2c_\hbar)\cdot f(z_1,z_2) &= \zeta^{-1}\int \frac{e^{2\pi i t(u_1+v_2-v_1+c_\hbar)}}{\varphi(t-c_\hbar)}f(z_1,z_2)dt\\
&=\zeta^{-1}\int \frac{e^{\pi i t^2}e^{2\pi i t(u_1+c_\hbar)}}{\varphi(t-c_\hbar)}e^{2\pi i t(v_2-v_1)}\cdot f(z_1,z_2)dt\\
&=\zeta^{-1}\int \frac{e^{\pi i t^2}e^{2\pi i t(u_1+c_\hbar)}}{\varphi(t-c_\hbar)}\cdot f(z_1+t,z_2-t)dt\\
&=e^{-\pi ic_\hbar^2}(\zeta_{inv}\zeta)^{-1}\int e^{-2\pi i t(u_1+2c_\hbar)}{\varphi(t+c_\hbar)}\cdot f(z_1-t,z_2+t)dt
\end{align*}
Hence we see that up to multiplicative phase constants,
\begin{align*}
\sigma\cdot \psi_{\underline{\mathbf{i}}_2} &\equiv \int {e^{-2\pi i t(z_1+2c_\hbar)}e^{-2\pi i s(z_2+t+2c_\hbar)}}{\varphi(t+c_\hbar)\varphi(s+c_\hbar)} dsdt\\
&\equiv \int \frac{e^{-2\pi i t(z_1+2c_\hbar)}\varphi(t+c_\hbar)}{\varphi(t+z_2+c_\hbar)}dt,
\end{align*}
where the Fourier integral over $s$ is again performed using~\eqref{eq:qdl-fourier}. On the other hand, the resulting integral over $t$ may be computed by means of the `pentagon' integral evaluation~\eqref{beta-2}, with the result
\begin{align*}
\sigma\cdot \psi_{\underline{\mathbf{i}}_2} &\equiv \frac{\varphi(z_1+z_2+3c_\hbar)}{\varphi(z_2+c_\hbar)\varphi(z_1+c_\hbar)}.
\end{align*}
\end{proof}
Hence we conclude that 
$$
(\tau\circ\sigma)\cdot \psi_{\underline{\mathbf{i}}_2} \equiv \psi_{\underline{\mathbf{i}}_3}.
$$

\begin{example}
\label{sec:cube}
An interesting example, explored at the semiclassical level in Section 5.4 of~\cite{TZ}, is the genus 3 cubic graph obtained as the 1-skeleton of the cube shown in Figure~\ref{fig:cube}. 
\begin{figure}
\begin{tikzpicture}
\pgfmathsetmacro{\Ei}{1.2}
\pgfmathsetmacro{\Eo}{3}
\pgfmathsetmacro{\shv}{.25}
\pgfmathsetmacro{\shh}{.6}

\draw[thick ] (-\Ei,-\Ei)--(-\Ei,\Ei);
\draw[thick] (-\Ei,\Ei)--(\Ei,\Ei);
\draw[thick] (\Ei,-\Ei)--(\Ei,\Ei);
\draw[thick] (-\Ei,-\Ei)--(\Ei,-\Ei);

\draw[thick ] (-\Eo,-\Eo)--(-\Eo,\Eo);
\draw[thick] (-\Eo,\Eo)--(\Eo,\Eo);
\draw[thick] (\Eo,-\Eo)--(\Eo,\Eo);
\draw[thick] (-\Eo,-\Eo)--(\Eo,-\Eo);

\draw[thick] (\Ei,\Ei)--(\Eo,\Eo);
\draw[thick] (-\Ei,\Ei)--(-\Eo,\Eo);
\draw[thick] (-\Ei,-\Ei)--(-\Eo,-\Eo);
\draw[thick] (\Ei,-\Ei)--(\Eo,-\Eo);

\node[red] at (0,\Ei-\shv) {\small 3};

\node[red] at (\Ei-\shv,0) {\small 2};

\node[red] at (-\Ei+\shv,0) {\small 4};

\node[red] at (0,-\Ei+\shv) {\small 1};

\node[red] at (0,\Eo-\shv) {\small 7};

\node[red] at (\Eo-\shv,0) {\small 6};

\node[red] at (0,-\Eo+\shv) {\small 5};

\node[red] at (-\Eo+\shv,0) {\small 8};

\node[red] at (\Ei/2+\Eo/2+\shv,\Ei/2+\Eo/2-\shv) {\small $11$};

\node[red] at (-\Ei/2-\Eo/2-\shv,\Ei/2+\Eo/2-\shh/2)  {\small $12$};

\node[red] at (-\Ei/2-\Eo/2-\shv,-\Ei/2-\Eo/2+\shv/2)  {\small $9$};

\node[red] at (\Ei/2+\Eo/2+\shh/2,-\Ei/2-\Eo/2+\shv/2) {\small $10$};

\end{tikzpicture}
\caption{The cubic graph of genus $3$ given by the 1-skeleton of the cube.}
\label{fig:cube}
\end{figure}
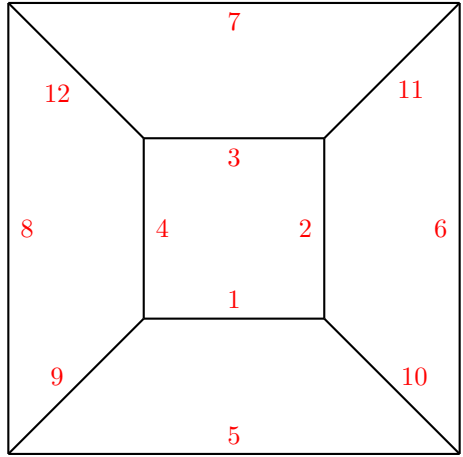
In our setting, the phase and framing considered in~\cite{TZ} correspond to the object $\underline{\mathbf{i}}_{\mathrm{cube}}$ of the framed seed groupoid in which the edges $e_i$ are labelled by Heisenberg algebra elements
\begin{align*}
\nonumber\underline{\mathbf{i}}_{\mathrm{cube}}:~&e_1\mapsto u_2-v_2-3c_\hbar \qquad e_2\mapsto u_2+c \qquad e_3\mapsto v_2-c_\hbar \qquad  e_4\mapsto c_\hbar-u_1-u_2\\
&e_5\mapsto c_\hbar-u_1-u_2 \qquad e_6\mapsto v_3-c_\hbar \qquad e_7\mapsto u_3+c_\hbar \qquad e_8\mapsto u_1-v_3-3c_\hbar\\
\nonumber &e_9\mapsto v_1-c_\hbar \qquad e_{10}\mapsto u_3+v_2-v_1+3c_\hbar \qquad e_{11}\mapsto v_1-v_2-u_2-u_3-3c_\hbar\\
 & e_{12}\mapsto u_2+v_3-v_1+c_\hbar
\end{align*} 
Let us now explain how to derive the analytic wavefunction associated to this framed seed, and verify its semiclassical limit reproduces the prediction for the holomorphic disk invariants given in~\cite{TZ}.

We begin with the standard $g=3$ necklace framed seed $\underline{\mathbf{i}}_{\mathrm{neck}}$, and perform positive mutations at the strands labelled $u_1-c_\hbar,u_2-c_\hbar,u_3-c_\hbar$, as well as a negative mutation at the strand labelled $5c_\hbar-u_1-u_2-u_3$. The resulting cubic graph is isomorphic to the 1-skeleton of the cube, and the corresponding framed seed $\underline{\mathbf{i}}_{\mathrm{cube}}'$ is given by
\begin{align*}
\nonumber\underline{\mathbf{i}}_{\mathrm{cube}}':~&e_1\mapsto c_\hbar-u_1 \qquad e_2\mapsto u_2+v_2-v_1-2c_\hbar \qquad e_3\mapsto c_\hbar-u_2 \qquad  e_4\mapsto u_1+v_1-v_2-2c_\hbar\\
&e_5\mapsto u_1+u_2+u_3-5c_\hbar \qquad e_6\mapsto c_\hbar-v_3 \qquad e_7\mapsto c_\hbar-u_3 \qquad e_8\mapsto v_3-u_1-u_2+3c_\hbar\\
\nonumber &e_9\mapsto -v_1-c_\hbar \qquad e_{10}\mapsto v_1-u_2-u_3+3c_\hbar \qquad e_{11}\mapsto u_3+v_3-v_2-2c_\hbar\\
 & e_{12}\mapsto u_2+v_2-v_3-2c_\hbar
\end{align*} 
The wavefunction associated to $\underline{\mathbf{i}}_{\mathrm{cube}}'$ is given by
$$
\psi_{\mathrm{cube}}' \equiv \frac{\varphi(z_1-c_\hbar)\varphi(z_2-c_\hbar)\varphi(z_3-c_\hbar)}{\varphi(z_1+z_2+z_3-5c_\hbar)}.
$$

 To pass from the framed seed $\underline{\mathbf{i}}_{\mathrm{cube}}'$ to the desired one $\underline{\mathbf{i}}_{\mathrm{cube}}$, we first apply the coordinate rescaling and framing shift operators $e^{-6\pi ic_\hbar(v_1+v_2)}e^{\pi i (v_1-v_2)^2}$. By Lemma~\ref{lem:23}, their effect on the wavefunction is given by
 $$
 e^{-6\pi ic_\hbar(v_1+v_2)}e^{\pi i (v_1-v_2)^2}\cdot \psi_{\mathrm{cube}}' \equiv \frac{\varphi(z_1+z_2+3c_\hbar)\varphi(z_3-c_\hbar)}{\varphi(z_1+z_2+z_3+c_\hbar)\varphi(z_1+c_\hbar)\varphi(z_2+c_\hbar)}.
 $$
In the resulting framed seed, the edge numbered $2$ in Figure~\ref{fig:cube} now carries the label $u_2+c_\hbar$.  We now perform two \emph{consecutive} positive mutations at this edge. Note that in the setting of algebraic representations and wavefunctions, we can never perform two such consecutive mutations of the same sign at an edge, as it is impossible for both mutations to be admissible in the sense of Section~\ref{sec:groupoid-rep}. In the analytic setting, on the other hand, the action of these two mutations on the wavefunction yields
\begin{align*}
\psi_{\mathrm{cube}}'' &\equiv \varphi(-u_2-c_\hbar)\varphi(u_2+c_\hbar)\cdot\psi_{\mathrm{cube}}'\\
&\equiv
 \frac{\varphi(z_1+z_2+3c_\hbar)\varphi(z_3-c_\hbar)\varphi(-z_2-c_\hbar)}{\varphi(z_1+z_2+z_3+c_\hbar)\varphi(z_1+c_\hbar)}. 
\end{align*}
After applying the symplectic lift $\tau$ of the change of basis $(u_1,u_2,u_3)\mapsto (-u_1-u_2, ~u_2,~-u_3)$, the resulting framed seed  can be identified with $\underline{\mathbf{i}}_{\mathrm{cube}}$, and the corresponding wavefunction is given by
\begin{align}
\nonumber \psi_{\mathrm{cube}} &\equiv \tau\cdot\psi_{\mathrm{cube}}''\\
&\equiv \frac{\varphi(-z_1+3c_\hbar)\varphi(-z_2-c_\hbar)\varphi(-z_3-c_\hbar)}{\varphi(-z_1-z_3+c_\hbar)\varphi(-z_1-z_2+c_\hbar)}.
\end{align}
In this setting, the semiclassical limit is realized by sending $\hbar\rightarrow 0$ while rescaling $z_i\mapsto (2\pi\hbar)^{-1}z_i$. Using the asymptotics for $\varphi$ given in the Appendix and setting $Z_i = e^{-z_i}$ as in~\cite{TZ}, we find
$$
2\pi i \hbar^2\log(\psi_{\mathrm{cube}})\sim \mathrm{Li}_2(Z_1)+\mathrm{Li}_2(Z_2)+\mathrm{Li}_2(Z_3) - \mathrm{Li}_2(Z_1Z_2)-\mathrm{Li}_2(Z_1Z_3),
$$
which coincides with the expression for the superpotential $W$ in Section 5.4 of~\cite{TZ}.
\end{example}

%% file: framing_and_quiver.tex
\section{Framing Duality}

In this section we observe a curious relationship we call \emph{framing duality}.  For the class of Legendrian surfaces
which generalize the Clifford torus to arbitrary genus, the wavefunctions for different framings correspond
to Donaldson-Thomas/Hall-Algebra generating functions of different quivers, as computed in \cite{KS}.  The framings are defined by
$g\times g$ symmetric integer matrices, $A$.  When all entries are non-negative,
the matrix determines a quiver $Q_A$ with
adjacency matrix $A$.  The quiver invariants are proven to be integers in \cite{E}, and are thus
equal to the Ooguri-Vafa integers for the corresponding brane. 
They are also conjectured in \cite{HRV} to count the dimensions of isotypic components of the middle cohomology of
twisted character varieties.

\subsection{Wavefunction for Clifford Surfaces}
\label{sec:framingwavefunction}

We define the \emph{Chekanov surface} of genus $g$ to be the 
Legendrian defined by the necklace graph $\Gamma_{\rm neck}^g$,
and the \emph{Clifford surface} to be
the defined by the canoe graph $\Gamma_{\rm canoe}^g$ ---
see Figure \ref{fig:necklacecanoe}. 

We now define the standard necklace framed seed $\underline{\bf i}^{\rm neck}_{g,0},$
generalizing the $g=1$ case
of Example \ref{eg:intertwining} (see also Figure \ref{fig:g1-necklace-alg}).
To fix notation and the cyclic structure, we embed $\Gamma_{\rm neck}^g$
in the plane $\bR^2$ with its standard orientation.  Let $I = \{0,1,...,2g+1\}$
and define the vertex set be $I\times \{0\}.$
Define the strand edges $s_1,...,s_g$ by $s_k = [2k-1,2k]\times \{0\}$ and set
$s_{g+1}$ to be a big loop in the upper half plane connecting $(0,0)$ and $(2g+1,0).$\footnote{As we are working on the Riemann sphere $S^2,$ there is no difference
between placing $s_{g+1}$ in the upper or lower half-plane, so the necklaces of Figures
\ref{fig:necklacecanoe} and \ref{fig:g1-necklace-alg} are in fact consistent.}
The two edges of the $k$th bead, $k = 1,...,g+1,$ are taken to lie on a circle of radius $1/2$ centered at $(2k-3/2,0),$ with the upper hemisphere called edge $a_i$ and lower hemisphere $b_i$, in the upper and lower half-planes, respectively.
We parametrize the edge variables $X_{e_i}$ with the quantum torus $\cD_{2g}$
as follows:
$$X_{s_k} = \begin{cases}
-q^{-1}U_k & 1 \leq k \leq g\\
-q^{2g-1}U_1^{-1}\cdots U_g^{-1} & k = g+1
\end{cases}
\qquad
X_{a_k} = \begin{cases}
-q^{-1}V_{1}^{-1} & k = 1\\
-q^{-1}V_{k-1}V_k^{-1} & 2 \leq k \leq g\\
-q^{-1}V_g & k = g+1
\end{cases}
\qquad 
X_{b_k} = \frac{q^{-2}}{X_{a_k}}
$$
It is straightforward to check that this assignment satisfies
Equations 
\eqref{quantum.face.0}
and 
\eqref{quantum.global}.
The quantized chromatic Lagrangian is the ideal defined by the face relations
of Equation
\eqref{quantum.face.2h}.
These impose $V_i = 1$, and nothing further, giving rise to the wavefunction
$$\Psi_{\underline{{\bf i}}^{\rm neck}_{g,0}}\equiv 1.$$

The canoe $\Gamma_{\rm canoe}^g$
is obtained by performing $g$ positive mutations at strands $1,...,g$,
similar to Examples \ref{eg:intertwining} and \ref{eg:canoe}.
These mutations are all admissible and mutually commuting.
The edge variable on the $k$th strand is $-q^{-1}U_k,$ so the
mutation is effected by conjugation by $\Phi(-q^{-1}U_k)^{-1} = (U_k;q^2)_\infty,$
since $\Phi(x) = \prod_{n\geq 0}(1+q^{2n+1}X)^{-1}.$ 
What results, then is the wavefunction
$$\Psi_{\underline{{\bf i}}^{\rm canoe}_{g,0}} = 
\prod_{k=1}^g (X_k,q^2)_\infty = 
\sum_{{\bf v}\in \bZ_{\geq 0}^g} \frac{1}{(q^2)_{\bf v}} X^{\bf v},
\qquad (q^2)_{\bf v}:=\prod_{i=1}^g \prod_{k=1}^{v_i}(1-q^{2k}),
$$
where on the right we have used the power-series expression of the
infinite Pochhammer symbol from Lemma \ref{lem:canoepowerseries}.
Now let $A$ be an $n\times n$, symmetric matrix with non-negative entries.
According to Equation \ref{eq:frame-change-action},
a frame-changing transformation by $A$ takes us to the seed
${\underline{{\bf i}}^{\rm canoe}_{g,A}}$
and 
\begin{equation}
\label{eq:wavefunctioncanoeA}\Psi_{\underline{{\bf i}}^{\rm canoe}_{g,A}} =
\sum_{{\bf v} \in \bZ_{\geq 0}^g}
\frac{q^{{\bf v}^t A {\bf v}}}{(q^2)_{\bf v}} X^{\bf v}.
\end{equation}
We define the Ooguri-Vafa invariants $n^{(A)}_{{\bf v},k}$ by passing to an infinite product expansion and setting
\[
\Psi_{\underline{{\bf i}}^{\rm canoe}_{g,A}}
=\prod_{{\bf v}\in \mathbb{Z}_{\geq 0}^n\setminus\{0\}}\prod_{k\in \mathbb{Z}}((-q)^{k}X^{\bf v};q^2)_\infty^{n^{(A)}_{{\bf v},k}}.
\]

Framing duality is the observation that, using the computation of
Kontsevich-Soibelman \cite{KS}, we can identify $\Psi_{\underline{{\bf i}}^{\rm canoe}_{g,A}}$
of Equation \eqref{eq:wavefunctioncanoeA}
as the DT series
of the quiver with adjacency matrix $A$.  This statement
will be made precise after defining
these terms in the next section.

\subsection{DT series for symmetric quivers.} 
\label{sec:dtsersforsymmquivs}

Let $A=(a_{ij})$ be an $n\times n$ symmetric matrix with nonnegative integral entries and let $Q_A$ be its corresponding symmetric quiver.
The generating function for the COHA $\mathcal{H}$ of $Q_A,$ also called the DT series, is  
\begin{align}
\label{DT.A}
\mathrm{DT}_{A} (t^{1/2}, X)&= \sum_{{\bf v}\in \mathbb{Z}_{\geq 0}^n, k\in \mathbb{Z}} (-1)^k \dim (\mathcal{H}_{{\bf v}, k}) t^{k/2} X^{\bf v} \\
\nonumber
\end{align}
where $\chi_A({\bf v},{\bf w}) := {\bf v}^t(I-A){\bf w}$ and
\[
\mathcal{H}_{{\bf v},k} = H^{k - {\chi_A({\bf v}, {\bf v})}}(BG_{{\bf v}}),\qquad G_{\bf v} = \prod_i GL_{v_i}(\bC).
\]
These were computed in \cite[\S 5.6]{KS} for symmetric quivers $Q_A,$ giving the result
\begin{equation}
    \label{eq:dtinv}
\mathrm{DT}_{A} (t^{1/2}, X) =: \sum _{{\bf v}\in \mathbb{Z}_{\geq 0}^n} DT^A_{\bf v} X^{\bf v} =   \sum _{{\bf v}\in \mathbb{Z}_{\geq 0}^n} \frac{(-t^{\frac{1}{2}})^{\chi_A({\bf v}, {\bf v})}}{(t)_{\bf v}} X^{\bf v}  \hskip 1cm \in~~~~~1+\mathfrak{m}
\end{equation}
where above we have defined the coefficient functions $DT^A_{\bf v}(t^{\frac{1}{2}})$, and once again
$(t)_{\bf v} := \prod_{i=1}^n(1-t)(1-t^2)\cdots (1-t^{v_i}).$

We now prove a lemma to be used in the next section.  First define $\sigma({\bf v}) = \sum_{i=1}^n v_i.$
Then we have:
\begin{lemma}
    \label{lem:coefficient-change}
$$DT_A(t^{\frac{1}{2}},-t^{-\frac{1}{2}}X) = DT_{I-A}(t^{\frac{1}{2}},X)$$
and
$$DT_{\bf v}^{I-A}(t^{\frac{1}{2}}) = (-t^{\frac{1}{2}})^{\sigma({\bf v})}DT^A_{\bf v}(t^{-\frac{1}{2}}).$$
\end{lemma}
\begin{proof}
To see this, note
$$(y^{-1})_{\bf v} = (y)_{\bf v}(-y^{\frac{1}{2}})^{-{\bf v}^t{\bf v}}(-y^{-\frac{1}{2}})^{\sigma({\bf v})}.$$
As a result, we have
$$\frac{(-t^{\frac{1}{2}})^{{\bf v}^tA{\bf v}}}{(t)_{\bf v}}X^{\bf v} =
\frac{(-t^{-\frac{1}{2}})^{{\bf v}^t(I-A){\bf v}}}{(-t^{-1})_{\bf v}}(-t^{-\frac{1}{2}}X)^{\bf v}$$
The second equation follows, and then the first.
\end{proof}
 
 \subsection{Framing Duality}
 \label{sec:framing-duality}
We now come to the main point of this section:  to compare
wavefunctions for canoe graphs with $\mathrm{DT}$ series of symmetric quivers.
 
Recall Equation \ref{eq:wavefunctioncanoeA} for the genus-$g$ canoe graph in framing $A$
from Section \ref{sec:framingwavefunction} above.
Comparing the form of its wavefunction given in Equation~\eqref{eq:canoe-general-framing} with
that of Equation~\eqref{eq:dtinv}, we have the following.
\begin{proposition}[Framing Duality]
For any integral $g\times g$ symmetric matrix $A$ with non-negative entries, the wavefunction $\Psi_{\mathbf{i}_{g,A}^\mathrm{canoe}}$ associated to the framed seed $\mathbf{i}_{g,A}^\mathrm{canoe}$ of the genus-$g$ canoe graph coincides with the $\mathrm{DT}$ series for the symmetric quiver with adjacency matrix $A$ under the identification $q = -t^{\frac{1}{2}}$:
 $$\Psi_{\mathbf{i}_{g,A}^\mathrm{canoe}} = \mathrm{DT}_{A} (-q, X).$$
 \end{proposition}
\begin{proof}
Comparing Equations \eqref{eq:dtinv} and \eqref{eq:wavefunctioncanoeA}, the proposition follows from
Lemma \ref{lem:coefficient-change}. 
\end{proof}

\subsection{Integer Invariants}

We define the quiver invariants $N^{(A)}_{{\bf v},k}$ by setting
\begin{equation}
\label{eq:dt-definition}
\mathrm{DT}_{A} (t^{\frac{1}{2}}, X) =\prod_{{\bf v}\in \mathbb{Z}_{\geq 0}^n\setminus\{0\}}\prod_{k\in \mathbb{Z}}(t^{\frac{k}{2}}X^{\bf v};t)_\infty^{N^{(A)}_{{\bf v},k}}
\end{equation}
We can rewrite this in another form using the plethystic exponential $\mathrm{Exp}$ and its inverse $\mathrm{Log}.$  Recall that for a power series vanishing at the origin, $f\in x\bC[[x]],$ we have
$\mathrm{Exp}(f) = \exp\left(\sum_{n\geq 1}\frac{f(x^n)}{n}\right).$  Then $\mathrm{Exp}(f+g) = \mathrm{Exp}(f)\mathrm{Exp}(g)$
and note $\mathrm{Exp}(x) = \frac{1}{1-x}.$  Then it is straightforward
to show
\[
(t-1){\rm Log}\,\mathrm{DT}_{A} (t^{1/2}, X) = \sum_{{\bf v}\in \mathbb{Z}_{\geq 0}^n\setminus\{0\}}
\sum_{k\in \mathbb{Z}}
N^{(A)}_{{\bf v}, k} t^{\frac{k}{2}} X^{\bf v}
\]
By \cite[Corollary 4.1]{E}, the values
$(-1)^{k-1}N^{(A)}_{{\bf v}, k}$ are non-negative integers and are nonzero only for finitely many $k\in \mathbb{Z}$.
In terms of integer invariants, framing duality says $N^{(A)}_{{\bf v},k} = n^{(A)}_{{\bf v},k}.$
{Note, however, that the $n^{(A)}_{{\bf v},k}$ are well-defined for non-positive $A$.}


\begin{remark} The paper \cite{HLRV} gives a cohomological interpretation of DT-invariants of quivers. Let $\Gamma$ be a quiver with $r$ vertices and let ${\bf v}=(v_1,..., v_r)\in \Z^{r}_{\geq 0}$ be a dimension vector. Associate to $(\Gamma, {\bf v})$ a new quiver $\tilde{\Gamma}$ by attaching a leg of length $v_i-1$ at the vertex $i$. We extend the dimension vector ${\bf v}$ to $\tilde{\bf v}$ by placing decreasing dimensions $v_i-1, v_i-2,..., 1$ at the extra leg.
Let $W_{\bf v}$ be the Weyl group of type $A_{v_1-1}\times \ldots \times A_{v_r-1}$ that is generated by the reflections at the extra vertices. 
Let $\mathcal{Q}_{\tilde{\bf v}}$ be the smooth generic complex quiver variety associated to $(\tilde{\Gamma}, \tilde{\bf v})$. The Weyl group $W_{\bf v}$ acts on $H^\ast_c(\mathcal{Q}_{\tilde{\bf v}}, \CC)$ and hence gives a natural decomposition of the latter into isotypical components. According to \cite[Cor 1.5]{HLRV}, after a slight renormalization, we have 
\[
{\rm DT}_{\Gamma, {\bf v}}(t^{\frac{1}{2}})= \sum_i \dim \big(H_c^{2i}(\mathcal{Q}_{\tilde{\bf v}}, \CC)^{W_{\bf v}}\big) t^{i-d_{\tilde{\bf v}}}.
\]
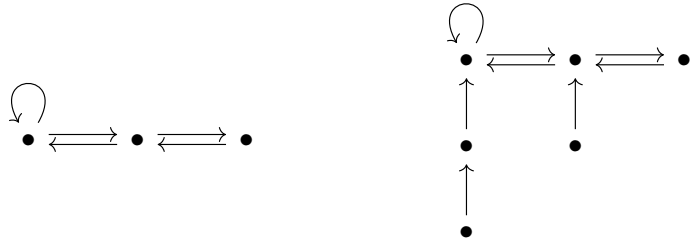
\begin{figure}[H]
\begin{tikzcd}
\bullet \arrow[r,shift left=.4ex] \arrow[out=60,in=120,loop] & \bullet \arrow[r, shift left=.4ex] \arrow[l, shift left=.4ex] & \bullet \arrow[l, shift left=.4ex]
\end{tikzcd}
\hskip 2cm 
\begin{tikzcd}
\bullet \arrow[r,shift left=.4ex] \arrow[out=60,in=120,loop] & \bullet \arrow[r, shift left=.4ex] \arrow[l, shift left=.4ex] & \bullet \arrow[l, shift left=.4ex]\\
\bullet \arrow[u]&\bullet\arrow[u]&\\
\bullet \arrow[u]& &
\end{tikzcd}
\caption{A new quiver associated to ${\bf v}=(3,2,1)$}
\end{figure}
\end{remark}

\subsection{Classical Limit}
\label{sec:classical-limit}

We define the quantity $W^{(A)}$ from the classical limit of the logarithm of the $\mathrm{DT}$ series.
We then relate it to the superpotential \`a la Aganagic-Vafa, through framing duality.

Comparing \eqref{eq:dt-definition} with \eqref{log}, we have 
\[
W^{(A)}(X):=\lim_{t^{\frac{1}{2}}\rightarrow 1}(1-t^{-1})\cdot \log {\rm DT}_{A} = \sum_{{\bf v}\in  \mathbb{Z}_{\geq 0}^n\backslash \{0\}} {\rm DT}_{\bf v}(1){\rm Li}_2(X^{\bf v})
\]
Then
\[
d W^{(A)}(X) = \sum_{i=1}^n \log Y_i \,  d\log X_i
\]
where
\[
Y_i:= \lim_{t^{1/2}\rightarrow 1} \frac{{\mathrm{DT}_{A}}(t^{1/2}; X_1,..., t X_i, ..., X_n)}{{\mathrm{DT}_{A}}(t^{1/2}; X_1,...,
X_i, ..., X_n)}.
\]

By \cite[Theorem 5.3]{KS},  $Y_1, ..., Y_n$ are solutions to the system of equations 
\[
X_i (-Y_i)^{1-a_{ii}}\big(\prod_{j\neq i} Y_j^{-a_{ij}}\big) + Y_i =1, \hskip 2cm i=1,..., n.
\]
Recalling that framing duality equates $q$ with $-t^{\frac{1}{2}}$, we should compare this
equation with the $q\to -1$ limit of the wavefunction of the genus-$g$ canoe. 
Upon setting $X_i = U_i$ and $Y_i = V_i$, we recognize this as 
the $q\to -1$ limit of the genus-$g$ canoe wavefunction of Equation \eqref{eq:genus-g-necklace}, after
changing frame using Equation \eqref{frame.change.matrix} with change-of-frame matrix $\Omega = {\bf 1} - A$.
The classical limit has a geometric interpretation in terms of moduli spaces as well, which we now discuss.

The canoe graph $\Gamma^{\rm canoe}_g$ is an iterated $g$-fold blow-up of the $\Theta$-graph with two nodes and three edges:
$\Gamma^{\rm canoe}_{g}$ is obtained from blowing up $\Gamma^{\rm canoe}_{g-1}$ at either of the two vertices at the ends of the canoe.
Recall from \cite[Section 5.2]{TZ} that if $\widehat{\Gamma}$ is the blow-up of $\Gamma$ at a vertex, then the moduli of objects of the corresponding sheaf
categories are related by $\cM_{\widehat \Gamma} = H \times \cM_\Gamma,$
where $H$ is the pair of pants $\bP^1\setminus\{0,1,\infty\}$.  So since $\cM_\Theta$ is a point, $\cM_{\Gamma^{\rm canoe}_g} \cong H^g,$ i.e.~$g$ copies of the tetrahedron moduli space.

After choosing a framing, we can define compatible coordinates on the torus in which $\cM_{\Gamma^{\mathrm canoe}_g}$ lives,
then lifting from $(\bC^\times)^{2g}$ to its half-universal cover $T^*(\bC^\times)^g$, we can write
$\cM_{\widehat{\Gamma}}$ as the graph of the differential of a superpotential.

In zero framing --- i.e., the one defined by mutation from the necklace as in Example \ref{eg:canoe} ---
we have
\begin{equation}
\label{eq:zeroframing}
W^{(0)}_{\Gamma^{\rm canoe}_g} = \sum_{i=1}^g \Li_2({X_i}),
\end{equation}
where $X_i = e^{u_i}$ --- and we will also need the conjugate logarithmic coordinates $v_i$ on the cotangent fibers.
A framing shift by a $g\times g$ symmetric integral matrix $A$,
as in Equation \eqref{frame.change.matrix} and Remark \ref{rmk:framing-BCH}, defines new coordinates
$v'_i = v_i, u'_i = u_i + A_{i,j}v_j.$  Then the lift of $\cM_{\Gamma^{\rm canoe}_g}$ is cut out from $\cP_{\Gamma^{\rm canoe}_g}$
in these coordinates as the graph of
the associated superpotential $W^{(A)}_{\Gamma^{\rm canoe}_g},$
which is the classical limit of the frame-shifted wavefunction.
We write
\begin{equation}
\label{eq:wg}
W^{(A)}_{\Gamma^{\rm canoe}_g} = \sum_{d \in (\bZ_{\geq 0})^g\setminus\{0\}} n^{(A)}_d \Li_2(X'^d),
\end{equation}
where $X'^d = \prod_{i=1}^g {X'_i}^{d_i}$ and $X'_i = e^{u'_i},$
and $n^{(A)}_d = \sum_{k} n^{(A)}_{d,k}.$
For $g=1$ and $A = (p),$ these integers appear (with a slightly different sign convention) in \cite[Section 6.1]{AKV}.

\subsection{Kac polynomial of a quiver}

We recall the Kac polynomial of a quiver.
Let $B$ be an $n\times n$ integal matrix with non-negative entries,
and let $Q_B$ be the quiver with $n$ nodes labeled $1,...,n,$ and $B_{i,j}$ arrows between node $i$ and node $j$.
The Kac polynomial of $Q_B$ is defined as follows.  Let $d\neq 0\in (\bZ_{\geq 0})^n$ be a dimension vector.
Then $$A^{(B)}_d(q) = \#\{\text{absolutely irreducible representations of $Q_B$ over $\bF_q$ modulo isomorphism}\}.$$

\begin{remark} The Kac polynomials are DT-invariants for quivers with potential. Let $Q$ be a quiver with arrows $a_{ij}$ at vertex $i$. Let $\widehat{Q}$ be the double quiver by adding arrows $a_{ij}^*$ of opposite direction and a new loop $c_i$ for each vertex. Take the quiver potential
\[
\widehat{W}=\sum_{i} c_i \sum_{j} [a_{ij}, a_{ij}^*].
\]
The Kac polynomials  for $Q$ are  the DT-invariants for the quiver-with-potential $(\widehat{Q}, \widehat{W})$.
\end{remark}

In \cite{RV}, a refinement of the Kac polynomial was introduced, in which the label is not simply a counting number $d$ but a 
partition $\lambda$.  Then $A_d(q) = \sum_{|\lambda|=d} A_\lambda(q).$
We will be interested in the special case $\lambda = 1^d = (1,1,...,1)$.

\begin{proposition}
\label{prop:kac-prop}
Let $h \geq 1$ and let $B = (2-2h)$ be the one-by-one matrix with single entry $2-2h,$
considered as a framing of the genus-one canoe.
Let $Q$ be the quiver with one node and $h$ arrows.
Then $$A_{1^d}(1) = n^{(B)}_d.$$
\end{proposition}
\begin{proof}
Consider first the case when $h=1,$ so $B=0.$
Then the quiver has no arrows and there is a unique irreducible representation for each $d,$ thus $A_d(q)=1$
and it is shown in \cite{RV} that this corresponds to the partition $(d).$  
When $d=1$ this equals $1^d,$ but not otherwise --- so we require $n^{(0)}_d = \delta_{d,1},$ which
agrees with Equation \eqref{eq:zeroframing} when $g=1.$
More generally, we refer to Equation (4.3.1) and Proposition 4.2.1 of \cite{RV}, where the notations $V, x, z,$ and $N$ are here
$W, X, e^v,$ and $g$, respectively.
In our notation, Equation (4.3.1) says
$dW^{(B)} = \sum_{i=1}^g v_i du_i.$
(Recall $g=1$ here.)
Writing $W^{(B)}_{\Gamma^{\rm canoe}_g}$ as in Equation \eqref{eq:wg}, this says
$e^{v_i} = \prod_{d}(1-X^d)^{d_i n_d^{(B)}}.$
Comparison with Equation (4.3.1) of \cite{RV} gives {$n^{(B)}_d = A_{1^d}(1),$ as claimed.}
We note that in \cite{RV} the function $W(X)$ is called Schl\"afli's differential by analogy with the volume of hyperbolic polyhedra,
which is part of a dual superpotential computation in \cite{DGGo}.

\end{proof}

\begin{example} Let us illustrate Proposition~\ref{prop:kac-prop} in the some examples where the relevant integer invariants have been recorded elsewhere in the literature.

Consider first the quiver with one node and two arrows.  The polynomials $A_{1^d}(q)$ are listed in \cite[Appendix II]{H} for $d = 1, 2, 3, 4$,
giving $A_{1^d}(1) = 1, 1, 3, 10,$ respectively.
On the other hand, we may compare these integers with the disk invariants for framing $p = 2 - 2\cdot 2 = -2$ obtained in the formulas of \cite[Section 6.1]{AKV} after their Equation 6.4, where they find $n^{(2)}_d = 1, 1, 3, 10$ for these same values of $d$.

For another class of examples, consider the quiver with one node and $g > 1$ arrows.
The polynomials $A_{1^2}(q)$ were computed in \cite[Section 3]{H}, giving $A_{1^2}(1) = g-1$.
This agrees with $n^{(2-2g)}_2$ as computed in \cite{AKV}.

For the reader's convenience, we record the following table of the integers $A_{1^d}(1)$ for the quiver with a single vertex and $h$ loops, which by Proposition~\ref{prop:kac-prop} coincide with $n_d^{(2-2h)}$:


\begin{figure}
\begin{tabular}{|c||ccccccc|}
\hline 
 $h$ & $d$=1 & 2 & 3 & 4 & 5 & 6 & 7   \\
   \hline
  2 & 1 & 1 & 3 & 10 & 40 & 171 & 791 \\

3 & 1 & 2 & 10 & 60 & 425 & 3296 & 27447 \\
4 & 1 & 3 & 21 & 182 & 1855 & 20811 & 250439 \\
5 & 1 & 4 & 36 & 408 & 5430 & 79704 & 1254582 \\
6 & 1 & 5 & 55 & 770 & 12650 & 229427 & 4461611 \\
7 & 1 & 6 & 78 & 1300 & 25415 & 548808 & 12706421 \\
8 & 1 & 7 & 105 & 2030 & 46025 & 1152963 & 30966971\\
\hline
\end{tabular}
\caption{The integers $n_d^{(2-2h)}$ for $1\leq d\leq 7$ and $2\leq h \leq 8$.}
\end{figure}

\end{example}

\begin{remark}
\label{rmk:hrv-conjecture}
It is a conjecture of Hausel and Rodriguez Villegas \cite[Remark 4.4.6]{HRV} that for the one-node quiver with $h$ arrows,
we have that
$A^{(h)}_d(1)$ is the dimension of the middle cohomology of the twisted $GL_d$-character variety $\cM_h$ of a genus-$h$ surface.
{
Given that $A_d(1) = \sum_{|\lambda|=d}A_\lambda(1),$ it would be interesting to find a relationship
between other refined Kac polynomials and invariants
of topological strings \cite{LMV}.
Curiously, such results for various genera $h$ would correspond to different
framings of the same genus-\emph{one} Legendrian surface.}
\end{remark}

\appendix

%% file: appendix.tex
\section{Non-compact quantum dilogarithms}
\label{sec:qdl-appendix}
In this appendix, we recall some important properties of the non-compact quantum dilogarithm that we use in the paper. For further background and details regarding this function, we refer the reader to~\cite{FKV,Ka,V}. We assume that $\hbar\in\mathbb{C}$ is such that $\hbar+\hbar^{-1} \in\mathbb{R}$, and lies in the first quadrant $\Re(\hbar)>0,\Im(\hbar)\geq0$. Let us also write
$$
c_\hbar= \frac{i(\hbar + \hbar^{-1})}{2} 
$$
as well as
$$
\zeta = e^{\pi i(1-4c_\hbar^2)/12} \qquad\text{and}\qquad \zeta_{\mathrm{inv}} = \zeta^{-2} e^{-\pi i c_\hbar^2}.
$$

\subsection{The non-compact quantum dilogarithm}

\begin{definition}
Let $C$ be the contour going along the real line from $-\infty$ to $+\infty$, surpassing the origin in a small semi-circle from above. The \emph{non-compact quantum dilogarithm function} $\varphi_\hbar(z)$ is defined in the strip $|\Im(z)| < c_\hbar$ by the following formula \cite{Ka}:
$$
\varphi_\hbar(z) = \mathrm{exp}\left(\frac{1}{4} \int_C\frac{e^{-2izt}}{\mathrm{sinh}(t\hbar)\mathrm{sinh}(t\hbar^{-1})}\frac{dt}{t}\right).
$$
\end{definition}

The non-compact quantum dilogarithm can be analytically continued to the entire complex plane as a meromorphic function with an essential singularity at infinity. The resulting function $\varphi_\hbar(z)$ enjoys the following properties \cite{Ka}:

\begin{itemize}

\item[]{\bf Relation with the compact quantum dilogarithm:}
For $\Im(\hbar^2)>0$, setting $\widetilde{q}  =e^{-\pi i \hbar^{-2}} $ we have
$$
\varphi_\hbar(z) = \frac{(e^{2\pi b(z+c_\hbar)};~q^2)_\infty}{(e^{2\pi b^{-1}(z-c_\hbar)};~\widetilde{q}^2)_\infty}.
$$

\item[]{\bf Poles and zeros:}
$$
\varphi_\hbar(z)^{\pm1} = 0 \quad\Leftrightarrow\quad z = \mp \left( c_\hbar + im\hbar + in\hbar^{-1} \right) \quad\text{for}\quad m,n \in \mathbb{Z}_{\geq 0};
$$

\item[]{\bf Behavior around poles and zeros:}
$$
\varphi_\hbar(z\pm c_\hbar) \sim \pm \zeta^{-1} (2 \pi i z)^{\mp1} \qquad\text{as}\qquad z \to 0;
$$

\item[]{\bf Asymptotic behavior:}
\label{eq:asymp}
$$
\varphi_\hbar(z)\big |_{z\rightarrow\infty} \sim
\begin{cases}
\zeta_{\mathrm{inv}} e^{\pi i z^2}, & |\arg(z)|<\frac{\pi}{2}-\arg(\hbar), \\
1, & |\arg(z)|>\frac{\pi}{2}+\arg(\hbar);
\end{cases}
$$
while we have the following asymptotic behaviour as $\hbar\rightarrow 0$:
$$
\varphi_\hbar\left(\frac{z}{2\pi\hbar}\right)\big |_{\hbar\rightarrow0} \sim \exp\left(\frac{\mathrm{Li}_2(-e^z)}{2\pi i \hbar^2}\right).
$$
\item[]{\bf Symmetry:}
$$
\varphi_\hbar(z) = \varphi_{-\hbar}(z) = \varphi_{\hbar^{-1}}(z);
$$

\item[]{\bf Inversion formula:}
\begin{equation}
\label{inv}
\varphi_\hbar(z) \varphi_\hbar(-z) = \zeta_{\mathrm{inv}} e^{\pi i z^2};
\end{equation}

\item[]{\bf Functional equations:}
\begin{equation}
\label{eqn-func}
\varphi_b\left(z - i\hbar^{\pm1}/2\right) = \left(1 + e^{2\pi \hbar^{\pm1}z}\right) \varphi_b\left(z + i\hbar^{\pm1}/2\right);
\end{equation}

\item[]{\bf Unitarity:}
$$
\overline{\varphi_\hbar(z)} \varphi_\hbar(\overline{z}) = 1;
$$

\end{itemize}

In what follows we will drop the subscript $\hbar$ from the notation for the quantum dilogarithm, and simply write $\varphi(z)$.

\subsection{Integral identites for \texorpdfstring{$\varphi(z)$}{varphi(z)}.}
The quantum dilogarithm function $\varphi(z)$ satisfies many important integral identities. Before describing some of them, let us fix a useful convention.
\begin{remark}
\label{contour-convention}
We will often consider contour integrals of the form 
$$
\int_{C} \prod_{j,k}\frac{\varphi(t-a_j)}{\varphi(t-b_k)}f(t)dt,
$$
where $f(t)$ is some entire function. Unless otherwise specified, the contour $C$ in such an integral is always chosen to be passing below the poles of $\varphi(t-a_j)$ for all $j$, above the poles of $\varphi(t-b_k)^{-1}$ for all $k$, and escaping to infinity in such a way that the integrand is rapidly decaying. 
\end{remark}

The Fourier transform of the quantum dilogarithm can be calculated explicitly by the following integrals:
\begin{align}
\label{eq:qdl-fourier}
\zeta \varphi(w) &= \int\frac{e^{2\pi i x(w-c_b)}}{\varphi(x-c_b)} dx, \\
\frac{1}{\zeta \varphi(w)} &= \int \frac{\varphi(x+c_b)}{e^{2\pi i x(w+c_b)}} dx.
\end{align}

It was shown in~\cite{FKV} that $\varphi$ satisfies the following integral analogs of Ramanujan's $_1\psi_1$ summation formula:
\begin{align}
\label{beta-1}
\frac{\varphi(a) \varphi(w)}{\varphi(a+w-c_b)} &= \zeta^{-1} \int \frac{\varphi(x+a)}{\varphi(x-c_b)} e^{2\pi i x(w-c_b)} dx, \\
\label{beta-2}
\frac{\varphi(a+w+c_b)}{\varphi(a) \varphi(w)} &= \zeta \int \frac{\varphi(x+c_b)}{\varphi(x+a)} e^{-2\pi i x(w+c_b)} dx.
\end{align}
Each of these integral evaluations is equivalent to the non-commutative pentagon identity for $\varphi$ -- for further details, see~\cite{FKV}.